\documentclass[11pt]{article}

\usepackage[a4paper,margin=1in]{geometry}
\usepackage{amsmath,amssymb,amsthm,mathtools}
\usepackage{bm}
\usepackage{enumitem}
\usepackage{xcolor}
\usepackage{needspace}
\usepackage{booktabs,tabularx}
\usepackage{hyperref}

\hypersetup{
  colorlinks=true,
  linkcolor=blue!55!black,
  citecolor=blue!55!black,
  urlcolor=blue!55!black,
  pdftitle={Hydrodynamic and Stochastic Limits of Space-Time Hawkes Processes}
}

\numberwithin{equation}{section}
\setlist[itemize]{leftmargin=2em,itemsep=0.3em,topsep=0.4em}
\setlist[enumerate]{leftmargin=2.2em,itemsep=0.35em,topsep=0.4em}

\newtheorem{theorem}{Theorem}[section]
\newtheorem{proposition}[theorem]{Proposition}
\newtheorem{lemma}[theorem]{Lemma}
\newtheorem{corollary}[theorem]{Corollary}
\theoremstyle{definition}

\theoremstyle{remark}
\newtheorem{remark}[theorem]{Remark}

\newcommand{\R}{\mathbb R}
\newcommand{\Rd}{\mathbb R^d}
\newcommand{\cA}{\mathcal A}
\newcommand{\E}{\mathbb E}
\newcommand{\dd}{\mathrm d}
\newcommand{\cM}{\mathcal M}
\newcommand{\cL}{\mathcal L}
\newcommand{\cN}{\mathcal N}
\newcommand{\cY}{\mathcal Y}
\newcommand{\ip}[2]{\left\langle #1,#2\right\rangle}
\newcommand{\norm}[1]{\left\lVert #1\right\rVert}

\title{Hydrodynamic and Stochastic Limits
of Space-Time Hawkes Processes}
\author{
Fuqing Gao\thanks{School of Mathematics and Statistics, Wuhan University, Wuhan, China; \texttt{fqgao@whu.edu.cn}}
\and
Yingli Wang\thanks{School of Mathematical Sciences, Fudan University, Shanghai, People's Republic of China; \texttt{yingliwang@fudan.edu.cn}}
\and
Lingjiong Zhu\thanks{Department of Mathematics, Florida State University, Tallahassee, Florida, United States of America; \texttt{zhu@math.fsu.edu}}}
\date{\today}

\begin{document}

\maketitle

\begin{abstract}
We study the large-scale behavior of critical and nearly critical
space-time Hawkes processes initialized by a prescribed Poisson field
of past events.  The density of this initial field separates three
macroscopic regimes.  In the high-density regime, fluctuations average
out and the event measure converges to a hydrodynamic limit which 
is a deterministic diffusion profile.  At the fluctuation threshold, random branching remains
visible and produces a measure-valued stochastic limit.  Below the threshold, the
rescaled event measure vanishes locally even though its mean generally
remains nontrivial, with rare large clusters accounting for the surviving mean.
For reproduction delays with finite mean, the critical limits are
described by heat flow and the occupation measure of a branching
diffusion.  A nearly critical perturbation introduces linear decay
without changing the three-regime structure.  Heavy-tailed delays lead
to a different balance between time and space: the mean follows a
fractional diffusion law, the stochastic limit is a superprocess occupation measure on
time-space, characterized by a nonlinear fractional equation, and the past persists as a temporal
memory term.  We also determine how the scaling of the Poisson past must be matched to the observation
scales and show that, in one spatial dimension, a
shorter observation window recovers stochastic behavior hidden in the
sparse regime.  The results allow bounded initial profiles with
possibly infinite total mass and describe the limits directly at the
level of locally finite event measures.
\end{abstract}




\section{Introduction}

A Hawkes process is a simple point process whose intensity depends on its
history; the resulting dynamics are generally non-Markovian.
Hawkes processes have self-exciting and clustering properties
that are popular in applications.
Applications of Hawkes processes include seismology
\cite{ogata1988earthquake,ogata1998spacetime}, finance
\cite{bacry2015finance}, insurance \cite{Stabile,ZhuRuin}, queueing theory \cite{GZ3},
neuroscience \cite{PerniceI}, genomics \cite{Reynaud}, criminology \cite{Mohler}, sociology
\cite{Crane}
and machine learning \cite{Mei2017}; see
\cite{Zhu2013Thesis,HawkesBook2021} for surveys
and the references therein and \cite{reinhart2018review} for a survey on applications
of spatial self-exciting systems.

Self-exciting systems exhibit two competing macroscopic effects.  Repeated
offspring production spreads activity through space and time, while the
random genealogy creates fluctuations that may persist even when many
initial events are present.  The balance becomes especially delicate at
criticality, where every event has one direct child on average.  This paper
studies that balance for critical and nearly critical linear Hawkes
processes on continuous space--time, initialized by a Poisson field
before time zero.  We determine the observation scales and
dimension-dependent intensities of the prescribed past for finite-mean and heavy-tail
reproduction delays.  The main object is the rescaled space--time event
measure.  Its limits range from deterministic diffusion profiles to
random measures and local convergence to zero with a nonzero limiting
first moment.  The temporal tail determines whether the past becomes
an initial condition or persists as a memory term.

\subsection{The linear cluster model and related limits}

A linear temporal Hawkes process is a random counting measure $N(\dd t)$ that has the predictable intensity
\begin{equation}\label{eq:intro-temporal-Hawkes}
  \lambda_t:=\mu+\int_{(-\infty,t)}h(t-s)N(\dd s),
  \qquad \mu>0,
\end{equation}
where $h:\mathbb{R}_{+}\rightarrow\mathbb{R}_{\geq0}$ is the reproduction \textit{kernel function} (or exciting function)
\cite{hawkes1971spectra} and $\mu>0$ is the \textit{baseline intensity}.  The linear Hawkes process has an immigration--birth representation
\cite{hawkes1974cluster}; immigrants arrive at rate $\mu$ and each
event independently produces a Poisson process of children with
intensity $h(s)\dd s$. The total number of immigrants and the descents on an time interval $(0,t]$ gives $N(0,t]$. The branching ratio $\int_0^\infty h(s)\dd s$
is the mean number of direct children.  At branching ratio one, each
cluster is finite almost surely but has infinite expected total size;
with positive stationary immigration the mean intensity cannot be finite.
This paper instead prescribes a Poisson past and follows its descendants
after time zero, with no subsequent immigration.

A space--time Hawkes process replaces event times by event locations
$(t,x)$ in $\mathbb R_+\times\mathbb R^d$.  It is a random counting measure
$N(\dd t,\dd x)$ with predictable intensity density $\lambda(t,x)$: for
each bounded Borel set $A\subset\mathbb R^d$,
\begin{equation}\label{eq:intro-spacetime-compensator}
  N((0,t]\times A)
  -\int_0^t\int_A\lambda(s,x)\dd x\dd s
\end{equation}
is a local martingale. In the translation-invariant linear model,
\begin{equation}\label{eq:intro-Hawkes}
  \lambda(t,x)
  :=\mu(t,x)+\int_{(-\infty,t)\times\mathbb R^d}
       h(t-s,x-y)N(\dd s,\dd y).
\end{equation}
Here $\mu:\mathbb{R}\times\mathbb{R}^{d}\rightarrow\mathbb{R}_{\geq 0}$ is the baseline intensity and $h:\mathbb{R}\times\mathbb{R}^{d}\rightarrow\mathbb{R}_{\geq 0}$ is the kernel function (or exciting function) such that $h(r,\cdot)=0$ for $r\leq0$; its first coordinate describes the
reproduction delay and its second the spatial displacement.  The
factorized case $h(r,z)=\phi(r)p(z)$ separates the temporal response from
the displacement density, but our paper allows dependence between the two
coordinates.  Unlike a multivariate Hawkes process with finitely many event
types \cite{hawkes1971spectra,BacryEtAl2013}, \eqref{eq:intro-Hawkes} uses a continuum of physical locations and
therefore permits an actual spatial rescaling.

The cluster construction has an equally concrete space--time form.  A
parent at $(s,y)$ produces a Poisson random measure of children at
$(s+r,y+z)$ with intensity $h(r,z)\dd r\dd z$.  Consequently,
\begin{equation}\label{eq:intro-spacetime-branching-ratio}
  \norm{h}_1
  :=\int_0^\infty\int_{\mathbb R^d}h(r,z)\dd z\dd r
\end{equation}
is both the mean number of direct children and the space--time branching
ratio.  Conditional on a child being born, $(S,Z)$ has joint density
$h/\norm{h}_1$.  At the critical value $\norm{h}_1=1$, a cluster is at once
a critical Crump--Mode--Jagers age-dependent branching process and a
branching random walk.  The age $S$ records reproduction time and $Z$
records displacement in $\Rd$.  When $m_1=\E[S]<\infty$, this mean age
converts genealogical generations into physical time, whereas the centered
displacement covariance $\Sigma=\E[ZZ^{\top}]$ produces spatial Brownian motion.  The finite-mean theory requires no independence between
$S$ and $Z$; the heavy-tail theorem uses a factorized kernel.

For (linear) temporal Hawkes processes,
the asymptotics have developed along several directions.
The large-time functional law of large numbers
and central limit theorems are obtained by \cite{BacryEtAl2013}.
The large-time large deviations are treated by
\cite{bordenave2007large}.
The large-time precise deviations are studied in \cite{GaoZhu2021Precise}.
Large initial intensity leads to laws of large
numbers and Gaussian fluctuations for temporal Markovian Hawkes processes
\cite{gao2018largeinitial} as well as the large deviations \cite{GaoZhu2018InitialLDP}.
Nearly unstable temporal processes converge
to integrated Cox-Ingersoll-Ross (CIR) type limits, with fractional analogues for heavy-tailed
kernels \cite{jaisson2015limit,jaisson2016rough}.
Limit theorems for (linear) temporal Hawkes processes with random marks
are first studied
for the univariate case in \cite{KarabashZhu2015}.
For the multivariate case,
functional law of large numbers and central limit theorem are obtained in
\cite{HorstXu2021}, and large deviations are established in \cite{Karim2025}.
When the mutual excitation rates are heavy-tailed,
sample path large deviations are obtained in \cite{Blanchet2025}.

Recent work has refined the critical and nearly critical temporal theory.
Functional limit theorems with CIR and integrated-CIR behavior, Feller
random-measure limits, and stochastic Volterra limits are obtained in
\cite{horst2026functional}, \cite{pace2026feller}, and
\cite{xu2024volterra}, respectively.  Connections with general
age-dependent branching include the excursion limits of
\cite{lambert2013excursion}, the diffusion approximations of
\cite{xu2024diffusion}, and the many-ancestor laws of large numbers and
central limit theorems of \cite{drame2025cmj}.  In particular, the cluster
Laplace approach and random-measure formulation used here have precedents
in the temporal Hawkes literature.  Our additional spatial scaling leads
to a heat operator in the limiting log-Laplace equation, while the
prescribed Poisson past determines its initial measure.

There are also close precedents outside the Hawkes literature.
\cite{meleard2012age} prove superprocess limits for an interacting
birth--death population with trait and age structure, using separation
of fast age and slow trait scales.  \cite{fleischmann2005clustering}
study multiscale clustering in a non-Markovian spatial branching system
with finite-mean lifetimes, including increasing Poisson initial
populations at the critical dimension.  These are precedents for age averaging, superprocess limits, and
clustering at different initial population densities.  Here the Hawkes
mechanism produces offspring at multiple ages, allows dependent delay
and displacement, and counts births throughout space--time.  The
proofs treat this observable and its prescribed Poisson past directly.

Despite the numerous works on temporal Hawkes processes,
the theory for spatial or space-time Hawkes processes is less
developed.
Spatial Hawkes asymptotics include both network limits and long-time
limits on a fixed infinite graph.  \cite{delattre2016networks} construct
Hawkes processes on possibly infinite networks and study both mean-field
limits and long-time behavior on $\mathbb Z^d$.  The spatial mean-field
models of \cite{chevallier2019spatial} lead to neural-field equations,
and \cite{chevallier2020fluctuations} study Gaussian fluctuations about
such a limit.  \cite{maclaurin2026adaptive} obtain a hydrodynamic limit
as the size of an adaptive network grows.  More recently,
\cite{belmabrouk2026longrange} study long-time laws of large numbers for
an infinite lattice with power-law interactions in the subcritical and
supercritical branching regimes.

In this paper, we are
interested in studying
the limiting behavior
of space-time Hawkes processes
by jointly rescaling
physical time and space in a continuum event measure, first at branching
ratio one and then under near-critical perturbations.  The tail of the
reproduction delay determines the temporal balance, while the past
intensity provides a separate parameter governing the fluctuation threshold.

\subsection{Scaling mechanism and contributions}
\label{subsec:intro-scaling}

Unless explicitly stated otherwise, every limit involving an
$\epsilon$-indexed quantity is taken as $\epsilon\downarrow0$.
Let $E:=(0,\infty)\times\Rd$, and let $\cM_{\mathrm{loc}}(E)$ denote
the space of nonnegative Radon measures on $E$: each measure has finite
mass on every compact subset of $E$, but its total mass may be infinite.
Let $C_c(E)$ be the space of continuous real-valued functions with compact
support in $E$, and let $C_c^+(E)$ be its nonnegative members.  For
$\mu\in\cM_{\mathrm{loc}}(E)$ and $\varphi\in C_c(E)$, write
$\ip{\mu}{\varphi}:=\int_E\varphi(t,x)\mu(\dd t,\dd x)$.
We equip $\cM_{\mathrm{loc}}(E)$ with the vague topology, characterized by
\[
 \mu_n\xrightarrow{v}\mu
 \quad\text{i.e.}\quad
 \ip{\mu_n}{\varphi}\longrightarrow\ip{\mu}{\varphi}
 \quad\text{for every }\varphi\in C_c(E).
\]
For random measures, convergence in law or in probability is understood
with respect to this topology.
For a random measure $\xi$, let
$L_\xi(f):=\E[e^{-\ip{\xi}{f}}]$, $f\in C_c^+(E)$, be its Laplace
functional.  The corresponding convergence criterion is
\[
 \xi_n\Rightarrow\xi\quad\text{in }\cM_{\mathrm{loc}}(E)
 \quad\text{i.e.}\quad
 L_{\xi_n}(f)\longrightarrow L_\xi(f)
 \quad\text{for every }f\in C_c^+(E).
\]
Lemma~\ref{lem:Laplace-vague} proves this criterion and the fact that
$L_\xi$ uniquely determines the law of $\xi$.  Our stochastic-limit
proofs compute the Laplace functional from the Poisson-cluster
representation and identify its limit with that of a locally finite
random measure.  If the limiting measure is deterministic, the same
criterion also gives convergence in probability in the vague topology.

\paragraph*{Finite-mean delays}
In our paper, a hydrodynamic limit means convergence of the normalized
space--time event measure to a deterministic density.  The high-density
proof uses tightness, a martingale identity, and uniqueness of the
limiting weak equation.  For this branching model, the drift closes
through the offspring operator
\begin{equation}\label{eq:intro-forward-operator}
  K_\epsilon\varphi(t,x)
    :=\int_0^\infty\int_{\Rd}
        h(s,z)\varphi\left(t+\epsilon^2s,x+\epsilon z\right)\dd z\dd s.
\end{equation}
Under a finite mean delay $m_1$, centered displacement, and finite
spatial covariance $\Sigma$,
\begin{equation}\label{eq:intro-generator-expansion}
  \epsilon^{-2}(K_\epsilon\varphi-\varphi)
    \xrightarrow[\epsilon\downarrow0]{} m_1\partial_t\varphi+\frac12\Sigma:D_x^2\varphi
    \quad\text{in }L^1.
\end{equation}
Here $D_x^2\varphi=(\partial_{x_i}\partial_{x_j}\varphi)_{i,j=1}^d$
is the spatial Hessian, and the colon denotes the entrywise matrix
inner product:
\[
 \Sigma:D_x^2\varphi
 :=\sum_{i=1}^d\sum_{j=1}^d\Sigma_{ij}
                              \partial_{x_i}\partial_{x_j}\varphi.
\]
The limit in \eqref{eq:intro-generator-expansion} identifies the
diffusion matrix $D=\Sigma/(2m_1)$.
The past intensity determines whether the branching fluctuations
survive on the same scale.

Two notions of criticality must be separated.  In Section~\ref{sec:critical}, $\norm{h}_1=1$ fixes the microscopic genealogy at critical
branching.  A density parameter $\alpha>0$ controls the density
$\epsilon^{-(\alpha+1-d)}\lambda_0(\epsilon x)$ of the prescribed Poisson
past in $\Rd$.  Its microscopic density exponent is
$\alpha+1-d$, which can be negative; the comparison with the fluctuation
threshold, rather than divergence of microscopic density, defines the
three regimes below.  After division by the macroscopic equation scale, the
martingale bracket is of order $\epsilon^{\alpha-1}$.  This creates a
phase transition:
\begin{itemize}
  \item for $\alpha>1$, the noise vanishes and the hydrodynamic limit is the
  deterministic heat flow;
  \item at $\alpha=1$, the noise remains of order one and the limit is the
  occupation measure of a Dawson--Watanabe superprocess;
  \item for $0<\alpha<1$, the original long-scale measure vanishes locally
  in probability even though its expectation follows the heat flow;
  for a nontrivial initial profile this is a loss of uniform integrability.
\end{itemize}
Thus the heat equation is only the first-moment description below the
threshold, while at the threshold it is the drift part of a genuinely
stochastic measure-valued evolution.
The critical stochastic limit belongs to the Dawson--Watanabe universality
class.  Superprocesses arise as high-density limits of critical branching
particle systems, beginning with \cite{watanabe1968limit}; see also
\cite{dynkin1991branching,etheridge2000superprocesses,perkins2002superprocess}.
The canonical measure of critical branching limits is discussed in
\cite{perkins2004critical}; it provides the Laplace exponent in our sparse
expansion.  The finite-mean results retain the age-dependent clock,
the prescribed past, and the space--time event observable: the established superprocess
$X$ is seen through $m_1^{-1}X_t(\dd x)\dd t$, with initial measure
$X_0(\dd x)=m_1\lambda_0(x)\dd x$.

\paragraph*{Heavy-tailed delays}
For a reproduction delay with tail index $\eta\in(0,1)$, the renewal
scale becomes $\delta_\epsilon=\epsilon^{2\eta}$.  Over a physical time
window of order $\epsilon^{-2}$, an ancestral line contains order
$\epsilon^{-2\eta}$ generations, so centered finite-variance spatial
displacements accumulate on the scale $\epsilon^{-\eta}$.  This leads
to the space--time rescaling
$(t,x)\mapsto(\epsilon^2t,\epsilon^\eta x)$ and to the near-critical
balance $1-a_\epsilon$ of order $\epsilon^{2\eta}$.  The normalized
martingale bracket is of order $\epsilon^{\eta(\alpha-1)}$, so the
fluctuation threshold remains $\alpha=1$.  On this scale the mean is
governed by a Caputo fractional diffusion equation rather than the heat
equation.  At the operator level, the first-order temporal term
$m_1\partial_t$ is replaced by the one-sided stable generator
$-c_\eta D_+^\eta$, while the spatial part remains diffusive.  Moreover,
the rescaled past remains distributed over positive physical times
instead of collapsing to an initial measure at time zero.

\paragraph*{Contributions}

Our main results are as follows; see also Table~\ref{tab:results-overview}.
\begin{enumerate}
  \item For factorized kernels with temporal tail index
  $\eta\in(0,1)$ and finite spatial covariance,
  Theorem~\ref{thm:ht-limits} establishes the mean, deterministic,
  stochastic, and sparse event-measure limits when
  $(1-a_\epsilon)/\epsilon^{2\eta}\to\kappa\geq0$.
  The past survives as $c_Ht^{-\eta}\lambda_0(x)\dd t\dd x$.
  Proposition~\ref{prop:ht-lifted} realizes the stochastic limit as
  the occupation measure of a quadratic superprocess on space--time,
  with stable temporal motion and Brownian spatial motion.
  Its mean solves the Caputo equation \eqref{eq:ht-caputo}.
  This representation also shows that applying one common inverse
  stable clock to an ordinary superprocess reproduces the mean behavior
  but not the branching law of the limiting random measure.

  \item In the finite-mean setting,
  Proposition~\ref{prop:annealed-heat} and
  Theorems~\ref{thm:main}, \ref{thm:critical}, and
  \ref{thm:sparse-long} establish the classical heat-flow,
  Dawson--Watanabe, and canonical-measure scaling picture for the
  Hawkes event observable and its prescribed Poisson past.
  All three regimes require only a first temporal moment and a
  centered second spatial moment, and allow dependence between
  delay and displacement and bounded profiles of infinite total mass.
  Theorem~\ref{thm:critical-varying} gives stability under locally
  uniform changes of the initial profile.
  Theorem~\ref{thm:nearly-critical} includes the perturbation
  $(1-a_\epsilon)/\epsilon^2\to\kappa\geq0$, with limiting decay
  rate $\kappa/m_1$ and initial measure $m_1\lambda_0(x)\dd x$.
  In dimension one, the stochastic limits admit densities satisfying
  stochastic heat equations, with the near-critical perturbation
  appearing as linear decay.

  \item Proposition~\ref{prop:gamma-phase} matches the scaling of the Poisson past
  to the observation scales.  In dimension one,
  Theorem~\ref{thm:local-critical} and Corollary~\ref{cor:nc-local}
  recover a critical occupation limit on a shorter window below the
  original fluctuation threshold.  Matching the past intensity to the observation scale retains the
  full profile, while local observation freezes it at the observation
  point.

  \item Lemma~\ref{lem:renewal-framework} gives a common
  weighted-resolvent framework for the finite-mean and heavy-tail
  nonlinear cluster equations.  Corollary~\ref{cor:clock-unification}
  places the two stochastic limits within a single space--time
  occupation construction and identifies the corresponding limiting
  source measures generated by the prescribed past.
\end{enumerate}

\begin{table}[p]
\centering
\caption{Overview of the scaling limits.  All limits concern locally
finite event measures with the vague topology.  The branching ratio
and the density parameter $\alpha$ play separate roles.}
\label{tab:results-overview}
\small
\setlength{\tabcolsep}{5pt}
\renewcommand{\arraystretch}{1.18}

\begin{tabularx}{\textwidth}{@{}>{\raggedright\arraybackslash}p{.21\textwidth}
 *{3}{>{\raggedright\arraybackslash}X}@{}}
\toprule
\multicolumn{4}{@{}l}{\textbf{A. Kernel, observation scale, and normalization}}\\
\midrule
Quantity & Critical & Nearly critical:\newline finite mean
         & Nearly critical:\newline heavy tail\\
\midrule
Branching ratio
 & $a_\epsilon=1$
 & $\dfrac{1-a_\epsilon}{\epsilon^2}\to\kappa\geq0$
 & $\dfrac{1-a_\epsilon}{\epsilon^{2\eta}}\to\kappa\geq0$\\
Observed $(t,x)$
 & $(\epsilon^2t,\epsilon x)$
 & $(\epsilon^2t,\epsilon x)$
 & $(\epsilon^2t,\epsilon^\eta x)$\\
Past prefactor $q_\epsilon$
 & $\epsilon^{-(\alpha+1-d)}$
 & $\epsilon^{-(\alpha+1-d)}$
 & $\epsilon^{2-\eta(\alpha+3-d)}$\\
Mass per event
 & $\epsilon^{\alpha+3}$
 & $\epsilon^{\alpha+3}$
 & $\epsilon^{\eta(\alpha+3)}$\\
Kernel assumptions
 & \eqref{eq:H1}--\eqref{eq:H4}
 & \eqref{eq:H1}--\eqref{eq:H4}
 & \eqref{eq:ht-kernel}--\eqref{eq:ht-spatial-moments};\newline $0<\eta<1$\\
\bottomrule
\end{tabularx}

\medskip
\begin{tabularx}{\textwidth}{@{}>{\raggedright\arraybackslash}p{.21\textwidth}
 *{3}{>{\raggedright\arraybackslash}X}@{}}
\toprule
\multicolumn{4}{@{}l}{\textbf{B. Limits on the scales in A}}\\
\midrule
Density regime & Critical & Finite mean & Heavy tail\\
\midrule
$\alpha>1$\newline (in probability)
 & $\rho(t,x)\dd t\dd x$
 & $\rho^\kappa(t,x)\dd t\dd x$
 & $\rho^{\eta,\kappa}(t,x)\dd t\dd x$\\
$\alpha=1$\newline (in law)
 & $m_1^{-1}X_t(\dd x)\dd t$
 & $m_1^{-1}X_t^\kappa(\dd x)\dd t$
 & $\nu_{\lambda_0}^{\eta,\kappa}$\\
$0<\alpha<1$\newline (zero in probability)
 & Exponent $\Psi_{\lambda_0}$;\newline
   $r_\epsilon=\epsilon^{1-\alpha}$
 & Exponent $\Psi_{\lambda_0}^\kappa$;\newline
   $r_\epsilon=\epsilon^{1-\alpha}$
 & Exponent $\Psi_{\lambda_0}^{\eta,\kappa}$;\newline
   $r_\epsilon=\epsilon^{\eta(1-\alpha)}$\\
Results
 & Proposition~\ref{prop:annealed-heat};\newline
   Theorems~\ref{thm:main}, \ref{thm:critical}, \ref{thm:sparse-long}
 & Theorem~\ref{thm:nearly-critical}
 & Theorem~\ref{thm:ht-limits}\\
\bottomrule
\end{tabularx}

\medskip
\begin{tabularx}{\textwidth}{@{}>{\raggedright\arraybackslash}p{.21\textwidth}
 >{\raggedright\arraybackslash}p{.18\textwidth}
 >{\raggedright\arraybackslash}X
 >{\raggedright\arraybackslash}p{.19\textwidth}@{}}
\toprule
\multicolumn{4}{@{}l}{\textbf{C. Alternative observation scales (finite-mean setting)}}\\
\midrule
Setting & Parameter & Limit & Result\\
\midrule
Matched, critical
 & $\vartheta>0$
 & Heat flow, in probability
 & Proposition~\ref{prop:gamma-phase}\\
Matched, critical
 & $\vartheta=0$
 & Dawson--Watanabe occupation, in law
 & Proposition~\ref{prop:gamma-phase}\\
Matched, critical
 & $\vartheta<0$
 & Zero in probability; heat-flow mean
 & Proposition~\ref{prop:gamma-phase}\\
Local, $d=1$:\newline critical or nearly critical
 & $0<\alpha<1$;\newline $\gamma=\alpha$
 & Critical occupation, in law;\newline
   $X_0=m_1\lambda_0(x_\star)\dd x$
 & Theorem~\ref{thm:local-critical};\newline
   Corollary~\ref{cor:nc-local}\\
\bottomrule
\end{tabularx}

\medskip
\begin{minipage}{\textwidth}
\footnotesize
In A, the past density is $q_\epsilon\lambda_0(\epsilon y)$ for the
first two columns and $q_\epsilon\lambda_0(\epsilon^\eta y)$ for the
third.  The profile $\lambda_0\geq0$ is bounded and uniformly continuous.

\smallskip
In B, $\rho=P_t^D\lambda_0$, $\rho^\kappa=P_t^{D,\kappa}\lambda_0$,
and $\rho^{\eta,\kappa}$ is the fractional mean in
\eqref{eq:ht-mean-profile}.  For every $\alpha>0$, the mean measures
converge to the deterministic measure in the first row.
The exponent notation means
$-r_\epsilon^{-1}\log\E[e^{-\ip{r_\epsilon\cN_\epsilon}{g}}]\to\Psi(g)$,
using the event measure of the corresponding column.

\smallskip
In C, the matched family has past intensity
$\epsilon^{-\beta}\lambda_0(\epsilon^\gamma y)$, observed coordinates
$(\epsilon^{2\gamma}t,\epsilon^\gamma x)$, and mass per event
$\epsilon^{\beta+(d+2)\gamma}$; assume $\gamma>0$ and
$\beta+(d-1)\gamma>0$, and set $\vartheta=\beta+(d-2)\gamma$.
The local row uses the original Poisson past, centers at $x_\star/\epsilon$,
and has event mass $\epsilon^{4\alpha}$.
\end{minipage}
\end{table}

\section{Model, Assumptions, and Scaling}
\label{sec:model}

\subsection{The offspring kernel}

Fix an integer $d\geq1$.  Spatial vectors carry the Euclidean norm
$|\cdot|$, and we write $B_R:=\{x\in\Rd:|x|\leq R\}$.
Let
\[
  h:\R_+\times\Rd\longrightarrow\R_+
\]
be a Borel function.  We impose the following assumptions:
\begin{align}
  \int_0^\infty\int_{\Rd} h(s,z)\dd z\dd s &=1,
  \label{eq:H1}\\
  \int_0^\infty\int_{\Rd} z h(s,z)\dd z\dd s &=0,
  \label{eq:H2}\\
  m_1:=\int_0^\infty\int_{\Rd} s h(s,z)\dd z\dd s
  &\in(0,\infty),
  \label{eq:H3}\\
  \Sigma:=\int_0^\infty\int_{\Rd} zz^{\top}h(s,z)\dd z\dd s
  &\quad\text{ is finite and positive definite}.
  \label{eq:H4}
\end{align}
Condition \eqref{eq:H1} is criticality.  Condition \eqref{eq:H2} centers
the displacement of a typical child.  Since \eqref{eq:H1} and
\eqref{eq:H4} imply integrability of $|z|h(s,z)$, the centered first moment
is well defined.  Equivalently, the finite second-moment condition is
$\int_0^\infty\int_{\Rd} |z|^2h(s,z)\dd z\dd s=\operatorname{tr}(\Sigma)<\infty$.
In dimension one, we also write $m_2=\Sigma_{11}$.  No product structure
between reproduction time and spatial displacement is assumed.
Set
\begin{equation}\label{eq:diffusion-coefficient}
  D:=\frac{\Sigma}{2m_1},
  \qquad \cA\psi:=D:D_x^2\psi.
\end{equation}
Let $(P_t^D)_{t\geq0}$ denote the heat semigroup with generator $\cA$.

\paragraph*{Interpretation and scope of the kernel assumptions}

These are conditions on a fixed microscopic kernel, which can be checked
before taking a scaling limit.  They specify a critical, centered,
finite-moment regime; they are not claimed to cover all Hawkes models.
Their roles are as follows.
\begin{enumerate}
  \item \emph{Criticality, \eqref{eq:H1}.}  Each event has a Poisson number
  of children with mean one, so the genealogy is a critical branching
  tree.  Unit mass is a substantive choice of regime, not a normalization
  that can be imposed without changing the dynamics.  Subcritical and
  supercritical kernels have different long-time behavior; see, for
  example, \cite{delattre2016networks,belmabrouk2026longrange}.
  \item \emph{Centering, \eqref{eq:H2}.}  The mean spatial displacement is
  zero.  Without this condition the forward-operator expansion contains
  the term
  $\epsilon^{-1}(\int_0^\infty\int_{\Rd} z h(s,z)\dd z\dd s)
  \cdot\nabla_x\varphi$ on the diffusive
  equation scale.  Centering removes that leading transport term; it does
  not require symmetry or independence of the two coordinates.
  \item \emph{Finite mean reproduction time, \eqref{eq:H3}.}  Along a
  lineage, the sum of reproduction delays divided by the number of
  generations converges to $m_1$.  Thus genealogical time has a
  deterministic first-order conversion into physical time.  The same
  moment is the total weight of direct offspring crossing time zero,
  $\int_{s<0}\int_{q>0}\int_{\Rd} h(q-s,z)\dd z\dd q\dd s=m_1$.
  None of the finite-mean limits requires a second temporal moment.
  \item \emph{Finite spatial variance, \eqref{eq:H4}.}  This is the usual
  finite-variance setting for Brownian spatial scaling.  Combining the
  spatial covariance per generation with the mean time per generation
  gives the diffusion matrix $D=\Sigma/(2m_1)$.  The $L^1$ proof of this
  operator limit is given in
  Lemma~\ref{lem:generator}, without a third spatial or mixed moment
  assumption.
\end{enumerate}
The same moment assumptions suffice for all three finite-mean regimes.
Lemma~\ref{lem:nonlinear-renewal} uses a resolvent contraction, and
Lemma~\ref{lem:spatial-localization} provides a uniform second-moment
bound for cluster costs.  The contribution of the prescribed past is controlled by its
finite crossing mass $m_1$.  Consequently, exponential temporal and
Gaussian spatial kernels are included, as are dependent delay and
displacement laws satisfying the stated moments.

The assumptions have explicit nontrivial examples.  Let $a$ be a smooth
probability density supported in $(1,2)$ and let $p$ be a smooth,
symmetric, compactly supported probability density on $\Rd$ with positive
definite covariance.  Then $h(s,z)=a(s)p(z)$ satisfies
\eqref{eq:H1}--\eqref{eq:H4}.  Dependence between time and displacement is
also allowed: take
$h(s,z)=\frac12 a_1(s)p_1(z)+\frac12 a_2(s)p_2(z)$,
where $a_1,a_2$ are smooth compactly supported probability densities with
disjoint positive-time supports and $p_1,p_2$ are smooth symmetric compactly
supported densities on $\Rd$ with different positive definite covariance
matrices.  The conditional
spatial covariance then depends on the reproduction delay, while all four
assumptions still hold.

\subsection{Prescribed Poisson past and positive-time dynamics}

Let $\alpha>0$ and let $\lambda_0:\Rd\to\R_+$ be bounded and uniformly
continuous.  For every
$\epsilon\in(0,1]$, let $\Pi^\epsilon$ be a Poisson random measure on
$(-\infty,0)\times\Rd$ with intensity
\begin{equation}\label{eq:past-intensity}
  \epsilon^{-(\alpha+1-d)}\lambda_0(\epsilon y)\dd s\dd y.
\end{equation}
The factor $\lambda_0(\epsilon y)$, rather than $\lambda_0(y)$, is
essential: $y$ is the microscopic coordinate and $\epsilon y$ is the
macroscopic one.  With an unshifted microscopic density exponent $\beta$,
the first-moment normalization would be $\epsilon^{\beta+d+2}$ and
the martingale bracket on the equation scale would be of order
$\epsilon^{\beta+d-2}$.  We choose $\beta=\alpha+1-d$ to keep the
threshold at $\alpha=1$ in every dimension.  The past intensity is thus of
constant microscopic order at criticality when $d=2$, and tends to zero
microscopically when $d>2$.

On $(0,\infty)\times\Rd$, we define the space-time Hawkes process $N^\epsilon$ as a locally finite point process,
with no two events at the same time, and with the predictable intensity
\begin{align}
  \lambda^\epsilon(t,x)
  ={}&\int_{(-\infty,0)\times\Rd}
       h(t-s,x-y)\Pi^\epsilon(\dd s,\dd y)
  \nonumber\\
  &\qquad\qquad+\int_{(0,t)\times\Rd}
       h(t-s,x-y)N^\epsilon(\dd s,\dd y),
  \qquad t>0.
  \label{eq:hawkes-intensity}
\end{align}
The filtration contains the prescribed past history at time zero.  We work in the
class of solutions whose mean intensity is uniformly bounded in space
on each finite time interval:
\begin{equation}\label{eq:Hawkes-mean-class}
  \operatorname*{ess\,sup}_{(t,x)\in(0,T]\times\Rd}
       \E\left[\lambda^\epsilon(t,x)\right]<\infty
  \qquad\text{for every }T<\infty.
\end{equation}
The cluster construction below verifies this condition from boundedness
of $\lambda_0$.  This class gives uniqueness in every dimension without
requiring higher spatial moments to control a translated integrable weight.

\begin{proposition}[Construction and uniqueness with bounded mean intensity]
\label{prop:construction}
Under \eqref{eq:H1}--\eqref{eq:H4}, for every fixed $\epsilon>0$,
equations \eqref{eq:past-intensity}--\eqref{eq:hawkes-intensity} have a
locally finite Poisson-cluster solution.  It is unique in law among
solutions satisfying \eqref{eq:Hawkes-mean-class}, and pathwise unique
when realized with the same prescribed past and Poisson driving measure.
\end{proposition}

\begin{proof}
First restrict the past field to a bounded space--time box.  For each past
parent $(r,y)$ with $r<0$, generate only its direct children born at
positive times, with intensity
\begin{equation}\label{eq:past-direct-offspring}
  \mathbf 1_{\{q>0\}}h(q-r,w-y)\dd q\dd w.
\end{equation}
Each such child then generates its full positive-time descendant cluster
with offspring intensity $h(s,z)\dd s\dd z$.  No additional negative-time
events are generated: the entire negative-time configuration is the given
Poisson field.  The first-generation offspring number has mean
\[
\int_{-r}^{\infty}\int_{\Rd} h(s,z)\dd z\dd s\leq1,
\]
and all subsequent
offspring numbers are Poisson with mean one.  Thus a parent clan is
dominated by a full critical Poisson Galton--Watson tree, which is finite
almost surely.  Coupling the children and their descendants consistently
as the past boxes increase gives a monotone sequence of point measures.

Here is a direct bound that also justifies the construction.  Write
\[
c_\epsilon:=\epsilon^{-(\alpha+1-d)}\norm{\lambda_0}_\infty,
\]
and
\begin{align}\label{H:s:function}
H(s):=\int_{\Rd} h(s,z)\dd z
\quad\text{and}\quad
\overline H(t):=\int_t^\infty H(s)\dd s.
\end{align}
The mean density of direct
positive-time children is bounded by $c_\epsilon\overline H(t)$.
Let $b_\epsilon^0$ be the mean density of direct positive-time children
of the past, and let $\bar\lambda_{\epsilon,n}$ count their descendants
through $n$ generations, with $\bar\lambda_{\epsilon,0}=0$.
Integrating the direct-offspring intensity \eqref{eq:past-direct-offspring}
against \eqref{eq:past-intensity} gives $b_\epsilon^0$.
The two terms in \eqref{eq:hawkes-intensity} then separate children
of past parents from children of positive-time parents, giving
\begin{align}
 b_\epsilon^0(t,x)
 &=\epsilon^{-(\alpha+1-d)}\int_0^\infty\int_{\Rd}
       h(t+a,z)\lambda_0(\epsilon(x-z))\dd z\dd a
       \leq c_\epsilon\overline H(t),\nonumber\\
 \bar\lambda_{\epsilon,n+1}(t,x)
 &=b_\epsilon^0(t,x)+\int_0^t\int_{\Rd}
       h(t-s,x-y)\bar\lambda_{\epsilon,n}(s,y)\dd y\dd s
       \nonumber\\
 &\leq c_\epsilon\overline H(t)
       +c_\epsilon\int_0^t H(s)\dd s=c_\epsilon,
       \qquad n\geq0,
 \label{eq:construction-generation-bound}
\end{align}
where the last inequality follows inductively from
$\bar\lambda_{\epsilon,n}\leq c_\epsilon$.
Taking increasing limits in generations and past boxes gives the
measure bound
\[
  \E\left[N^\epsilon(\dd t,\dd x)\right]\leq c_\epsilon\dd t\dd x.
\]
Consequently the limit is locally finite, satisfies
\eqref{eq:Hawkes-mean-class}, and, for every bounded Borel set $B\subset\Rd$,
\begin{equation}\label{eq:Hawkes-local-construction-bound}
  \E\left[ N^\epsilon((0,T]\times B)\right]\leq c_\epsilon T|B|.
\end{equation}
The independent offspring processes give
\eqref{eq:hawkes-intensity}; their absolutely continuous time laws also
give the asserted absence of simultaneous events.

For uniqueness, use a Poisson random measure $\mathcal P$ of intensity
$\dd t\dd x\dd a$, independent of $\Pi^\epsilon$, and the thinning equation
\begin{equation}\label{eq:Hawkes-Poisson-equation}
  N^\epsilon(\dd t,\dd x)
  =\int_0^\infty\mathbf 1_{\{a\leq\lambda^\epsilon(t,x)\}}
       \mathcal P(\dd t,\dd x,\dd a).
\end{equation}
This representation also applies to any solution in
\eqref{eq:Hawkes-mean-class} on an enlarged probability space.  To see
this, enumerate its events as $(T_i,X_i)_{i\in I}$, $I\subset\mathbb N$,
using a countable partition into bounded boxes.  Independently attach uniform marks
$U_i\in(0,1)$, revealed at $T_i$, and let
\[
 \mathcal P^-(\dd t,\dd x,\dd a)
 :=\sum_{i\in I}\delta_{(T_i,X_i,U_i\lambda^\epsilon(T_i,X_i))}
                         (\dd t,\dd x,\dd a).
\]
The compensator formula for independently marked events gives, for
nonnegative predictable $F$,
\begin{align*}
 \E\left[\int_{E\times\R_+}F(t,x,a)\,\dd\mathcal P^-\right]
 &=\E\left[\int_E\lambda^\epsilon(t,x)
       \int_0^1F(t,x,u\lambda^\epsilon(t,x))\dd u\dd t\dd x\right]\\
 &=\E\left[\int_E\int_0^{\lambda^\epsilon(t,x)}
                                  F(t,x,a)\dd a\dd t\dd x\right].
\end{align*}
The compensator identity gives zero event mass on
$\{\lambda^\epsilon=0\}$.  By \eqref{eq:Hawkes-mean-class},
$\{\lambda^\epsilon=\infty\}$ is null for
$\mathbb P\otimes\dd t\dd x$ and also carries zero event mass.
Thus all event marks lie in $\{0<\lambda^\epsilon<\infty\}$ almost surely.  Take an independent Poisson measure $\mathcal P^0$ of
intensity $\dd t\dd x\dd a$ and define
\[
 \mathcal P:=\mathcal P^-+
       \mathbf1_{\{a>\lambda^\epsilon(t,x)\}}\mathcal P^0.
\]
Its compensator is
$\left(\mathbf1_{\{a\leq\lambda^\epsilon\}}+
\mathbf1_{\{a>\lambda^\epsilon\}}\right)\dd t\dd x\dd a=\dd t\dd x\dd a$.
On bounded boxes it has no simultaneous times.  For deterministic
bounded $F\geq0$ supported in such a box with time horizon $T$, the
point-process exponential martingale is
\[
 Z_t=\exp\left\{-\int_{(0,t]\times\Rd\times\R_+}F(s,x,a)\,\dd\mathcal P
       +\int_0^t\int_{\Rd}\int_0^\infty\left(1-e^{-F(s,x,a)}\right)
                                          \dd a\dd x\dd s\right\}.
\]
It is bounded on $[0,T]$, so $\E[Z_T\mid\mathcal F_0]=1$.  Thus,
\begin{equation}\label{eq:Poisson-embedding-Laplace}
 \E\left[e^{-\int_{E\times\R_+}F(t,x,a)\,\dd\mathcal P}
                   \Big|\mathcal F_0\right]
 =\exp\left\{-\int_E\int_0^\infty\left(1-e^{-F(t,x,a)}\right)
                                  \dd a\dd t\dd x\right\}.
\end{equation}
This conditional Laplace functional proves that $\mathcal P$ is a
Poisson measure independent of $\mathcal F_0$, hence of $\Pi^\epsilon$.
Its restriction below $\lambda^\epsilon$ is $\mathcal P^-$, which
proves \eqref{eq:Hawkes-Poisson-equation}.  This is the marked-space
embedding corresponding to \cite[Proposition~3]{delattre2016networks}.

Starting with the zero positive-time measure, monotone Picard iteration
in \eqref{eq:Hawkes-Poisson-equation} constructs a solution as a measurable
function of $(\Pi^\epsilon,\mathcal P)$.  Each iterate has mean intensity bounded
by $c_\epsilon$ by the preceding induction.  Its increasing limit is
locally finite, satisfies \eqref{eq:Hawkes-local-construction-bound}, and solves
the equation by monotone convergence.  Let $N,\widetilde N$ be two
solutions in \eqref{eq:Hawkes-mean-class} driven by the same
$(\Pi^\epsilon,\mathcal P)$, with intensities $\lambda,\widetilde\lambda$, and set
\[
  q(t,x):=\E\left[\left|\lambda(t,x)-\widetilde\lambda(t,x)\right|\right],
\]
and recall $H(s):=\int_{\Rd} h(s,z)\dd z$ from \eqref{H:s:function}.
The common-noise coupling gives
\[
\E\left[\left|N-\widetilde N\right|(\dd s,\dd y)\right]=q(s,y)\dd s\dd y.
\]
The identical prescribed past cancels, and Tonelli's theorem yields,
for almost every $(t,x)$,
\begin{equation}\label{eq:Hawkes-Volterra-uniqueness}
  q(t,x)\leq\int_0^t\int_{\Rd}
       h(t-s,x-y)q(s,y)\dd y\dd s.
\end{equation}
Fix $T<\infty$.  For $b>0$, let
\[
q_T^*:=\operatorname*{ess\,sup}_{(t,x)\in(0,T]\times\Rd}e^{-bt}q(t,x),
\]
which is finite by \eqref{eq:Hawkes-mean-class}.  Since $H$ is integrable
and has no atom at zero, choose $b$ large enough such that
\[
k_b:=\int_0^T e^{-bs}H(s)\dd s<1.
\]
Equation \eqref{eq:Hawkes-Volterra-uniqueness} gives $q_T^*\leq k_bq_T^*$,
so $q=0$ almost everywhere on $(0,T]\times\Rd$.  The common-noise
identity then gives $N=\widetilde N$ almost surely on every bounded box.
For $Q_n:=(0,n]\times B_n$, this gives
\[
 \mathbb P\left(\bigcap_{n=1}^\infty
          \left\{\left|N-\widetilde N\right|(Q_n)=0\right\}\right)=1,
 \qquad \bigcup_{n=1}^\infty Q_n=E,
\]
Thus $N=\widetilde N$ on all of $E$ on a single probability-one event,
which proves pathwise uniqueness.
  Any embedded solution in
\eqref{eq:Hawkes-mean-class} therefore equals the Picard solution on its
driving space, which also proves uniqueness in law and identifies it
with the cluster construction in the class \eqref{eq:Hawkes-mean-class}.
This completes the proof.
\end{proof}

Define the rescaled random measure $\cN_\epsilon$ on $E$ by
\begin{equation}\label{eq:scaled-measure}
  \ip{\cN_\epsilon}{\varphi}
  :=\epsilon^{3+\alpha}
  \int_E \varphi(\epsilon^2t,\epsilon x)
  N^\epsilon(\dd t,\dd x),
  \qquad \varphi\in C_c(E).
\end{equation}
We also write $\overline E=[0,\infty)\times\Rd$ when the initial
time is included.
The exponent is $(\alpha+1-d)+(d+2)=\alpha+3$: the microscopic
density exponent is added to the temporal and spatial Jacobian
exponents.  Thus both the first moment and the fluctuation threshold are
normalized consistently in every dimension.  We use the vague topology
on $\cM_{\mathrm{loc}}(E)$ defined in Section~\ref{subsec:intro-scaling}; see
\cite{kallenberg2017random}.  This is the appropriate topology because
$\lambda_0$ need not be integrable and the total mass on all of space may
be infinite.

\section{Critical Case}
\label{sec:critical}

We now assume that the branching ratio is exactly one and use the
finite-moment kernel conditions of Section~\ref{sec:model}.
The density parameter $\alpha$ distinguishes deterministic,
stochastic, and sparse limits within this critical model; it does not
change the branching ratio.  Table~\ref{tab:results-overview} compares
these limits with the nearly critical results.

\subsection{Limit theorems and observation scales}
\label{subsec:critical-results}

The event-measure normalization \eqref{eq:scaled-measure} is used for
every $\alpha>0$ and fixed $d\geq1$.  All three regimes require only
\eqref{eq:H1}--\eqref{eq:H4}.

Let $M_\epsilon(\varphi)$ denote the compensated event martingale:
\begin{equation}\label{eq:Q-direct}
  M_\epsilon(\varphi)
  :=\epsilon^{3+\alpha}
  \int_E \varphi(\epsilon^2t,\epsilon x)
  \left[N^\epsilon(\dd t,\dd x)
       -\lambda^\epsilon(t,x)\dd t\dd x\right].
\end{equation}
On the macroscopic equation scale this is $\epsilon^{-2}M_\epsilon(\varphi)$.
For the homogeneous profile $\lambda_0\equiv c>0$, the mean microscopic
intensity is exactly $c\epsilon^{-(\alpha+1-d)}$.  Consequently,
\begin{equation}\label{eq:homogeneous-bracket}
  \E\!\left[\left\langle \epsilon^{-2}M_\epsilon(\varphi)
                  \right\rangle_{\mathrm{pred}}\right]
  =c\epsilon^{\alpha-1}\int_E(\varphi(t,x))^2\dd t\dd x.
\end{equation}

The mean event measure provides a common reference for all density
regimes.  The next proposition identifies this reference as heat flow,
independently of whether the fluctuations vanish, persist, or become sparse.

\begin{proposition}[Annealed heat-flow limit]
\label{prop:annealed-heat}
Assume \eqref{eq:H1}--\eqref{eq:H4}.  For every $\alpha>0$, the mean
measures satisfy
\begin{equation}\label{eq:annealed-heat}
  \E\left[\cN_\epsilon\right]
  \xrightarrow[\epsilon\downarrow0]{v}
  \rho(t,x)\dd t\dd x,
  \qquad \rho(t,\cdot)=P_t^D\lambda_0,
\end{equation}
where $D:=\Sigma/(2m_1)$ and the convergence is vague.
\end{proposition}

\subsubsection{Above the fluctuation threshold: the deterministic
hydrodynamic limit
(\texorpdfstring{$\alpha>1$}{alpha > 1})}

Above the fluctuation threshold, the compensated noise vanishes on the
macroscopic equation scale.  The following theorem shows that the heat
flow describing the mean also describes the limiting event measure.

\begin{theorem}[Hydrodynamic limit]\label{thm:main}
Suppose \eqref{eq:H1}--\eqref{eq:H4} hold, and let
$\lambda_0:\Rd\to\R_+$ be bounded and uniformly continuous.  If
$\alpha>1$, then
\begin{equation}\label{eq:main-convergence}
  \cN_\epsilon
  \xrightarrow[\epsilon\downarrow0]{\mathbb{P}}
  \rho(t,x)\dd t\dd x,
\end{equation}
in the vague topology on $\cM_{\mathrm{loc}}(E)$, where
\begin{equation}\label{eq:rho-semigroup}
  \rho(t,\cdot)=P_t^D\lambda_0.
\end{equation}
Equivalently, $\rho$ is the unique bounded weak solution of
\begin{equation}\label{eq:heat-equation}
  \begin{cases}
  m_1\partial_t\rho(t,x)=\frac12\Sigma:D_x^2\rho(t,x),
    & t>0,\ x\in\Rd,\\[0.4em]
  \rho(0,x)=\lambda_0(x),
    & x\in\Rd.
  \end{cases}
\end{equation}
In particular,
\begin{equation}\label{eq:heat-kernel}
  \rho(t,x)
  =\int_{\Rd}
  \frac{\exp\!\left[-\dfrac{m_1}{2t}
       (x-y)^{\top}\Sigma^{-1}(x-y)\right]}
       {(2\pi t/m_1)^{d/2}\sqrt{\det\Sigma}}
  \lambda_0(y)\dd y.
\end{equation}
\end{theorem}

\begin{remark}[The diffusion matrix]
Along a genealogical line, reproduction produces i.i.d.\ joint increments
$(S,Z)$ with law $h(s,z)\dd s\dd z$.  The mean spatial displacement is
zero, its covariance per generation is $\Sigma$, and the mean physical
time per generation is $m_1$.  The covariance per unit time is therefore
$\Sigma/m_1$, giving generator $\cA=(2m_1)^{-1}\Sigma:D_x^2$.
Independence of $S$ and $Z$ is unnecessary.  If $\Sigma=\sigma^2 I_d$,
then $\cA=(\sigma^2/(2m_1))\Delta$; for $d=1$ this is
$(m_2/(2m_1))\partial_{xx}$.
\end{remark}

\subsubsection{At the fluctuation threshold: the stochastic
occupation-measure limit
(\texorpdfstring{$\alpha=1$}{alpha = 1})}

At $\alpha=1$, \eqref{eq:homogeneous-bracket} stays of order one.  The
limiting object is most cleanly stated as a measure-valued process.  Let
$X=(X_t)_{t\geq0}$ denote the canonical locally finite Dawson--Watanabe
process constructed in Lemma~\ref{lem:locally-finite-SBM}, with initial measure
\begin{equation}\label{eq:SBM-initial}
  X_0(\dd x)=m_1\lambda_0(x)\dd x,
\end{equation}
and martingale problem
\begin{align}
  \cM_t^X(\psi)
  &:={\ip{X_t}{\psi}}-{\ip{X_0}{\psi}}
    -\int_0^t\ip{X_s}{\cA\psi}\dd s,
  \label{eq:SBM-martingale}\\
  \left\langle\cM^X(\psi)\right\rangle_t
  &=\frac1{m_1}\int_0^t\ip{X_s}{\psi^2}\dd s,
  \qquad \psi\in C_c^2(\Rd).
  \label{eq:SBM-bracket}
\end{align}
For the bounded initial density used here, the construction has continuous
paths in $\cM_{\mathrm{loc}}(\Rd)$ and finite weighted mass
$\ip{X_t}{w}$ for $w(x)=(1+|x|^2)^{-p}$, $p>d/2$, although its total
mass may be infinite.  Lemma
\ref{lem:locally-finite-SBM} proves the exhaustion construction and its
independence of the exhaustion, and derives the occupation Laplace
functional.  The finite-measure starting point is
\cite[Theorem~II.5.1]{perkins2002superprocess}.

At the threshold $\alpha=1$, branching fluctuations persist in the limit.
Because the observable counts events over physical time, we show
in the next theorem that the limiting
object is the occupation measure of the superprocess $X$ defined above.

\begin{theorem}[Critical occupation-measure limit]
\label{thm:critical}
Assume \eqref{eq:H1}--\eqref{eq:H4} and let $\alpha=1$.  Then
\begin{equation}\label{eq:critical-convergence}
  \cN_\epsilon
  \xRightarrow[\epsilon\downarrow0]{}
  \nu_X(\dd t,\dd x):=\frac1{m_1}X_t(\dd x)\dd t,
\end{equation}
in the vague topology on $\cM_{\mathrm{loc}}(E)$.
\end{theorem}

The proof also gives a triangular-array extension in every dimension;
its one-dimensional case is used in the $\alpha<1$ local theorem.  Let
$\ell_\epsilon,\ell:\Rd\to\R_+$ be bounded continuous functions such that
\begin{equation}\label{eq:varying-profile-assumption}
  \sup_{0<\epsilon\leq1}\norm{\ell_\epsilon}_\infty<\infty,
  \qquad
  \ell_\epsilon\xrightarrow[\epsilon\downarrow0]{}\ell
  \quad\text{locally uniformly on }\Rd.
\end{equation}
Construct the critical process using the negative-time intensity
\begin{equation}\label{eq:varying-past-intensity}
  \epsilon^{d-2}\ell_\epsilon(\epsilon y)\dd s\dd y,
\end{equation}
and denote its rescaled event measure, with the normalization
$\epsilon^4$, by $\cN_\epsilon^{\ell_\epsilon}$.

We also need stability under changes in the initial profile.
The extension in the following theorem allows locally uniform profile convergence
and will be used in the local observation limit below.

\begin{theorem}[Critical limit for varying locally finite profiles]
\label{thm:critical-varying}
Under \eqref{eq:H1}--\eqref{eq:H4} and
\eqref{eq:varying-profile-assumption},
\begin{equation}\label{eq:varying-critical-convergence}
  \cN_\epsilon^{\ell_\epsilon}
  \xRightarrow[\epsilon\downarrow0]{}
  \frac1{m_1}X_t^\ell(\dd x)\dd t,
\end{equation}
in $\cM_{\mathrm{loc}}(E)$, where $X^\ell$ is the locally finite
Dawson--Watanabe process with initial measure
$m_1\ell(x)\dd x$ and characteristics
\eqref{eq:SBM-martingale}--\eqref{eq:SBM-bracket}.
\end{theorem}

When $d=1$, write $m_2=\Sigma_{11}$.  Lemma
\ref{lem:locally-finite-SBM} constructs, on a filtered probability
space, a nonnegative predictable density $\rho$ with
$X_t(\dd x)=m_1\rho(t,x)\dd x$ for $t>0$, together with an adapted
space--time white noise $W$ with control measure $\dd s\dd x$.
The martingale problem
\eqref{eq:SBM-martingale}--\eqref{eq:SBM-bracket} then takes the
following form: for each $\psi\in C_c^2(\R)$, almost surely for all
$t\geq0$,
\begin{align*}
 m_1\ip{\rho_t}{\psi}
 &=m_1\ip{\lambda_0}{\psi}
   +\frac{m_2}{2}\int_0^t\ip{\rho_s}{\psi''}\dd s\\
 &\qquad+\int_0^t\int_{\R}
           \psi(x)\sqrt{\rho(s,x)}\,W(\dd s,\dd x).
\end{align*}
This is the spatially weak formulation: the second spatial derivative
acts on the test function, and the stochastic integral is an It\^o
integral against white noise.  Its
square integrability follows from \eqref{eq:SBM-local-mean}:
\[
 \E\left[\int_0^t\int_{\R}\psi(x)^2\rho(s,x)\dd x\dd s\right]
 \leq t\|\lambda_0\|_\infty\|\psi\|_{L^2(\R)}^2<\infty.
\]
The quadratic variation of this stochastic integral is
$\int_0^t\ip{\rho_s}{\psi^2}\dd s$, in agreement with
\eqref{eq:SBM-bracket} under $X_t(\dd x)=m_1\rho(t,x)\dd x$.
The initial condition means
$\ip{\rho_t}{\psi}\to\ip{\lambda_0}{\psi}$ as $t\downarrow0$,
almost surely, for every $\psi\in C_c(\R)$, by the vague continuity
of $X$.

We write the test-function identity in stochastic heat equation
notation as
\begin{equation}\label{eq:candidate-spde}
  m_1\dd\rho_t
  =\frac{m_2}{2}\partial_{xx}\rho_t\dd t
   +\sqrt{\rho_t}\,W(\dd t,\dd x),
  \qquad \rho_0=\lambda_0.
\end{equation}
Thus \eqref{eq:candidate-spde} is the density
representation of the limit in Theorem~\ref{thm:critical} when $d=1$.
For Theorem~\ref{thm:critical-varying}, replace $\lambda_0$ by $\ell$.
The finite-measure density and white-noise representation used in the
construction are provided by
\cite[Theorem~III.4.2(a)--(b)]{perkins2002superprocess}.

\subsubsection{Below the fluctuation threshold: sparse cluster limit
(\texorpdfstring{$\alpha<1$}{alpha < 1})}

When $\alpha<1$, the bracket in \eqref{eq:homogeneous-bracket}
diverges for every nonzero test function.
This proves that the compensated term does not vanish in $L^2$ under the
present normalization and that the high-density law-of-large-numbers proof
cannot extend across the threshold.

Equation \eqref{eq:homogeneous-bracket} is a noise-scaling statement, not
an identification of the limit.  A divergent second moment alone neither
proves non-tightness nor specifies a replacement limit.  Rare large
critical clusters can instead produce a zero limit with loss of uniform
integrability.  The next results make this mechanism precise.
Let
\begin{equation}\label{eq:sparse-parameter}
  r_\epsilon:=\epsilon^{1-\alpha},
  \qquad
  \cY_\epsilon:=r_\epsilon\cN_\epsilon.
\end{equation}
Thus $\cY_\epsilon$ assigns mass $\epsilon^4$ to each event.  Let
$\mathbb N_x$ be the canonical measure of the Dawson--Watanabe process
with characteristics \eqref{eq:SBM-martingale}--\eqref{eq:SBM-bracket},
started canonically at $x$
\cite[Theorem~II.7.3]{perkins2002superprocess}, and write
\begin{equation}\label{eq:canonical-occupation}
  \nu_X(\dd t,\dd x):=m_1^{-1}X_t(\dd x)\dd t.
\end{equation}
Define the sigma-finite measure on
$\cM_{\mathrm{loc}}(E)\setminus\{0\}$
\begin{equation}\label{eq:canonical-mixture}
  \mathbb Q_{\lambda_0}(B)
  :=m_1\int_{\Rd}\lambda_0(x)
       \mathbb N_x(\nu_X\in B)\dd x.
\end{equation}
Equivalently, for $g\in C_c^+(E)$,
\begin{equation}\label{eq:canonical-Laplace}
  \int_{\cM_{\mathrm{loc}}(E)\setminus\{0\}}
       \left(1-e^{-\ip{\nu}{g}}\right)\mathbb Q_{\lambda_0}(\dd\nu)
  =:\Psi_{\lambda_0}(g)
  =m_1\int_{\Rd}\lambda_0(x)u_g(0,x)\dd x,
\end{equation}
where $u_g$ solves \eqref{eq:log-Laplace-PDE} with source $g$.
The identity follows from the canonical Poisson decomposition in
\cite[Theorem~II.7.3(c)]{perkins2002superprocess} and the finite-measure
occupation calculation in Lemma~\ref{lem:locally-finite-SBM}.
The first-moment bound in the proof of Lemma
\ref{lem:canonical-parent-clan} also verifies sigma-finiteness: for a
countable compact exhaustion and nonnegative test functions $a_j$ equal to one on
its members, the sets $\{\nu:\ip{\nu}{a_j}>1/k\}$ have finite
$\mathbb Q_{\lambda_0}$-mass and cover every nonzero Radon measure.

Below the threshold, local convergence to zero coexists with a
nontrivial limiting mean.  The canonical occupation exponent defined
above captures the first-order contribution of rare clans to the
Laplace functional as in the following theorem.

\begin{theorem}[Sparse limit and canonical Laplace expansion]
\label{thm:sparse-long}
Assume \eqref{eq:H1}--\eqref{eq:H4}, let $0<\alpha<1$, and use the
original Poisson past \eqref{eq:past-intensity} and the original parabolic
scale \eqref{eq:scaled-measure}.  Then:
\begin{enumerate}
  \item for every $g\in C_c^+(E)$,
  \begin{equation}\label{eq:sparse-Laplace-first-order}
    -\frac1{r_\epsilon}
      \log\E\left[e^{-\ip{\cY_\epsilon}{g}}\right]
    \xrightarrow[\epsilon\downarrow0]{} \Psi_{\lambda_0}(g);
  \end{equation}
  \item in the vague topology,
  \begin{equation}\label{eq:sparse-zero}
    \cN_\epsilon\xrightarrow[\epsilon\downarrow0]{\mathbb{P}}0,
    \qquad
    \E\left[\cN_\epsilon\right]\xrightarrow[\epsilon\downarrow0]{v}\rho(t,x)\dd t\dd x.
  \end{equation}
\end{enumerate}
\end{theorem}

If $f\in C_c^+(E)$ satisfies
$\int_E f(t,x)\rho(t,x)\dd t\dd x>0$, then
$\ip{\cN_\epsilon}{f}\to0$ in probability while its expectation has a
strictly positive limit.  The family is therefore not uniformly
integrable.  Section~\ref{subsec:critical-sparse-proof} proves this conclusion by estimates for growing
Laplace sources.

\subsubsection{Observation scales and local recovery}

\paragraph*{Dimension-dependent observation scales}
\label{subsec:phase-diagram}

Let $\beta\in\R$ be a microscopic density exponent and let
$\gamma>0$ be a spatial observation exponent.  Write
$\Pi^{\epsilon;\beta,\gamma}$ for a Poisson random measure on
$(-\infty,0)\times\Rd$ with intensity
\begin{equation}\label{eq:gamma-past}
  \epsilon^{-\beta}\lambda_0(\epsilon^\gamma y)\dd s\dd y.
\end{equation}
Let $N^{\epsilon;\beta,\gamma}$ be the corresponding Hawkes process and
normalize its event measure by
\begin{equation}\label{eq:gamma-measure}
  \ip{\cN_\epsilon^{\beta,\gamma}}{\varphi}
  :=\epsilon^{\beta+(d+2)\gamma}
  \int_E\varphi\left(\epsilon^{2\gamma}t,\epsilon^\gamma x\right)
       N^{\epsilon;\beta,\gamma}(\dd t,\dd x).
\end{equation}
The spatial scaling of the initial profile is matched to the observation
scale. Changing $\gamma$ therefore changes the prescribed Poisson past
and the observation window.  After division by the generator scale $\epsilon^{2\gamma}$,
the expected martingale bracket is bounded by a constant times
\begin{equation}\label{eq:gamma-bracket}
  \epsilon^{\beta+(d-2)\gamma}.
\end{equation}
For a nonzero homogeneous profile this is its exact order.

The sign of $\beta+(d-2)\gamma$ determines the fluctuation regime
for the matched family.  Rescaling the three limit theorems gives
the phase diagram in the following proposition.

\begin{proposition}[Phase diagram with matched past scaling]
\label{prop:gamma-phase}
Assume \eqref{eq:H1}--\eqref{eq:H4}, $\gamma>0$, and
$\beta+(d-1)\gamma>0$.  For the family
\eqref{eq:gamma-past}--\eqref{eq:gamma-measure}:
\begin{enumerate}
  \item if $\beta+(d-2)\gamma>0$, the measure converges in probability
  to the deterministic heat-flow measure from Theorem~\ref{thm:main};
  \item if $\beta+(d-2)\gamma=0$, it converges
  in law to $m_1^{-1}X_t(\dd x)\dd t$ from Theorem~\ref{thm:critical};
  \item if $\beta+(d-2)\gamma<0$, it converges
  vaguely to zero in probability, while its mean converges to the
  heat-flow measure.
\end{enumerate}
\end{proposition}

\begin{proof}
Let $\delta=\epsilon^\gamma$ and
$a=\beta/\gamma+d-1>0$.  Then the past intensity and event mass
are respectively
\[
  \delta^{-(a+1-d)}\lambda_0(\delta y)\dd s\dd y,
  \qquad
  \epsilon^{\beta+(d+2)\gamma}=\delta^{a+3}.
\]
Also $a-1=[\beta+(d-2)\gamma]/\gamma$.  The three statements follow
from Theorems~\ref{thm:main}, \ref{thm:critical}, and
\ref{thm:sparse-long} with small parameter $\delta$ and density
parameter $a$.
\end{proof}

For the original Poisson past, $\beta=\alpha+1-d$.  Changing only the
observation window gives the same bracket exponent in the homogeneous
case, namely $\alpha+1-d+(d-2)\gamma$.  For $0<\alpha<1$, its zero is
$\gamma=\alpha<1$ when $d=1$; when $d=2$ the exponent is $\alpha-1$
for every $\gamma$; and when $d\geq3$ its zero is
$\gamma=1+(1-\alpha)/(d-2)>1$.  The last value corresponds to a longer
window and does not yield the local profile freezing used below.
Thus recovery by a shorter parabolic window is a one-dimensional result.

\paragraph*{One-dimensional case: matched past scaling and local
profile freezing}
\label{subsec:one-dimensional-local}

For the local recovery results below assume $d=1$, so that $\beta=\alpha$ in
\eqref{eq:gamma-past}.  A nontrivial stochastic limit in the sparse regime
is recovered by shortening the observation exponent from one to $\alpha$.

In particular, when the original density exponent satisfies $\alpha<1$,
the natural noise-preserving choice is
\begin{equation}\label{eq:alpha-natural-scaling}
  \gamma=\alpha,
  \qquad
  t_{\mathrm{mac}}=\epsilon^{2\alpha}t,
  \qquad
  x_{\mathrm{mac}}=\epsilon^\alpha x,
  \qquad
  \text{event mass}=\epsilon^{4\alpha}.
\end{equation}
To retain the full initial profile on this new scale, the intensity of the prescribed Poisson past must
be $\epsilon^{-\alpha}\lambda_0(\epsilon^\alpha y)\dd s\dd y$.
Thus Proposition~\ref{prop:gamma-phase} proves a critical limit for the
family with a modified Poisson past, not for the original family observed through
a shorter window.  The original family nevertheless has a local limit on
the same noise-preserving scale.  Fix a macroscopic observation point
$x_\star\in\R$ and define the recentered measure
\begin{equation}\label{eq:local-alpha-measure}
  \ip{\widehat{\cN}^{\alpha,x_\star}_\epsilon}{\varphi}
  :=\epsilon^{4\alpha}
  \int_E
  \varphi\!\left(\epsilon^{2\alpha}t,
          \epsilon^\alpha\left(x-\frac{x_\star}{\epsilon}\right)\right)
  N^\epsilon(\dd t,\dd x).
\end{equation}

On this shorter window, the initial profile becomes
locally constant around $x_\star$.  The next theorem identifies the
recovered critical occupation law and its initial density
$m_1\lambda_0(x_\star)$.

\begin{theorem}[Local critical limit below the threshold]
\label{thm:local-critical}
Assume $d=1$, \eqref{eq:H1}--\eqref{eq:H4}, and $0<\alpha<1$.  Use the
original Poisson past \eqref{eq:past-intensity}.  For every $x_\star\in\R$,
\begin{equation}\label{eq:local-alpha-limit}
  \widehat{\cN}^{\alpha,x_\star}_\epsilon
  \xRightarrow[\epsilon\downarrow0]{}
  \frac1{m_1}X_t^{x_\star}(\dd x)\dd t,
\end{equation}
in the vague topology on $\cM_{\mathrm{loc}}(E)$, where
$X^{x_\star}$ is the locally finite Dawson--Watanabe process with the
martingale characteristics \eqref{eq:SBM-martingale}--\eqref{eq:SBM-bracket}
and homogeneous initial measure
\begin{equation}\label{eq:local-alpha-initial}
  X_0^{x_\star}(\dd x)=m_1\lambda_0(x_\star)\dd x.
\end{equation}
\end{theorem}

\begin{proof}
Let $\delta=\epsilon^\alpha$ and translate the microscopic coordinate by
$x_\star/\epsilon$.  Writing $y=x_\star/\epsilon+\widetilde y$ in
\eqref{eq:past-intensity}, with $d=1$, gives the past intensity
\begin{align}
 \epsilon^{-\alpha}\lambda_0\left(x_\star+\epsilon\widetilde y\right)
       \dd s\dd\widetilde y
 &=\delta^{-1}\lambda_{0,\delta}^{x_\star}\left(\delta\widetilde y\right)
       \dd s\dd\widetilde y,\nonumber\\
 \lambda_{0,\delta}^{x_\star}(x)
 &:=\lambda_0\!\left(x_\star+\delta^{1/\alpha-1}x\right).
 \label{eq:frozen-profile}
\end{align}
Here the varying family, sometimes called a triangular array of
profiles, means that the function $\lambda_{0,\delta}^{x_\star}$ is
allowed to depend on the small parameter $\delta$.  For every $R<\infty$,
\begin{equation}\label{eq:frozen-profile-convergence}
 \left\|\lambda_{0,\delta}^{x_\star}\right\|_\infty\leq\|\lambda_0\|_\infty,
 \qquad
 \sup_{|x|\leq R}
 \left|\lambda_{0,\delta}^{x_\star}(x)-\lambda_0(x_\star)\right|
 \leq\omega_{\lambda_0}\left(R\delta^{1/\alpha-1}\right)\longrightarrow0,
\end{equation}
where $\omega_{\lambda_0}(a):=\sup_{|x-y|\leq a}
|\lambda_0(x)-\lambda_0(y)|\to0$ as $a\downarrow0$.
Since $1/\alpha-1>0$, \eqref{eq:frozen-profile-convergence} verifies
\eqref{eq:varying-profile-assumption} with limiting profile
$\ell(x)=\lambda_0(x_\star)$.  Translation invariance of $h$ preserves
the offspring law, while the observed time, space, and event mass are
$\delta^2$, $\delta$, and $\delta^4$, respectively, by
\eqref{eq:local-alpha-measure}.  Theorem~\ref{thm:critical-varying},
applied with parameter $\delta$, proves \eqref{eq:local-alpha-limit}.
\end{proof}

Thus, changing only the observation window does not retain the whole
function $\lambda_0$; it zooms into one
macroscopic point and freezes the profile there.  Matching the past to the observation scale and
local observation answer two different questions.

Keeping the original $\gamma=1$ window when $\alpha<1$ answers a different
question from the local theorem.  Theorem~\ref{thm:sparse-long} shows that
this longer window has vanishing local mass in probability despite its
annealed heat profile.  The shorter scale
\eqref{eq:alpha-natural-scaling} gives the local
critical limits of Theorem~\ref{thm:local-critical} and Proposition~\ref{prop:gamma-phase}.

\subsection{Proof of the annealed and deterministic limits}
\label{subsec:critical-deterministic-proof}

\subsubsection{Uniform moment bounds and tightness}

Write
\[
  \bar\lambda^\epsilon(t,x):=\E[\lambda^\epsilon(t,x)].
\]
The generation bound from the cluster construction yields a uniform
first-moment estimate for the rescaled event measure in the following lemma, which will control
both compact-set masses and errors in test functions.

\begin{lemma}[Uniform mean bound]\label{lem:mean}
For all $t\geq0$ and $x\in\Rd$,
\begin{equation}\label{eq:mean-bound}
  \bar\lambda^\epsilon(t,x)
  \leq \epsilon^{-(\alpha+1-d)}\norm{\lambda_0}_\infty.
\end{equation}
Consequently, for every nonnegative $f\in C_c(E)$,
\begin{equation}\label{eq:measure-mean-bound}
  \E\left[\ip{\cN_\epsilon}{f}\right]
  \leq \norm{\lambda_0}_\infty\int_E f(t,x)\dd t\dd x.
\end{equation}
\end{lemma}

\begin{proof}
With $b_\epsilon^0$ from \eqref{eq:construction-generation-bound}, define
\[
 \left(\mathsf H_+v\right)(t,x):=\int_0^t\int_{\Rd}
        h(t-s,x-y)v(s,y)\dd y\dd s.
\]
The canonical cluster version of the mean intensity is the nonnegative
Borel function
\begin{equation}\label{eq:mean-cluster-series}
 \bar\lambda^\epsilon(t,x)
 :=\sum_{n=0}^\infty\left(\mathsf H_+^n b_\epsilon^0\right)(t,x)
 =\lim_{n\to\infty}\bar\lambda_{\epsilon,n}(t,x)
 \leq\epsilon^{-(\alpha+1-d)}\|\lambda_0\|_\infty.
\end{equation}
Taking expectations in the Poisson equation
\eqref{eq:Hawkes-Poisson-equation} gives the mean compensator identity
\begin{equation}\label{eq:mean-compensator-density}
 \E\left[N^\epsilon(\dd t,\dd x)\right]
 =\E\left[\lambda^\epsilon(t,x)\right]\dd t\dd x.
\end{equation}
Summing the expected generation densities in
\eqref{eq:construction-generation-bound} identifies
\eqref{eq:mean-cluster-series} as a Borel version of the density in
\eqref{eq:mean-compensator-density}.
The bound follows from \eqref{eq:construction-generation-bound} and
proves \eqref{eq:mean-bound} for this version.  The change of variables
$(t,x)\mapsto(t/\epsilon^2,x/\epsilon)$ in
\eqref{eq:scaled-measure} gives
\begin{align*}
 \E\left[\ip{\cN_\epsilon}{f}\right]
 =\int_E f(t,x)\epsilon^{\alpha+1-d}
       \bar\lambda^\epsilon\left(t/\epsilon^2,x/\epsilon\right)\dd t\dd x
 \leq\|\lambda_0\|_\infty\int_E f(t,x)\dd t\dd x,
\end{align*}
which is \eqref{eq:measure-mean-bound}.
This completes the proof.
\end{proof}

The first-moment estimate controls the mass on each set of a compact
exhaustion of $E$. This gives tightness in the vague topology in the following proposition.

\begin{proposition}[Tightness]\label{prop:tightness}
The family
$\left\{\operatorname{Law}(\cN_\epsilon):0<\epsilon\leq1\right\}$ is tight in
$\cM_{\mathrm{loc}}(E)$.
\end{proposition}

\begin{proof}
Choose compact sets $K_j$ with $K_j\subset K_{j+1}^{\circ}$ and
$\bigcup_jK_j=E$.  For $\delta_0>0$, choose $a_j$ so large that
\[
  \frac{\norm{\lambda_0}_\infty|K_j|}{a_j}
  \leq\delta_0\,2^{-j}.
\]
Lemma~\ref{lem:mean}, Markov's inequality, and the union bound give
\[
  \sup_{\epsilon\leq1}  \mathbb{P}\left(\cN_\epsilon(K_j)>a_j\text{ for some }j\right)
  \leq\delta_0.
\]
The set
$\left\{\mu\in\cM_{\mathrm{loc}}(E):\mu(K_j)\leq a_j\text{ for all }j\geq1\right\}$
is relatively compact in the vague topology.  Indeed, from any
sequence in this set, diagonal extraction gives convergence on a
countable family uniformly dense among continuous test functions
supported on each $K_j$.  The estimate
\[
 |\ip{\mu}{f}-\ip{\mu}{g}|
 \leq a_j\|f-g\|_\infty,
 \qquad \operatorname{supp}f,\operatorname{supp}g\subset K_j,
\]
extends this convergence to every $f\in C_c(E)$.  The resulting
positive linear functional is bounded by $a_j\|f\|_\infty$ on
functions supported on $K_j$, and hence the Riesz--Markov
representation theorem represents it by a Radon measure on $E$.
Thus the closure of this set is vaguely compact.
The preceding probability estimate therefore proves tightness.
\end{proof}

The Laplace functional both determines the law of a random measure
and characterizes convergence in the vague topology.  We record the
criterion used to identify the stochastic limits below.

\begin{lemma}[Laplace functionals and vague convergence]
\label{lem:Laplace-vague}
Let $\xi_n$ and $\xi$ be random elements of $\cM_{\mathrm{loc}}(E)$.
Then
\[
 \xi_n\Rightarrow\xi
 \quad\text{if and only if}\quad
 \E\left[e^{-\ip{\xi_n}{f}}\right]\longrightarrow\E\left[e^{-\ip{\xi}{f}}\right]
 \quad\text{for every }f\in C_c^+(E).
\]
Moreover, equality of the Laplace functionals on $C_c^+(E)$ implies
equality in law.  When $\xi=\mu$ is deterministic, the equivalent
conditions also imply $\xi_n\to\mu$ in probability in the vague topology.
\end{lemma}

\begin{proof}
For $f\in C_c^+(E)$, the map
$\zeta\mapsto e^{-\ip{\zeta}{f}}$ is bounded and continuous in the
vague topology, which proves the forward implication.

We first prove that the Laplace functional determines the law.
For $f_1,\ldots,f_k\in C_c^+(E)$ and $a_1,\ldots,a_k\geq0$,
\[
 L_\xi\!\left(\sum_{i=1}^k a_if_i\right)
 =\E\left[\exp\!\left\{-\sum_{i=1}^k a_i\ip{\xi}{f_i}\right\}\right].
\]
Consequently, equality of two Laplace functionals implies equality of
the joint laws of all such pairings by uniqueness of multivariate
Laplace transforms.  To pass to random-measure laws, choose compact
sets $K_j\subset K_{j+1}^{\circ}$ exhausting $E$ and cutoffs
$\chi_j\in C_c^+(E)$ with $\chi_j=1$ on $K_j$.  Let
$\mathcal D=\{f_n:n\geq1\}\subset C_c^+(E)$ be a countable family
containing these cutoffs and uniformly dense among nonnegative test
functions supported on each $K_j$.  For measures $\zeta,\eta$ and
test functions $f,h$ supported on $K_j$,
\begin{equation}\label{eq:vague-countable-approximation}
 |\ip{\zeta-\eta}{f-h}|
 \leq\|f-h\|_\infty(\zeta(K_j)+\eta(K_j))
 \leq\|f-h\|_\infty
       \left(\ip{\zeta}{\chi_j}+\ip{\eta}{\chi_j}\right).
\end{equation}
Convergence of the coordinates indexed by $\mathcal D$ bounds the
compact-set masses through the $\chi_j$.  Approximation using
\eqref{eq:vague-countable-approximation}, followed by decomposition
into positive and negative parts, therefore gives convergence against
every test function in $C_c(E)$.  Thus these coordinates induce the
vague topology through the metric
\[
 d_{\mathcal D}(\zeta,\eta)
 =\sum_{n=1}^\infty2^{-n}
       \left(1\wedge\left|\ip{\zeta}{f_n}-\ip{\eta}{f_n}\right|\right)
\]
and generate its Borel sigma-field.  Equality of the joint laws of
the coordinates indexed by $\mathcal D$ proves equality in law.

Now assume convergence of the Laplace functionals.  For any compact
$K\subset E$, choose $\chi\in C_c^+(E)$ with $\chi\geq1$ on $K$.
For $\theta,R>0$,
\[
 \mathbb P(\xi_n(K)>R)
 \leq\mathbb P(\ip{\xi_n}{\chi}>R)
 \leq\frac{1-L_{\xi_n}(\theta\chi)}{1-e^{-\theta R}}.
\]
Since $\ip{\xi}{\chi}<\infty$ almost surely, the bounded convergence theorem
gives $L_\xi(\theta\chi)\to1$ as $\theta\downarrow0$.  Hence
\[
 \lim_{R\to\infty}\limsup_{n\to\infty}\mathbb P(\xi_n(K)>R)
 \leq1-L_\xi(\theta\chi)
 \xrightarrow[\theta\downarrow0]{}0.
\]
For every $\delta>0$, there are $R_0<\infty$ and $N\geq1$ such that
$\mathbb P(\xi_n(K)>R_0)<\delta$ for all $n\geq N$.
The finitely many random variables with $n<N$ are finite almost surely;
increasing $R_0$ also controls their tails.  Hence
$\lim_{R\to\infty}\sup_{n\geq1}\mathbb P(\xi_n(K)>R)=0$.
The compact-exhaustion argument of Proposition~\ref{prop:tightness}
therefore gives tightness of the random-measure laws.  Every
subsequential limit has Laplace functional $L_\xi$ by the forward
implication, and thus has the law of $\xi$ by the uniqueness just
proved.  This proves $\xi_n\Rightarrow\xi$.

Finally, if $\xi=\mu$ is deterministic, the Portmanteau theorem
applied to the closed set
$\{\zeta:d_{\mathcal D}(\zeta,\mu)\geq\delta\}$ gives
\[
 \limsup_{n\to\infty}
 \mathbb P\!\left(d_{\mathcal D}(\xi_n,\mu)\geq\delta\right)=0,
 \qquad\delta>0,
\]
which is convergence in probability.
\end{proof}

\subsubsection{The martingale and the forward operator}

\paragraph*{Vanishing of the compensated noise}

Define the rescaled intensity
\begin{equation}\label{eq:scaled-intensity}
  \lambda_\epsilon(t,x)
  :=\epsilon^{\alpha+1-d}
    \lambda^\epsilon\left(t/\epsilon^2,x/\epsilon\right).
\end{equation}
For $\varphi\in C_c^{1,2}([0,\infty)\times\Rd)$, the martingale
in \eqref{eq:Q-direct} can be written as
\begin{equation}\label{eq:martingale}
  M_\epsilon(\varphi)
  =\ip{\cN_\epsilon}{\varphi}
    -\int_E\varphi(t,x)\lambda_\epsilon(t,x)\dd t\dd x.
\end{equation}
This is a compensated point-process martingale integral.  Its predictable
quadratic variation is
\begin{equation}\label{eq:martingale-bracket}
  \left\langle M_\epsilon(\varphi)\right\rangle_{\mathrm{pred}}
  =\epsilon^{3+\alpha}
   \int_E\left(\varphi(t,x)\right)^2\lambda_\epsilon(t,x)\dd t\dd x.
\end{equation}
Lemma~\ref{lem:mean} therefore gives
\begin{equation}\label{eq:martingale-L2}
  \E\left[
    \left|\frac{M_\epsilon(\varphi)}{\epsilon^2}\right|^2
  \right]
  \leq C_\varphi\epsilon^{\alpha-1}.
\end{equation}
Thus, for $\alpha>1$,
\begin{equation}\label{eq:martingale-vanishes}
  \epsilon^{-2}M_\epsilon(\varphi)\xrightarrow[\epsilon\downarrow0]{}0
  \qquad\text{in }L^2.
\end{equation}

\paragraph*{The offspring operator}

Define
\begin{equation}\label{eq:forward-operator}
  K_\epsilon\varphi(t,x)
  :=\int_0^\infty\int_{\Rd}
    h(s,z)\varphi\left(t+\epsilon^2s,x+\epsilon z\right)\dd z\dd s,
\end{equation}
and
\begin{equation}\label{eq:limit-generator}
  \cL\varphi
  :=m_1\partial_t\varphi+\frac12\Sigma:D_x^2\varphi.
\end{equation}

The rescaled offspring operator determines the drift in the limit.
In the following lemma, we establish its expansion in $L^1(E)$ so that the remainder can be
controlled using the first-moment bound.

\begin{lemma}[$L^1$ generator expansion]\label{lem:generator}
For every $\varphi\in C_c^{1,2}([0,\infty)\times\Rd)$,
\begin{equation}\label{eq:L1-generator}
  \norm{
    \frac{K_\epsilon\varphi-\varphi}{\epsilon^2}
    -\cL\varphi
  }_{L^1(E)}
  \xrightarrow[\epsilon\downarrow0]{}0.
\end{equation}
\end{lemma}

\begin{proof}
By the definition of $K_\epsilon$ in \eqref{eq:forward-operator}, we have
\begin{align*}
\frac{K_\epsilon\varphi(t,x)-\varphi(t,x)}{\epsilon^2}
=\int_0^\infty\int_{\Rd}
    h(s,z)\left[\frac{\varphi\left(t+\epsilon^2s,x+\epsilon z\right)-\varphi(t,x)}{\epsilon^{2}}\right]\dd z\dd s.
\end{align*}
Split the increment as
\begin{align*}
&\varphi\left(t+\epsilon^2s,x+\epsilon z\right)-\varphi(t,x)\\
&\quad=
\left[\varphi\left(t+\epsilon^2s,x+\epsilon z\right)
      -\varphi(t,x+\epsilon z)\right]
+\left[\varphi(t,x+\epsilon z)-\varphi(t,x)\right].
\end{align*}
For $v\in L^1(E)$, let $\tau_{a,b}v(t,x):=v(t+a,x+b)$,
where $a\geq0$ and $b\in\Rd$. We have the translation continuity in $L^1$, that is,
\begin{equation}\label{eq:L1-translation-continuity}
 \|\tau_{a,b}v-v\|_{L^1(E)}\longrightarrow0
 \qquad\text{as }a+|b|\longrightarrow0.
\end{equation}
Indeed, extend $v$ by zero to $\R^{d+1}$ and approximate this
extension in $L^1$ by continuous compactly supported functions.
Translations preserve the full-space $L^1$ norm, and translations
of each approximating function converge in $L^1$ by uniform
continuity and a common compact support for small shifts.
Restriction to $E$ proves \eqref{eq:L1-translation-continuity}; also
$\|\tau_{a,b}v\|_{L^1(E)}\leq\|v\|_{L^1(E)}$ for $a\geq0$.

For fixed $(s,z)$, the fundamental theorem of calculus gives
\begin{equation}\label{eq:generator-time-quotient}
 Q_{\epsilon,s,z}(t,x)
 :=\frac{\varphi(t+\epsilon^2s,x+\epsilon z)
             -\varphi(t,x+\epsilon z)}{\epsilon^2}
 =s\int_0^1\tau_{r\epsilon^2s,\epsilon z}
                          \partial_t\varphi(t,x)\dd r.
\end{equation}
By \eqref{eq:L1-translation-continuity},
\begin{align*}
 \left\|Q_{\epsilon,s,z}-s\partial_t\varphi\right\|_1
 &\leq s\int_0^1\left\|\tau_{r\epsilon^2s,\epsilon z}\partial_t\varphi
                         -\partial_t\varphi\right\|_1\dd r\longrightarrow0,\\
 \|Q_{\epsilon,s,z}\|_1&\leq s\|\partial_t\varphi\|_1.
\end{align*}
The integrand in the first bound is at most
$2\|\partial_t\varphi\|_1$, and \eqref{eq:H3} makes
$s h(s,z)$ integrable. The dominated convergence theorem therefore
gives the time contribution $m_1\partial_t\varphi$ in $L^1(E)$.

For the spatial part, subtract the centered linear term:
$\varphi(t,x+\epsilon z)-\varphi(t,x)
  -\epsilon z\cdot\nabla_x\varphi(t,x)$.
Let
\[
 R_{\epsilon,z}(t,x):=
 \frac{\varphi(t,x+\epsilon z)-\varphi(t,x)
             -\epsilon z\cdot\nabla_x\varphi(t,x)}{\epsilon^2}
 =\int_0^1(1-r)z^\top D_x^2\varphi(t,x+r\epsilon z)z\dd r.
\]
Taylor's integral formula and \eqref{eq:L1-translation-continuity} give
\begin{align*}
 \left\|R_{\epsilon,z}-\frac12z^\top D_x^2\varphi z\right\|_1
 &\leq |z|^2\int_0^1(1-r)
 \left\|D_x^2\varphi(\cdot,\cdot+r\epsilon z)-D_x^2\varphi\right\|_1\dd r
 \longrightarrow0,\\
 \|R_{\epsilon,z}\|_1&\leq\frac12|z|^2\left\|D_x^2\varphi\right\|_1,
\end{align*}
where the matrix norm in these integrals is the operator norm.
By \eqref{eq:H2}, the linear term satisfies
\[
 \int_0^\infty\int_{\Rd}h(s,z)
           \epsilon^{-1}z\cdot\nabla_x\varphi(t,x)\dd z\dd s
 =\epsilon^{-1}\left(\int_0^\infty\int_{\Rd}zh(s,z)\dd z\dd s\right)
                  \cdot\nabla_x\varphi(t,x)=0.
\]

Since
$\int_0^\infty\int_{\Rd} |z|^2h(s,z)\dd z\dd s=\operatorname{tr}\Sigma$,
dominated convergence theorem gives the spatial limit
$\frac12\Sigma:D_x^2\varphi$ and proves \eqref{eq:L1-generator}.
\end{proof}

In the following corollary, we show that
combining the generator expansion with the mean bound makes the
remainder negligible when paired with the random event measure.

\begin{corollary}\label{cor:generator-pairing}
\begin{equation}\label{eq:generator-pairing}
  \ip{\cN_\epsilon}{
    \frac{K_\epsilon\varphi-\varphi}{\epsilon^2}
    -\cL\varphi}
  \xrightarrow[\epsilon\downarrow0]{}0
  \qquad\text{in }L^1(\mathbb{P}).
\end{equation}
\end{corollary}

\begin{proof}
Apply \eqref{eq:measure-mean-bound} to the absolute value of the integrand
and then use Lemma~\ref{lem:generator}.
\end{proof}

\subsubsection{The initial Poisson boundary layer}

The contribution to the rescaled compensator from parents before time zero
is
\begin{align}
  \mathsf I_\epsilon(\varphi)
  :={}&\epsilon^{\alpha+1-d}\int_E\varphi(t,x)
  \int_{(-\infty,0)\times\Rd}
  h\!\left(\frac{t}{\epsilon^2}-s,
           \frac{x}{\epsilon}-y\right)
  \Pi^\epsilon(\dd s,\dd y)\dd t\dd x.
  \label{eq:boundary-def}
\end{align}

The prescribed Poisson past enters the limiting equation through
an initial measure.  The estimate in the next lemma identifies this boundary term
and shows that its fluctuations vanish.

\begin{lemma}[Initial-layer law of large numbers]\label{lem:boundary}
For every $\varphi\in C_c^{1,2}([0,\infty)\times\Rd)$,
\begin{equation}\label{eq:boundary-convergence}
  \frac{\mathsf I_\epsilon(\varphi)}{\epsilon^2}
  \xrightarrow[\epsilon\downarrow0]{}
  m_1\int_{\Rd}\varphi(0,x)\lambda_0(x)\dd x
  \qquad\text{in }L^2.
\end{equation}
\end{lemma}

\begin{proof}
The Poisson intensity \eqref{eq:past-intensity} and the substitutions
$q=t/\epsilon^2$, $r=q-s$, and $z=x/\epsilon-y$ give
\begin{align}
  \E\left[\frac{\mathsf I_\epsilon(\varphi)}{\epsilon^2}\right]
  ={}&\int_0^\infty\int_{\Rd}
  \varphi(\epsilon^2q,x)
  \int_q^\infty\int_{\Rd}
  h(r,z)\lambda_0(x-\epsilon z)
  \dd z\dd r\dd x\dd q.
  \label{eq:boundary-mean}
\end{align}
Uniform continuity of $\lambda_0$, compact support of $\varphi$, and
dominated convergence theorem imply that the right-hand side of \eqref{eq:boundary-mean} converges to
\[
  \int_{\Rd}\varphi(0,x)\lambda_0(x)\dd x
  \int_0^\infty\int_q^\infty\int_{\Rd}
  h(r,z)\dd z\dd r\dd q.
\]
By Tonelli's theorem and \eqref{eq:H3},
\[
 \int_0^\infty\int_q^\infty\int_{\Rd}h(r,z)\dd z\dd r\dd q
 =\int_0^\infty\int_{\Rd}\left(\int_0^r\dd q\right)h(r,z)\dd z\dd r
 =\int_0^\infty\int_{\Rd}rh(r,z)\dd z\dd r=m_1.
\]
It remains to control the Poisson fluctuation.  Let
\[
  H(u):=\int_{\Rd} h(u,z)\dd z,
  \qquad
  \overline H(a):=\int_a^\infty H(u)\dd u.
\]
Writing $s=-a$ and changing $t=\epsilon^2q$, $x=\epsilon w$, the integrand
of the compensated past Poisson measure in
$\mathsf I_\epsilon(\varphi)/\epsilon^2$ equals $\epsilon^{\alpha+1}F_\epsilon(a,y)$, where
\[
  F_\epsilon(a,y)
  :=\int_0^\infty\int_{\Rd}
    \varphi\left(\epsilon^2q,\epsilon w\right)
    h(q+a,w-y)\dd w\dd q.
\]
Minkowski's inequality and Young's convolution inequality yield
\begin{align*}
  \norm{F_\epsilon(a,\cdot)}_2
  \leq
  \int_0^\infty
  \norm{\varphi\left(\epsilon^2q,\epsilon\,\cdot\right)}_2
  H(q+a)\dd q
  \leq C_\varphi\epsilon^{-d/2}\overline H(a).
\end{align*}
Since $\overline H\leq1$ and
$\int_0^\infty\overline H(a)\dd a=m_1$, the Poisson isometry gives
\begin{align}
  \operatorname{Var}\!\left(
    \frac{\mathsf I_\epsilon(\varphi)}{\epsilon^2}
  \right)
  \leq C\epsilon^{\alpha+1+d}
    \int_0^\infty
    \norm{F_\epsilon(a,\cdot)}_2^2\dd a
  \leq C_\varphi\epsilon^{\alpha+1}m_1
  \xrightarrow[\epsilon\downarrow0]{}0.
  \label{eq:boundary-variance}
\end{align}
Together with the mean convergence in \eqref{eq:boundary-mean} and the variance estimate \eqref{eq:boundary-variance}, this proves
\eqref{eq:boundary-convergence}. The proof is complete.
\end{proof}

\subsubsection{Identification of the mean and random limits}

Substituting the two intensity terms of \eqref{eq:hawkes-intensity}
into \eqref{eq:scaled-intensity} and applying Fubini's theorem gives,
by \eqref{eq:scaled-measure}, \eqref{eq:forward-operator}, and
\eqref{eq:boundary-def}, the exact compensator identity
\begin{equation}\label{eq:compensator}
  \int_E\varphi(t,x)\lambda_\epsilon(t,x)\dd t\dd x
  =\ip{\cN_\epsilon}{K_\epsilon\varphi}
   +\mathsf I_\epsilon(\varphi).
\end{equation}
Combining \eqref{eq:compensator} with \eqref{eq:martingale} yields
\begin{equation}\label{eq:master}
  \ip{\cN_\epsilon}{
    \frac{K_\epsilon\varphi-\varphi}{\epsilon^2}}
  +\frac{\mathsf I_\epsilon(\varphi)}{\epsilon^2}
  +\frac{M_\epsilon(\varphi)}{\epsilon^2}=0.
\end{equation}
Using test functions compactly supported in the entire positive half-line
also prevents an artificial terminal boundary term.

Let $\epsilon_n\downarrow0$.  By Proposition~\ref{prop:tightness}, a
subsequence converges in distribution to a random locally finite measure
$\nu$.  Test functions that meet $t=0$ require an additional truncation because the
vague topology is defined on the open set $E$.  Write
$C_0=\norm{\lambda_0}_\infty$.
The mean bound \eqref{eq:measure-mean-bound}, extended to
nonnegative Borel functions by monotone approximation, gives,
for $R,\delta>0$,
\begin{equation}\label{eq:initial-strip-bound}
  \sup_{0<\epsilon\leq1}\E\left[\cN_\epsilon((0,\delta)\times B_R)\right]
  \leq |B_R|C_0\delta.
\end{equation}
For every nonnegative $a\in C_c(E)$, convergence in law and truncation of
the random variable $\ip{\cN_\epsilon}{a}$ give
\[
\E\left[\ip{\nu}{a}\right]\leq C_0\int_E a(t,x)\dd t\dd x.
\]
Approximation of Borel sets and
monotone convergence therefore give the same initial-strip bound for
$\nu$.

Fix $\varphi\in C_c^{1,2}([0,\infty)\times\Rd)$ and let
\begin{equation}\label{eq:generator-remainder-definition}
 r_\epsilon^\varphi
 :=\frac{K_\epsilon\varphi-\varphi}{\epsilon^2}-\cL\varphi.
\end{equation}
Equations \eqref{eq:measure-mean-bound} and \eqref{eq:L1-generator} give
\begin{equation}\label{eq:generator-remainder-mean}
 \E\left[\left|\ip{\cN_\epsilon}{r_\epsilon^\varphi}\right|\right]
 \leq C_0\|r_\epsilon^\varphi\|_{L^1(E)}\longrightarrow0.
\end{equation}
Take a continuous cutoff $0\leq\zeta_\delta\leq1$ equal to zero on
$[0,\delta]$ and one on $[2\delta,\infty)$.  Applying this cutoff to
$\cL\varphi$ in \eqref{eq:generator-remainder-definition} gives the
exact decomposition
\begin{align}
 \ip{\cN_\epsilon}{\frac{K_\epsilon\varphi-\varphi}{\epsilon^2}}
 &=\ip{\cN_\epsilon}{\zeta_\delta\cL\varphi}
   +\ip{\cN_\epsilon}{(1-\zeta_\delta)\cL\varphi}
   +\ip{\cN_\epsilon}{r_\epsilon^\varphi}.
 \label{eq:generator-cutoff-decomposition}
\end{align}
If the spatial support of $\varphi$ is in $B_R$, then
\eqref{eq:initial-strip-bound} and its counterpart for $\nu$ imply,
for $\xi\in\{\cN_\epsilon,\nu\}$,
\begin{align}
 \E\left[\left|\ip{\xi}{(1-\zeta_\delta)\cL\varphi}\right|\right]
 &\leq\|\cL\varphi\|_\infty
          \E\left[\xi((0,2\delta)\times B_R)\right]\nonumber\\
 &\leq2|B_R|C_0\delta\|\cL\varphi\|_\infty.
 \label{eq:initial-cutoff-pairing-error}
\end{align}
For fixed $\delta$, the function $\zeta_\delta\cL\varphi$ belongs to
$C_c(E)$, so its pairing converges in law along the chosen subsequence.
Letting $\delta\downarrow0$ using
\eqref{eq:initial-cutoff-pairing-error} proves
\begin{equation}\label{eq:generator-limit-pairing}
  \ip{\cN_{\epsilon_n}}{\cL\varphi}
  \xRightarrow[\epsilon_n\downarrow0]{}\ip{\nu}{\cL\varphi}.
\end{equation}
In \eqref{eq:generator-cutoff-decomposition}, $\zeta_\delta$ multiplies
the already computed function $\cL\varphi$, while the generator acts
on $\varphi$ alone.  Thus no derivative of $\zeta_\delta$ occurs, and
the remainder estimate \eqref{eq:generator-remainder-mean} is independent
of $\delta$.
Substituting \eqref{eq:generator-remainder-definition} into
\eqref{eq:master}, and using \eqref{eq:generator-remainder-mean},
\eqref{eq:generator-limit-pairing}, \eqref{eq:boundary-convergence}, and
\eqref{eq:martingale-vanishes}, shows that every such limit satisfies,
almost surely,
\begin{equation}\label{eq:weak-limit}
  \int_E
  \left(m_1\partial_t\varphi+
        \frac12\Sigma:D_x^2\varphi\right)
  \nu(\dd t,\dd x)
  +m_1\int_{\Rd}\varphi(0,x)\lambda_0(x)\dd x=0,
\end{equation}
for all $\varphi\in C_c^{1,2}([0,\infty)\times\Rd)$.  A countable dense
family of test functions may be used to obtain a single probability-one
event on which all the required identities hold.

To identify the subsequential limit, we need uniqueness for the weak
equation in the growth class provided by the moment bounds.
Backward heat-equation test functions yield the characterization in the following proposition.

\begin{proposition}[Uniqueness of the limiting measure]\label{prop:uniqueness}
Among nonnegative Radon measures satisfying, for every $T<\infty$,
\begin{equation}\label{eq:limiting-measure-growth}
  \nu((0,T]\times B_R)\leq C_T(1+R)^{p_T},\qquad R\geq1,
\end{equation}
for some finite $C_T,p_T$, the weak equation
\eqref{eq:weak-limit} has the unique solution
\[
  \nu(\dd t,\dd x)=(P_t^D\lambda_0)(x)\dd t\dd x.
\]
\end{proposition}

\begin{proof}
Let $g\in C_c^\infty(E)$ and choose $T$ beyond its temporal support.
Define
\begin{equation}\label{eq:backward-test}
 \varphi(t,x):=\frac1{m_1}\int_t^T
       P_{s-t}^Dg(s,\cdot)(x)\dd s.
\end{equation}
Differentiating \eqref{eq:backward-test} gives
\begin{equation}\label{eq:backward-test-equation}
 \varphi(T,\cdot)=0,\qquad
 \cL\varphi=m_1\partial_t\varphi+\frac12\Sigma:D_x^2\varphi=-g.
\end{equation}
If the spatial support of $g$ is contained in $B_{R_g}$, the Gaussian
heat kernel in \eqref{eq:heat-kernel}, applied to $g$ and its spatial
derivatives, gives constants $c,C>0$ such that
\begin{equation}\label{eq:backward-test-tail}
 \sup_{0\leq t\leq T}\left(
 |\varphi(t,x)|+|\nabla_x\varphi(t,x)|+
 \left\|D_x^2\varphi(t,x)\right\|\right)\leq Ce^{-c|x|^2},
 \qquad |x|\geq2R_g+1.
\end{equation}
Indeed, spatial derivatives can be transferred to $g$ inside
$P_{s-t}^D$, and for $0<s-t\leq T$ the heat kernel between
$x$ and $B_{R_g}$ is bounded by $C_Te^{-c_T|x|^2}$ in this range.
Thus the function with Gaussian tails is precisely $\varphi$ in
\eqref{eq:backward-test}.

Let $\chi_R\in C_c^\infty(\Rd)$ equal one on $B_R$, vanish outside
$B_{2R}$, and satisfy $|\nabla\chi_R|\leq C/R$ and
$\|D_x^2\chi_R\|\leq C/R^2$.  By
\eqref{eq:backward-test-equation},
\begin{equation}\label{eq:backward-cutoff-error}
 \cL(\chi_R\varphi)=-\chi_R g+E_R,
 \qquad
 E_R:=\varphi\,\frac12\Sigma:D_x^2\chi_R
             +(\Sigma\nabla\chi_R)\cdot\nabla_x\varphi.
\end{equation}
The support of $E_R$ lies in $B_{2R}\setminus B_R$.
The polynomial growth condition \eqref{eq:limiting-measure-growth}
and the Gaussian bound \eqref{eq:backward-test-tail} imply
\begin{align*}
 \left|\int_E E_R(t,x)\nu(\dd t,\dd x)\right|
 &\leq C\left(R^{-1}+R^{-2}\right)e^{-cR^2}
            \nu((0,T]\times B_{2R})\\
 &\leq C C_T\left(R^{-1}+R^{-2}\right)(1+2R)^{p_T}e^{-cR^2}
       \longrightarrow0.
\end{align*}
Insert $\chi_R\varphi$ in \eqref{eq:weak-limit} and use
\eqref{eq:backward-cutoff-error}.  Since $\lambda_0$ is bounded and
$\varphi(0,\cdot)$ is integrable, letting $R\to\infty$ yields
\begin{align*}
 \ip{\nu}{g}
 =m_1\int_{\Rd}\varphi(0,x)\lambda_0(x)\dd x
 =\int_0^T\int_{\Rd}g(t,x)(P_t^D\lambda_0)(x)\dd x\dd t.
\end{align*}
The second equality follows from \eqref{eq:backward-test}, Fubini's
theorem, and symmetry of the heat kernel.  Smooth compactly supported
test functions determine Radon measures, which proves uniqueness.
\end{proof}

\begin{proof}[Proof of Proposition~\ref{prop:annealed-heat}]
Let $\mu_\epsilon:=\E[\cN_\epsilon]$ and $C_0:=\|\lambda_0\|_\infty$.
Equation \eqref{eq:measure-mean-bound}, extended by monotone
approximation to nonnegative Borel test functions, gives
\[
 \mu_\epsilon(A)\leq C_0|A|,\qquad(A\subset E\text{ Borel}).
\]
It follows that the deterministic measures $\mu_\epsilon$ are vaguely
relatively compact.  Let $\mu_{\epsilon_n}\to\mu$ vaguely along a
subsequence.  Testing with nonnegative functions in $C_c(E)$ and then
approximating Borel sets shows $\mu(A)\leq C_0|A|$.

For $\varphi\in C_c^{1,2}([0,\infty)\times\Rd)$, set
$G_\epsilon\varphi:=\epsilon^{-2}(K_\epsilon\varphi-\varphi)$.
Taking expectations in \eqref{eq:master} and using the zero mean of
\eqref{eq:Q-direct}, we obtain
\begin{align}
 \ip{\mu_\epsilon}{G_\epsilon\varphi}
       +\E\left[\epsilon^{-2}\mathsf I_\epsilon(\varphi)\right]=0,\nonumber
\end{align}
and
\begin{align}
 |\ip{\mu_\epsilon}{G_\epsilon\varphi-\cL\varphi}|
       \leq C_0\|G_\epsilon\varphi-\cL\varphi\|_{L^1(E)}\longrightarrow0.
 \label{eq:annealed-master-error}
\end{align}
The last limit is \eqref{eq:L1-generator} and holds for every
$\alpha>0$.
Let $\zeta_\delta$ be the initial-time cutoff used in
\eqref{eq:initial-strip-bound}, and choose $B_R$ containing the
spatial support of $\varphi$.  For
$\widetilde\mu\in\{\mu_\epsilon,\mu\}$,
\[
 |\ip{\widetilde\mu}{(1-\zeta_\delta)\cL\varphi}|
 \leq2C_0\delta|B_R|\|\cL\varphi\|_\infty.
\]
The pairing with $\zeta_\delta\cL\varphi\in C_c(E)$ converges vaguely.
Letting first $n\to\infty$ and then $\delta\downarrow0$ in
\eqref{eq:annealed-master-error}, and using the limit of
\eqref{eq:boundary-mean}, gives
\[
 \ip{\mu}{\cL\varphi}
       +m_1\int_{\Rd}\varphi(0,x)\lambda_0(x)\dd x=0.
\]
Moreover $\mu((0,T]\times B_R)\leq C_0T|B_R|$, so
Proposition~\ref{prop:uniqueness} identifies
$\mu=\rho(t,x)\dd t\dd x$.  Every subsequential limit is this measure;
hence the full family converges as asserted in \eqref{eq:annealed-heat}.
\end{proof}


\begin{proof}[Proof of Theorem~\ref{thm:main}]
Proposition~\ref{prop:tightness} provides subsequential limits in law.
The generator, boundary, and martingale limits in
\eqref{eq:generator-pairing}, \eqref{eq:boundary-convergence}, and
\eqref{eq:martingale-vanishes}, applied to \eqref{eq:master} with the
initial-strip truncation \eqref{eq:initial-strip-bound}, show that
every such limit $\nu$ satisfies \eqref{eq:weak-limit} almost surely.
We verify the growth condition needed for pathwise uniqueness.
The first-moment bound passes to $\nu$ as explained before
\eqref{eq:weak-limit}, and monotone approximation gives
\begin{align*}
 \E\left[\nu((0,n]\times B_n)\right]&\leq C_0n|B_n|\leq Cn^{d+1},\\
 \mathbb P\left(\nu((0,n]\times B_n)>n^{d+3}\right)&\leq Cn^{-2},
\end{align*}
where the second inequality follows from Markov's inequality.  Since
$\sum_{n=1}^\infty n^{-2}<\infty$, the Borel--Cantelli lemma gives a
finite random constant $C(\omega)$ such that
\[
 \nu((0,n]\times B_n)\leq C(\omega)n^{d+3}
 \qquad(n\geq1).
\]
For fixed $T<\infty$ and $R\geq1$, let
$n:=\lceil\max\{T,R,1\}\rceil$.  Monotonicity then gives
\[
 \nu((0,T]\times B_R)
 \leq C(\omega)n^{d+3}
 \leq C_T(\omega)(1+R)^{d+3}.
\]
Thus Proposition~\ref{prop:uniqueness} applies on a single
probability-one event and gives $\nu=\rho(t,x)\dd t\dd x$.
Tightness and uniqueness of all subsequential laws give
\[
 \cN_\epsilon\xRightarrow[\epsilon\downarrow0]{}
       \rho(t,x)\dd t\dd x,\qquad
 \rho(t,x)=(P_t^D\lambda_0)(x),
\]
for the full family $(\cN_\epsilon)_{0<\epsilon\leq1}$ defined in
\eqref{eq:scaled-measure}.

If $d_v$ is any metric inducing the vague topology, the Portmanteau
theorem applied to the closed set
$\left\{\xi:d_v(\xi,\rho(t,x)\dd t\dd x)\geq a\right\}$ gives, for every $a>0$,
\[
 \limsup_{\epsilon\downarrow0}
 \mathbb P\left(d_v(\cN_\epsilon,\rho(t,x)\dd t\dd x)\geq a\right)=0.
\]
This is \eqref{eq:main-convergence}.
The proof is complete.
\end{proof}


\subsection{Proof of the stochastic occupation limit}
\label{subsec:critical-stochastic-proof}

The bracket formulas \eqref{eq:martingale-bracket} and
\eqref{eq:homogeneous-bracket} suggest fluctuations of order one at
$\alpha=1$.  The exact Poisson-cluster Laplace functional will identify
their limiting law.
  Throughout this subsection
$\alpha=1$.  The cluster Laplace estimates
proved here will be applied to the sparse regime in Section~\ref{subsec:critical-sparse-proof}.

\subsubsection{The one-cluster nonlinear renewal equation}

Fix $T>0$ and $f\in C_c^\infty([0,T)\times\Rd)$ with $f\geq0$.  Let
$C^{r,y}$ denote the cluster consisting of an event at $(r,y)$ and all of
its descendants.  For a parent at the microscopic point
$(t/\epsilon^2,x/\epsilon)$, set
\begin{equation}\label{eq:cluster-Laplace}
  F_\epsilon(t,x)
  :=\E_{t/\epsilon^2,x/\epsilon}
  \left[\exp\!\left\{-\epsilon^4
  \int_{[0,\infty)\times\Rd}
       f(\epsilon^2s,\epsilon z)C^{t/\epsilon^2,x/\epsilon}
       (\dd s,\dd z)\right\}\right].
\end{equation}
For a Poisson random measure $\Pi$ with intensity $\Lambda$, a set
$A$ with $\Lambda(A)<\infty$, and a measurable function $0\leq L\leq1$,
we have the following Poisson exponential formula
\begin{equation}\label{eq:Poisson-exponential-formula}
 \E\left[\prod_{p\in\Pi\cap A}L(p)\right]
 =\exp\left\{\int_A(L(p)-1)\Lambda(\dd p)\right\},
\end{equation}
where the empty product is one.  For the direct children of the
cluster root, take $A=E$, $\Lambda(\dd s,\dd z)=h(s,z)\dd s\dd z$,
and $L(s,z)=F_\epsilon(t+\epsilon^2s,x+\epsilon z)$.
Conditional independence of their descendant clusters, followed by
\eqref{eq:Poisson-exponential-formula} and the root's own factor
$e^{-\epsilon^4f(t,x)}$, gives the exact identity
\begin{equation}\label{eq:cluster-fixed-point}
  F_\epsilon(t,x)
  =\exp\!\left\{-\epsilon^4f(t,x)
    +\int_0^\infty\int_{\Rd} h(s,z)
       \left[F_\epsilon\left(t+\epsilon^2s,x+\epsilon z\right)-1\right]
       \dd z\dd s\right\}.
\end{equation}
We take $F_\epsilon(t,\cdot)=1$ for $t\geq T$.  With
\begin{equation}\label{eq:u-epsilon}
  u_\epsilon:=\frac{1-F_\epsilon}{\epsilon^2},
\end{equation}
equation \eqref{eq:cluster-fixed-point} becomes
\begin{equation}\label{eq:nonlinear-renewal}
  \log\left(1-\epsilon^2u_\epsilon\right)
  =-\epsilon^4 f-\epsilon^2K_\epsilon u_\epsilon,
  \qquad u_\epsilon(T,\cdot)=0,
\end{equation}
where $K_\epsilon$ is the forward operator
\eqref{eq:forward-operator}.

For $0\leq r<\epsilon^{-2}$ and a bounded function $v$, extended by zero
for $t\geq T$, define
\begin{equation}\label{eq:scheme-operator}
  \mathsf S_\epsilon[r,v](t,x)
  :=\frac{-\log(1-\epsilon^2r)
           -\epsilon^2K_\epsilon v(t,x)}{\epsilon^4}-f(t,x).
\end{equation}
Thus \eqref{eq:nonlinear-renewal} is exactly
$\mathsf S_\epsilon[u_\epsilon(t,x),u_\epsilon](t,x)=0$.

A comparison principle provides a uniform bound for the nonlinear
cluster costs. This bound is the starting point for the renewal limit.

\begin{lemma}[Comparison and stability of the renewal scheme]
\label{lem:scheme-comparison}
Suppose $v,w:[0,T]\times\Rd\to[0,\epsilon^{-2})$ are bounded, vanish for
$t\geq T$, and satisfy
\begin{align}
 \mathsf S_\epsilon[v(t,x),v](t,x)&\leq0,
       \label{eq:scheme-subsolution}\\
 \mathsf S_\epsilon[w(t,x),w](t,x)&\geq0.
       \label{eq:scheme-supersolution}
\end{align}
Then
\begin{equation}\label{eq:scheme-comparison}
 v(t,x)\leq w(t,x)\qquad(0\leq t\leq T,\ x\in\Rd).
\end{equation}  Moreover, for all sufficiently small $\epsilon$, the
cluster solution satisfies
\begin{equation}\label{eq:u-uniform-bound}
  0\leq u_\epsilon(t,x)\leq
  U_f:=\sqrt{2\norm{f}_\infty}+1
  \qquad\text{on }[0,T]\times\Rd.
\end{equation}
It is continuous and vanishes on a fixed neighborhood of $t=T$.
\end{lemma}

\begin{proof}
Let $g_\epsilon(r):=-\log(1-\epsilon^2r)$.  If $a\geq b\geq0$, then
\begin{align}
  g_\epsilon(a)-g_\epsilon(b)-\epsilon^2(a-b)
  &=\int_b^a\frac{\epsilon^4r}{1-\epsilon^2r}\dd r
    \geq\frac{\epsilon^4}{2}(a-b)^2.
  \label{eq:strict-monotonicity}
\end{align}
Assume
$M:=\sup_{(t,x)\in[0,T]\times\Rd}(v(t,x)-w(t,x))>0$.
Choose $0<\delta_0<\min\left\{M/2,\epsilon^2M^2/16\right\}$ and
$(t_0,x_0)$ such that $\zeta:=v(t_0,x_0)-w(t_0,x_0)\geq M-\delta_0$.
Positivity and unit mass of $K_\epsilon$ imply
$K_\epsilon(v-w)(t_0,x_0)\leq M$.  By
\eqref{eq:scheme-operator} and \eqref{eq:strict-monotonicity},
\begin{align}
 &\mathsf S_\epsilon[v(t_0,x_0),v](t_0,x_0)
       -\mathsf S_\epsilon[w(t_0,x_0),w](t_0,x_0)\nonumber\\
 &\geq\frac{\zeta-M}{\epsilon^2}+\frac{\zeta^2}{2}
 \geq-\frac{\delta_0}{\epsilon^2}+\frac{M^2}{8}
 >\frac{M^2}{16}>0.
 \label{eq:scheme-positive-difference}
\end{align}
Equation \eqref{eq:scheme-positive-difference} contradicts
\eqref{eq:scheme-subsolution}--\eqref{eq:scheme-supersolution}, whose
difference is nonpositive.  This proves \eqref{eq:scheme-comparison}
using an approximating maximum, without requiring its attainment.

The zero function is a subsolution.  Let $\overline u=U_f$ for $t<T$ and
$\overline u=0$ for $t\geq T$.  Because
$K_\epsilon\overline u\leq U_f$ and
\[
  -\log(1-a)-a\geq\frac{a^2}{2},
  \qquad 0\leq a<1,
\]
we have
\[
  \mathsf S_\epsilon[U_f,\overline u]
  \geq \frac{U_f^2}{2}-\norm{f}_\infty>0,
\]
whenever $\epsilon^2U_f<1$.  Comparison gives
\eqref{eq:u-uniform-bound}.

The critical Poisson Galton--Watson cluster is finite almost surely, so
the cluster expectation \eqref{eq:cluster-Laplace} exists, is strictly
positive, and satisfies
the fixed-point equation \eqref{eq:cluster-fixed-point}.
  The right-hand side of
\eqref{eq:cluster-fixed-point} is continuous: convolution of a bounded
function with the $L^1$ kernel $h$ is continuous by translation
continuity of $h$ in $L^1$.  Hence $F_\epsilon$, and therefore
$u_\epsilon$, is continuous.  Finally, since $f$ has compact support in
$[0,T)\times\Rd$, there is $T_f<T$ such that $f=0$ for $t\geq T_f$.
Every descendant occurs after its parent, so a cluster rooted at
$t\geq T_f$ has zero cost in \eqref{eq:cluster-Laplace}.  Thus
$u_\epsilon=0$ on $[T_f,T]\times\Rd$. This completes the proof.
\end{proof}

We formulate the renewal argument for a general temporal clock so
that it applies to both finite-mean and heavy-tailed delays.
Exponential weights on a fixed physical-time horizon also cover
zero limiting killing.

\begin{lemma}[Renewal limits for a general temporal clock]
\label{lem:renewal-framework}
Let $(S,Z)$ have density $h$, with $S>0$, $\E[Z]=0$, and
$\E[ZZ^{\top}]=\Sigma$ finite.  Let
$\tau_\epsilon,\delta_\epsilon\downarrow0$ and suppose
\begin{equation}\label{eq:clock-assumption}
 \frac{1-\widehat H(q\tau_\epsilon)}{\delta_\epsilon}
       \xrightarrow[\epsilon\downarrow0]{}\chi(q)>0\quad(q>0),
 \qquad \widehat H(z):=\E\left[e^{-zS}\right],
\end{equation}
where $\chi$ is the Laplace exponent of a nonzero subordinator
$D_r^\chi$, with $\chi(0)=0$.  Let $0<a_\epsilon\leq1$ and
$\kappa_\epsilon=(1-a_\epsilon)/\delta_\epsilon\to\kappa\geq0$.
Define
\[
 \mathcal K_\epsilon v(t,x)
   :=\E\left[v\left(t+\tau_\epsilon S,x+\sqrt{\delta_\epsilon}Z\right)\right],\qquad
 \mathsf P_r^\chi v(t,x)
   :=\E\left[v\left(t+D_r^\chi,x+B_r^\Sigma\right)\right],
\]
where $D^\chi$ and $B^\Sigma$ are independent in the limiting motion.
For $b\geq0$, let
\begin{equation}\label{eq:renewal-resolvents}
 \mathcal R_{\epsilon,b}
  :=\frac{\delta_\epsilon}{1+b\delta_\epsilon}
      \sum_{n\geq0}\left(\frac{a_\epsilon}{1+b\delta_\epsilon}\right)^n
                    \mathcal K_\epsilon^n,
 \qquad
 \mathcal R_b^\chi:=\int_0^\infty e^{-br}\mathsf P_r^\chi\dd r.
\end{equation}
For $T<\infty$, let $\mathcal X_T$ consist of the functions in
$C_0(\overline E)$ vanishing for $t\geq T$.
For $q>0$ and $p=1,\infty$, write
$\|v\|_{q,p}:=\sup_{0\leq t\leq T}e^{qt}\|v(t,\cdot)\|_p$.
Fix $g\in C_c^+(\overline E)$ vanishing for $t\geq T$.
Then the following assertions hold.
\begin{enumerate}[label=(\roman*)]
\item \textit{Uniform estimates.}
The resolvents are bounded on $\mathcal X_T$, including at zero
killing, and
\begin{equation}\label{eq:clock-resolvent-bounds}
 \left\|\mathcal R_{\epsilon,b}\right\|_{q,p}
 \leq\frac{\delta_\epsilon}
     {1+b\delta_\epsilon-a_\epsilon\widehat H(q\tau_\epsilon)},
 \qquad
 \left\|\mathcal R_{\kappa+b}^\chi\right\|_{q,p}
 \leq\frac1{\kappa+b+\chi(q)}.
\end{equation}
These bounds also apply to measurable functions supported in
$0\leq t<T$ with
$\|v\|_{q,p}=\sup_{0\leq t\leq T}e^{qt}\|v(t,\cdot)\|_{L^p(\Rd)}<\infty$,
for $p=1$ or $p=\infty$.
  For all sufficiently small
$\epsilon$,
\begin{equation}\label{eq:clock-spatial-moment}
 \sup_{0\leq t\leq T}\int_{\Rd}\left(1+|x|^2\right)
           \mathcal R_{\epsilon,0}g(t,x)\dd x\leq C_{g,T}.
\end{equation}
\item \textit{Linear resolvent convergence.}
For every $b\geq0$,
\begin{equation}\label{eq:renewal-resolvent-strong}
 \left\|\mathcal R_{\epsilon,b}v-\mathcal R_{\kappa+b}^\chi v\right\|_\infty
       \xrightarrow[\epsilon\downarrow0]{}0\qquad(v\in\mathcal X_T).
\end{equation}
For the source $g$, convergence also holds in
$L^\infty([0,T];L^1(\Rd))$.
\item \textit{Nonlinear renewal limit.}
The cluster equation
\begin{equation}\label{eq:clock-cluster-renewal}
 \log\left(1-\delta_\epsilon u_\epsilon\right)
  =-\delta_\epsilon^2g-a_\epsilon\delta_\epsilon
                        \mathcal K_\epsilon u_\epsilon
\end{equation}
has a bounded nonnegative solution in $\mathcal X_T$, with
$0\leq u_\epsilon\leq\mathcal R_{\epsilon,0}g$.
It converges uniformly and in $L^\infty([0,T];L^1(\Rd))$ to the
unique bounded nonnegative solution in that space of
\begin{equation}\label{eq:clock-nonlinear-limit}
 u=\mathcal R_\kappa^\chi\left(g-u^2/2\right).
\end{equation}
\end{enumerate}
\end{lemma}

\begin{proof}
First, let us prove (i).
Fix $q>0$ and write $d_\epsilon=1-e^{-q\tau_\epsilon S}$.
Assumption \eqref{eq:clock-assumption} gives
$\E[d_\epsilon]=\mathcal O(\delta_\epsilon)$ and
$\E[d_\epsilon^2]\leq\E[d_\epsilon]$.
The geometric series proves \eqref{eq:clock-resolvent-bounds}.
To control spatial tails, let
\[
p_\epsilon:=\E\left[e^{-q\tau_\epsilon S}\right],
\quad
v_\epsilon:=\E\left[|Z|^2e^{-q\tau_\epsilon S}\right],
\quad\text{and}\quad
\mathbf m_\epsilon:=\E\left[Ze^{-q\tau_\epsilon S}\right].
\]
Centering and the Cauchy--Schwarz inequality imply
\[
 1-p_\epsilon\sim\chi(q)\delta_\epsilon,\qquad
 v_\epsilon\leq\operatorname{tr}\Sigma,\qquad
 |\mathbf m_\epsilon|^2
 \leq\operatorname{tr}\Sigma\,\E \left[d_\epsilon^2\right]=\mathcal O(\delta_\epsilon).
\]
For sums of $n$ independent increment pairs,
\[
 \E\left[\left|Z^{(n)}\right|^2e^{-q\tau_\epsilon S^{(n)}}\right]
  =n v_\epsilon p_\epsilon^{n-1}
   +n(n-1)|\mathbf m_\epsilon|^2p_\epsilon^{n-2},
\]
with the corresponding single terms for $n=0,1$.
A nonzero source term has $\tau_\epsilon S^{(n)}\leq T$.
Use $\mathbf1_{\{\tau_\epsilon S^{(n)}\leq T\}}
 \leq e^{qT-q\tau_\epsilon S^{(n)}}$ and translate the spatial
integral.  Summing the resulting geometric series bounds the
left-hand side of \eqref{eq:clock-spatial-moment} by
\begin{equation}\label{eq:clock-geometric-bound}
 C_g e^{qT}\left\{
 \frac{\delta_\epsilon}{1-p_\epsilon}
 +\frac{\delta_\epsilon^2v_\epsilon}{(1-p_\epsilon)^2}
 +\frac{2\delta_\epsilon^2|\mathbf m_\epsilon|^2}
                  {(1-p_\epsilon)^3}\right\}\leq C_{g,T}.
\end{equation}

Next, let us prove (ii).
With $d_\epsilon$ as in (i), split at $|Z|=L$ to obtain, for every $L>0$,
\begin{align}
 \frac{\left|\E[d_\epsilon
              (e^{\mathrm{i}\sqrt{\delta_\epsilon}\xi\cdot Z}-1)]\right|}
             {\delta_\epsilon}
             \leq C_{q,\xi}L\sqrt{\delta_\epsilon}
       +C_{q,\xi}\left(\E\left[|Z|^2\mathbf1_{\{|Z|>L\}}\right]\right)^{1/2}.
       \label{eq:clock-mixed-remainder}
\end{align}
Letting $\epsilon\downarrow0$ and then $L\to\infty$ makes the
mixed term vanish.  The centered spatial expansion therefore yields
\begin{equation}\label{eq:clock-joint-symbol}
 \E\left[e^{-q\tau_\epsilon S+\mathrm{i}\sqrt{\delta_\epsilon}\xi\cdot Z}\right]
  =1-\delta_\epsilon
       \left(\chi(q)+\frac12\xi^{\top}\Sigma\xi\right)
       +o(\delta_\epsilon).
\end{equation}
This uses no independence between the microscopic $S$ and $Z$.

For $b\geq0$ and $q>0$, let $Q_{\epsilon,b}$ and $Q_b$ be the
convolution measures of $\mathcal R_{\epsilon,b}$ and
$\mathcal R_{\kappa+b}^\chi$:
\begin{align*}
 Q_{\epsilon,b}
 &:=\frac{\delta_\epsilon}{1+b\delta_\epsilon}
   \sum_{n=0}^\infty\left(\frac{a_\epsilon}{1+b\delta_\epsilon}\right)^n
   \operatorname{Law}\left(\tau_\epsilon S^{(n)},
                              \sqrt{\delta_\epsilon}Z^{(n)}\right),\\
 Q_b&:=\int_0^\infty e^{-(\kappa+b)r}
                      \operatorname{Law}\left(D_r^\chi,B_r^\Sigma\right)\dd r.
\end{align*}
Their tilted versions $\widetilde Q_{\epsilon,b}(\dd s,\dd y)
:=e^{-qs}Q_{\epsilon,b}(\dd s,\dd y)$ and
$\widetilde Q_b(\dd s,\dd y):=e^{-qs}Q_b(\dd s,\dd y)$ are finite.
By \eqref{eq:clock-joint-symbol}, for $z\geq0$ and $\xi\in\Rd$,
\begin{align}
 \int_{\overline E}e^{-zs+\mathrm i\xi\cdot y}
                      \widetilde Q_{\epsilon,b}(\dd s,\dd y)
 &=\frac{\delta_\epsilon}{1+b\delta_\epsilon-a_\epsilon
       \E[e^{-(q+z)\tau_\epsilon S+\mathrm i\sqrt{\delta_\epsilon}\xi\cdot Z}]}
       \nonumber\\
 &\longrightarrow\frac1{\kappa+b+\chi(q+z)+\xi^\top\Sigma\xi/2}
 =\int_{\overline E}e^{-zs+\mathrm i\xi\cdot y}
                      \widetilde Q_b(\dd s,\dd y).
 \label{eq:clock-tilted-transform}
\end{align}
In particular their total masses converge.  With $p_\epsilon$,
$v_\epsilon$, and $\mathbf m_\epsilon$ from (i), and
$r_\epsilon=a_\epsilon/(1+b\delta_\epsilon)$, summing the spatial
second moments gives
\begin{align*}
 \int_{\overline E}|y|^2\widetilde Q_{\epsilon,b}(\dd s,\dd y)
=\frac{\delta_\epsilon^2}{1+b\delta_\epsilon}
 \left(\frac{r_\epsilon v_\epsilon}{(1-r_\epsilon p_\epsilon)^2}
       +\frac{2r_\epsilon^2|\mathbf m_\epsilon|^2}
                         {(1-r_\epsilon p_\epsilon)^3}\right)
 \leq C_q.
\end{align*}
Here $1-r_\epsilon p_\epsilon\geq1-p_\epsilon
\sim\chi(q)\delta_\epsilon$, $v_\epsilon\leq\operatorname{tr}\Sigma$,
and $|\mathbf m_\epsilon|^2=\mathcal O(\delta_\epsilon)$, so $C_q$
is independent of sufficiently small $\epsilon$.  Hence,
\begin{align}
 \widetilde Q_{\epsilon,b}\left([0,\infty)\times B_L^c\right)&\leq C_qL^{-2},
       \nonumber\\
 \widetilde Q_{\epsilon,b}\left((L,\infty)\times\Rd\right)
 &\leq\frac{1}{1-e^{-1}}
       \int_{\overline E}\left(1-e^{-s/L}\right)\widetilde Q_{\epsilon,b}(\dd s,\dd y),
       \nonumber\\
 \limsup_{\epsilon\downarrow0}
       \widetilde Q_{\epsilon,b}\left((L,\infty)\times\Rd\right)
 &\leq\frac1{1-e^{-1}}
 \left(\frac1{\kappa+b+\chi(q)}-
             \frac1{\kappa+b+\chi(q+L^{-1})}\right)
 \xrightarrow[L\to\infty]{}0.
 \label{eq:clock-tilted-tightness}
\end{align}
Thus the tilted measures are tight.  Every subsequential weak limit has
the transform \eqref{eq:clock-tilted-transform}, which determines a
finite measure, so $\widetilde Q_{\epsilon,b}\Rightarrow\widetilde Q_b$.

For $v\in\mathcal X_T$, let $v_q(t,x):=e^{qt}v(t,x)$ and
$F_{t,x}(s,y):=v_q(t+s,x+y)$.  These functions are uniformly bounded
and share the modulus of continuity of $v_q$.  On any compact set
$K\subset\overline E$, they have a finite uniform $\zeta$-net
$F_1,\ldots,F_N$ chosen from this family: a finite mesh in $K$ and
its common modulus of continuity provide such a net.  If
$\Delta_\epsilon=\widetilde Q_{\epsilon,b}-\widetilde Q_b$ and
$M_q$ bounds the sum of their total masses, then
\begin{align*}
 \sup_{0\leq t\leq T,\,x\in\Rd}
   \left|\int_{\overline E}F_{t,x}\,\dd\Delta_\epsilon\right|
   \leq\max_{1\leq j\leq N}
             \left|\int_{\overline E}F_j\,\dd\Delta_\epsilon\right|
       +\zeta M_q
       +2\|v_q\|_\infty
            \left(\widetilde Q_{\epsilon,b}+\widetilde Q_b\right)(K^c).
\end{align*}
The first term tends to zero by weak convergence.  Tightness makes the
last term arbitrarily small by enlarging $K$, and then $\zeta\downarrow0$
removes the middle term.  Since
\[
 e^{qt}\left(\mathcal R_{\epsilon,b}v-
              \mathcal R_{\kappa+b}^\chi v\right)(t,x)
       =\int_{\overline E}F_{t,x}\,\dd\Delta_\epsilon,
\]
this proves \eqref{eq:renewal-resolvent-strong}, including $b=\kappa=0$.

Let $v_\epsilon:=\mathcal R_{\epsilon,b}g$ and
$v:=\mathcal R_{\kappa+b}^\chi g$.  Positivity gives
$v_\epsilon\leq\mathcal R_{\epsilon,0}g$.  Fatou's lemma, uniform
convergence, and \eqref{eq:clock-spatial-moment} imply
\[
 \sup_{0\leq t\leq T}\int_{\Rd}\left(1+|x|^2\right)v(t,x)\dd x
 \leq C_{g,T}.
\]
Consequently, for $L>0$,
\begin{align}
 \sup_{0\leq t\leq T}\|v_\epsilon(t)-v(t)\|_1
 &\leq |B_L|\|v_\epsilon-v\|_\infty
       +\sup_{0\leq t\leq T}\int_{B_L^c}(v_\epsilon(t,x)+v(t,x))\dd x
       \nonumber\\
 &\leq |B_L|\|v_\epsilon-v\|_\infty+2C_{g,T}L^{-2}.
 \label{eq:clock-linear-L1-split}
\end{align}
Taking first $\epsilon\downarrow0$ and then $L\to\infty$ proves the
linear convergence in $L^\infty([0,T];L^1(\Rd))$.

Finally, let us prove (iii).
Define the scalar nonlinearity
\begin{equation}\label{eq:clock-Phi}
 \Phi_\epsilon(v):=
   \frac{-\log(1-\delta_\epsilon v)-\delta_\epsilon v}
             {\delta_\epsilon^2}.
\end{equation}
Choose $C>\sqrt{2\|g\|_\infty}$, $b\geq2C$, and
$\delta_\epsilon C\leq1/2$.
Then $0\leq\Phi_\epsilon'(v)=v/(1-\delta_\epsilon v)\leq b$
on $[0,C]$, and $\Phi_\epsilon(C)\geq C^2/2$.
Define the maps
\begin{align}
 \mathcal T_{\epsilon,b}v
 &:=\mathcal R_{\epsilon,b}\left(g+bv-\Phi_\epsilon(v)\right),\nonumber\\
 \mathcal T_b^\chi v
 &:=\mathcal R_{\kappa+b}^\chi\left(g+bv-v^2/2\right).
 \label{eq:renewal-fixed-point-maps}
\end{align}
For $0\leq r\leq C$, the functions $br-\Phi_\epsilon(r)$ and
$br-r^2/2$ are nondecreasing and vanish at zero, since their
derivatives are $b-\Phi_\epsilon'(r)\geq0$ and $b-r\geq0$.
Consequently, for $v\in\mathcal X_T$ with $0\leq v\leq C$,
\begin{align*}
 0\leq g+bv-\Phi_\epsilon(v)
 &\leq\|g\|_\infty+bC-\Phi_\epsilon(C)\leq bC,\\
 0\leq g+bv-v^2/2
 &\leq\|g\|_\infty+bC-C^2/2\leq bC,
\end{align*}
where the final inequalities are due to $\Phi_\epsilon(C)\geq C^2/2$
and $C^2/2>\|g\|_\infty$.  Both sources belong to
$\mathcal X_T$, which the resolvents preserve by part (i).
Moreover, $\mathcal K_\epsilon\mathbf1=\mathbf1$ and
$\mathsf P_r^\chi\mathbf1=\mathbf1$ for every $r\geq0$.
Since $b>0$, positivity and the convergent series and integral in
\eqref{eq:renewal-resolvents}, applied to
\eqref{eq:renewal-fixed-point-maps}, give
\begin{align*}
 0\leq\mathcal T_{\epsilon,b}v
 &\leq bC\mathcal R_{\epsilon,b}\mathbf1
 =\frac{bC\delta_\epsilon}{1+b\delta_\epsilon-a_\epsilon}
 =\frac{bC}{b+\kappa_\epsilon}\leq C,\\
 0\leq\mathcal T_b^\chi v
 &\leq bC\mathcal R_{\kappa+b}^\chi\mathbf1
 =bC\int_0^\infty e^{-(\kappa+b)r}\dd r
 =\frac{bC}{b+\kappa}\leq C.
\end{align*}
Thus, both maps preserve $\left\{v\in\mathcal X_T:0\leq v\leq C\right\}$.
By \eqref{eq:clock-resolvent-bounds}, the discrete Lipschitz constant
in $\|\cdot\|_{q,\infty}$ is bounded by
\begin{equation}\label{eq:renewal-weighted-contraction}
 \frac{b\delta_\epsilon}
 {1+b\delta_\epsilon-a_\epsilon\widehat H(q\tau_\epsilon)}
       \xrightarrow[\epsilon\downarrow0]{}\frac{b}{b+\kappa+\chi(q)}<1.
\end{equation}
For the limiting map $\mathcal T_b^\chi$ in
\eqref{eq:renewal-fixed-point-maps}, the derivative of
$r\mapsto br-r^2/2$ lies in $[0,b]$ on $[0,C]$.  Hence,
\[
 \left\|\mathcal T_b^\chi v-\mathcal T_b^\chi w\right\|_{q,\infty}
 \leq\frac{b}{b+\kappa+\chi(q)}\|v-w\|_{q,\infty},
 \qquad 0\leq v,w\leq C,
\]
with the same limiting contraction factor as in
\eqref{eq:renewal-weighted-contraction}.
Banach's fixed-point theorem gives the fixed points of the maps in
\eqref{eq:renewal-fixed-point-maps}.  If
$v_\epsilon=\mathcal T_{\epsilon,b}v_\epsilon$, the resolvent series
\eqref{eq:renewal-resolvents} and the definition \eqref{eq:clock-Phi} give
\begin{align*}
 &\left((1+b\delta_\epsilon)I-a_\epsilon\mathcal K_\epsilon\right)
                  v_\epsilon
 =\delta_\epsilon\left(g+bv_\epsilon-\Phi_\epsilon(v_\epsilon)\right),\\
 &\frac{v_\epsilon-a_\epsilon\mathcal K_\epsilon v_\epsilon}
      {\delta_\epsilon}+\Phi_\epsilon(v_\epsilon)=g,\\
 &\log(1-\delta_\epsilon v_\epsilon)
 =-\delta_\epsilon^2g-a_\epsilon\delta_\epsilon
                         \mathcal K_\epsilon v_\epsilon.
\end{align*}
Thus, this fixed point solves \eqref{eq:clock-cluster-renewal}.
The comparison argument of
Lemma~\ref{lem:scheme-comparison}, with $\delta_\epsilon$ replacing
$\epsilon^2$ and $a_\epsilon\mathcal K_\epsilon$ replacing
$K_\epsilon$, identifies it with the cluster cost $u_\epsilon$ in
\eqref{eq:clock-cluster-renewal}.
The semigroup property in \eqref{eq:renewal-resolvents} gives
\begin{align}
 b\mathcal R_\kappa^\chi\mathcal R_{\kappa+b}^\chi v
 &=\int_0^\infty\int_0^\infty
      b e^{-\kappa r-(\kappa+b)s}\mathsf P_{r+s}^\chi v\dd s\dd r
      \nonumber\\
 &=\int_0^\infty e^{-\kappa a}\left(1-e^{-ba}\right)
                                \mathsf P_a^\chi v\dd a
 =\left(\mathcal R_\kappa^\chi-\mathcal R_{\kappa+b}^\chi\right)v.
 \label{eq:clock-resolvent-identity}
\end{align}
For $v\in\mathcal X_T$, Fubini's theorem is justified by
\eqref{eq:clock-resolvent-bounds} applied to $|v|$, including
$\kappa=0$.  In particular,
$\left(I+b\mathcal R_\kappa^\chi\right)\mathcal R_{\kappa+b}^\chi
=\mathcal R_\kappa^\chi$.
Applying this identity to the fixed-point equation in
\eqref{eq:renewal-fixed-point-maps} yields
\[
 u+b\mathcal R_\kappa^\chi u
 =\mathcal R_\kappa^\chi\left(g-u^2/2+bu\right)
 =\mathcal R_\kappa^\chi\left(g-u^2/2\right)+b\mathcal R_\kappa^\chi u,
\]
which proves \eqref{eq:clock-nonlinear-limit}.  Conversely,
\eqref{eq:clock-resolvent-identity} and commutativity of the two
resolvents turn \eqref{eq:clock-nonlinear-limit} into the same
fixed-point equation.  Increasing $C$ gives uniqueness
among bounded nonnegative solutions in $\mathcal X_T$.
Uniform resolvent convergence, the common contraction constant, and
$\Phi_\epsilon(v)=v^2/2+\mathcal O(\delta_\epsilon)$ on $[0,C]$ prove uniform
convergence of the fixed points.

The inequality $1-e^{-Y}\leq Y$ gives
$u_\epsilon\leq\mathcal R_{\epsilon,0}g$.
Fatou's lemma transfers the spatial bound in (i) to $u$.
The same spatial decomposition as in (ii) proves the nonlinear
$L^1$ convergence.  It also gives continuity of the limiting
solutions from $[0,T]$ to $L^1(\Rd)$.
\end{proof}

For finite-mean delays, the clock exponent is $\chi(q)=m_1q$.
Specializing the preceding lemma therefore gives the backward
semilinear heat equation for the cluster cost.

\begin{lemma}[Finite-mean nonlinear renewal limit]\label{lem:nonlinear-renewal}
Under \eqref{eq:H1}--\eqref{eq:H4}, the cluster costs $u_\epsilon$
converge uniformly on $[0,T]\times\Rd$ to the unique bounded
nonnegative mild solution of
\begin{equation}\label{eq:log-Laplace-PDE}
 m_1\partial_tu+\frac12\Sigma:D_x^2u-\frac12u^2+f=0,
 \qquad u(T,\cdot)=0.
\end{equation}
For $0<a_\epsilon\leq1$ and
$(1-a_\epsilon)/\epsilon^2\to\kappa\geq0$, the corresponding costs
$u_\epsilon^a$ converge uniformly to $u_f^\kappa$, where
\begin{equation}\label{eq:shared-log-Laplace-mild}
 u_f^\kappa(t,\cdot)=\frac1{m_1}\int_t^T
    e^{-\kappa(s-t)/m_1}P_{s-t}^D
       \left[f(s,\cdot)-\left(u_f^\kappa(s,\cdot)\right)^2/2\right]\dd s.
\end{equation}
\end{lemma}

\begin{proof}
For $q>0$, \eqref{eq:H3} and the dominated convergence theorem give
\[
 \frac{1-\widehat H(q\epsilon^2)}{\epsilon^2}
 =\E\left[\frac{1-e^{-q\epsilon^2S}}{\epsilon^2}\right]
 \longrightarrow qm_1,
 \qquad 0\leq\frac{1-e^{-q\epsilon^2S}}{\epsilon^2}\leq qS.
\]
Thus Lemma~\ref{lem:renewal-framework} applies with
$\tau_\epsilon=\delta_\epsilon=\epsilon^2$ and $\chi(q)=m_1q$.
The limiting motion has $D_r^\chi=m_1r$ and spatial semigroup
$P_{m_1r}^D$.  Consequently \eqref{eq:clock-nonlinear-limit} becomes
\begin{align*}
 u_f^\kappa(t,\cdot)
 &=\int_0^\infty e^{-\kappa r}P_{m_1r}^D
       \left[f(t+m_1r,\cdot)-\left(u_f^\kappa(t+m_1r,\cdot)\right)^2/2\right]\dd r\\
 &=\frac1{m_1}\int_t^T e^{-\kappa(s-t)/m_1}P_{s-t}^D
       \left[f(s,\cdot)-\left(u_f^\kappa(s,\cdot)\right)^2/2\right]\dd s.
\end{align*}
This is \eqref{eq:shared-log-Laplace-mild}; setting $\kappa=0$ gives
the mild form of \eqref{eq:log-Laplace-PDE}.
If $u,v$ are bounded nonnegative solutions of
\eqref{eq:shared-log-Laplace-mild} and $U=\max\{\|u\|_\infty,
\|v\|_\infty\}$, semigroup contraction gives
\[
 \|u(t)-v(t)\|_\infty
 \leq\frac{U}{m_1}\int_t^T\|u(s)-v(s)\|_\infty\dd s.
\]
Gr\"onwall's inequality, with time reversed, gives $u=v$.

We verify the regularity used later in the exponential martingale.
Write $u:=u_f^\kappa$, $U:=\|u\|_\infty$, and
$P_r^{D,\kappa}:=e^{-\kappa r/m_1}P_r^D$.
For $1\leq i\leq d$, let $e_i$ be the $i$-th coordinate vector and,
for $h\ne0$, define
\begin{equation}\label{eq:log-Laplace-first-difference}
 D_i^h v(t,x):=\frac{v(t,x+he_i)-v(t,x)}h,
 \qquad A_i^h(t,x):=\frac{u(t,x+he_i)+u(t,x)}2.
\end{equation}
Spatial translations commute with $P_r^{D,\kappa}$ and
$D_i^h(u^2/2)=A_i^hD_i^h u$.  Applying $D_i^h$ to
\eqref{eq:shared-log-Laplace-mild} therefore gives the linear
Volterra equation
\begin{equation}\label{eq:log-Laplace-difference-Volterra}
 D_i^h u(t)=\frac1{m_1}\int_t^T P_{s-t}^{D,\kappa}
       \left[D_i^h f(s)-A_i^h(s)D_i^h u(s)\right]\dd s.
\end{equation}
Let $F_1:=\max_{1\leq i\leq d}\|\partial_i f\|_\infty$.
Since $0\leq A_i^h\leq U$ and $\|D_i^h f\|_\infty\leq F_1$,
semigroup contraction in \eqref{eq:log-Laplace-difference-Volterra}
and Gr\"onwall's inequality give
\begin{align}
 \left\|D_i^h u(t)\right\|_\infty
 &\leq\frac{(T-t)F_1}{m_1}
       +\frac U{m_1}\int_t^T\left\|D_i^h u(s)\right\|_\infty\dd s,\nonumber\\
 \sup_{h\ne0}\left\|D_i^h u\right\|_\infty
 &\leq L_1:=\frac{TF_1}{m_1}e^{UT/m_1}.
 \label{eq:log-Laplace-first-difference-bound}
\end{align}

Let $w_i$ be the bounded continuous solution of
\begin{equation}\label{eq:log-Laplace-first-Volterra}
 w_i(t)=\frac1{m_1}\int_t^T P_{s-t}^{D,\kappa}
                      \left[\partial_i f(s)-u(s)w_i(s)\right]\dd s.
\end{equation}
To construct this solution in $\mathcal X_T$, let
\[
 (\mathcal Vv)(t):=-\frac1{m_1}\int_t^T
 P_{s-t}^{D,\kappa}[u(s)v(s)]\dd s,
 \qquad
 a_i(t):=\frac1{m_1}\int_t^T P_{s-t}^{D,\kappa}\partial_i f(s)\dd s.
\]
The operator $\mathcal V$ preserves $\mathcal X_T$, and
$a_i\in\mathcal X_T$.
The ordered time integrals give
$\|\mathcal V^n\|\leq(UT/m_1)^n/n!$ for $n\geq0$.
Consequently $w_i=\sum_{n=0}^\infty\mathcal V^n a_i$
converges uniformly, and Gr\"onwall's inequality gives uniqueness
and $\|w_i\|_\infty\leq L_1$.
Since $u\in\mathcal X_T$ is uniformly continuous,
$\|A_i^h-u\|_\infty\to0$ as $h\to0$; smoothness of $f$ gives
$\|D_i^h f-\partial_i f\|_\infty\to0$.
Subtracting \eqref{eq:log-Laplace-first-Volterra} from
\eqref{eq:log-Laplace-difference-Volterra} and applying
Gr\"onwall's inequality yields
\begin{equation}\label{eq:log-Laplace-first-difference-limit}
 \left\|D_i^h u-w_i\right\|_\infty
 \leq\frac T{m_1}e^{UT/m_1}
  \left(\left\|D_i^h f-\partial_i f\right\|_\infty
             +\left\|A_i^h-u\right\|_\infty\|w_i\|_\infty\right)
 \longrightarrow0.
\end{equation}
Thus $\partial_i u=w_i$ is continuous and bounded.

For $1\leq j\leq d$, applying $D_j^h$ to
\eqref{eq:log-Laplace-first-Volterra}, with
$\tau_j^h u(t,x):=u(t,x+he_j)$, gives
\begin{equation}\label{eq:log-Laplace-second-difference-Volterra}
 D_j^h w_i(t)=\frac1{m_1}\int_t^T P_{s-t}^{D,\kappa}
 \left[D_j^h\partial_i f(s)
       -(\tau_j^h u(s))D_j^h w_i(s)-(D_j^h u(s))w_i(s)\right]\dd s.
\end{equation}
Let $F_2:=\max_{1\leq i,j\leq d}\|\partial_{ij}f\|_\infty$.
Equations \eqref{eq:log-Laplace-first-difference-bound} and
\eqref{eq:log-Laplace-second-difference-Volterra} imply
\[
 \sup_{h\ne0}\left\|D_j^h w_i\right\|_\infty
 \leq L_2:=\frac{T(F_2+L_1^2)}{m_1}e^{UT/m_1}.
\]
Let $w_{ij}$ solve
\begin{equation}\label{eq:log-Laplace-second-Volterra}
 w_{ij}(t)=\frac1{m_1}\int_t^T P_{s-t}^{D,\kappa}
       [\partial_{ij}f(s)-u(s)w_{ij}(s)-w_j(s)w_i(s)]\dd s.
\end{equation}
The convergent series for $\mathcal V$ constructs this solution,
and Gr\"onwall's inequality gives $\|w_{ij}\|_\infty\leq L_2$.
Subtracting \eqref{eq:log-Laplace-second-Volterra} from
\eqref{eq:log-Laplace-second-difference-Volterra} gives
\begin{align*}
 \left\|D_j^h w_i-w_{ij}\right\|_\infty
 \leq\frac T{m_1}e^{UT/m_1}\Big(&\left\|D_j^h\partial_i f-\partial_{ij}f\right\|_\infty
       +L_1\left\|D_j^h u-w_j\right\|_\infty\\
       &+L_2\left\|\tau_j^h u-u\right\|_\infty\Big)\longrightarrow0,
\end{align*}
by \eqref{eq:log-Laplace-first-difference-limit} and uniform continuity.
It follows that $\partial_{ij}u=w_{ij}$ is continuous and bounded.
The two limiting equations are therefore
\begin{align}
 \partial_i u(t)
 &=\frac1{m_1}\int_t^T P_{s-t}^{D,\kappa}
          [\partial_i f(s)-u(s)\partial_i u(s)]\dd s,\nonumber\\
 \partial_{ij}u(t)
 &=\frac1{m_1}\int_t^T P_{s-t}^{D,\kappa}
          [\partial_{ij}f(s)-u(s)\partial_{ij}u(s)
                       -\partial_i u(s)\partial_j u(s)]\dd s,
 \label{eq:log-Laplace-spatial-derivatives}
\end{align}
with $\|\partial_i u\|_\infty\leq L_1$ and
$\|\partial_{ij}u\|_\infty\leq L_2$.
Differentiating \eqref{eq:shared-log-Laplace-mild} in time, using
these continuous bounded spatial derivatives, yields
\[
 \partial_tu=-\cA u+\frac{\kappa}{m_1}u
                  +\frac{u^2}{2m_1}-\frac{f}{m_1}.
\]
Thus, $u\in C^{1,2}$ with bounded derivatives.  The source vanishes on a
neighborhood of $T$, and so does $u$, which also justifies its zero
extension beyond $T$. This completes the proof.
\end{proof}

Integrating cluster costs over a Poisson field on the entire negative-time
half-space requires control of their spatial tails. The bound in the next lemma gives the
$L^1$ convergence needed for that integration.

\begin{lemma}[Finite-mean spatial localization]\label{lem:spatial-localization}
With the notation of Lemma~\ref{lem:nonlinear-renewal}, let
\begin{equation}\label{eq:linear-cluster-majorant}
 b_\epsilon(t,x):=\epsilon^2\sum_{n\geq0}
    \E\left[f\left(t+\epsilon^2S^{(n)},x+\epsilon Z^{(n)}\right)\right].
\end{equation}
There is $C_{f,T}<\infty$, uniform for sufficiently small $\epsilon$,
such that
\begin{equation}\label{eq:u-weighted-tail}
 0\leq u_\epsilon^a\leq u_\epsilon\leq b_\epsilon,
 \qquad
 \sup_{0\leq t\leq T}\int_{\Rd}\left(1+|x|^2\right)b_\epsilon(t,x)\dd x\leq C_{f,T},
\end{equation}
and
\begin{equation}\label{eq:u-L1-convergence}
 \sup_{0\leq t\leq T}\left\|u_\epsilon^a(t,\cdot)-u_f^\kappa(t,\cdot)\right\|_1\to0.
\end{equation}
\end{lemma}

\begin{proof}
The bound for $b_\epsilon$ is \eqref{eq:clock-spatial-moment} with
$a_\epsilon=1$.  Thinning gives the bound for $u_\epsilon^a$,
and the last assertion is the $L^1$ conclusion of
Lemma~\ref{lem:renewal-framework}.
\end{proof}

\subsubsection{The negative-time Poisson field}

For a past parent $(r,y)$, $r<0$, its positive-time children initiate
independent clusters.  Define
\begin{equation}\label{eq:J-epsilon}
  J_\epsilon(r,y)
  :=\int_0^\infty\int_{\Rd}
    h(q-r,w-y)u_\epsilon\left(\epsilon^2q,\epsilon w\right)
    \dd w\dd q.
\end{equation}
The Laplace transform contributed by this parent is
$\exp\{-\epsilon^2J_\epsilon(r,y)\}$.  Applying the exponential formula
once more, now to the varying past Poisson field
\eqref{eq:varying-past-intensity}, gives
\begin{equation}\label{eq:full-Laplace}
  \log\E\left[e^{-\ip{\cN_\epsilon^{\ell_\epsilon}}{f}}\right]
  =\epsilon^{d-2}\int_{-\infty}^0\int_{\Rd}
    \ell_\epsilon(\epsilon y)
    \left(e^{-\epsilon^2J_\epsilon(r,y)}-1\right)
    \dd y\dd r.
\end{equation}

We now combine the cluster limit with the Laplace formula for the Poisson past.  The finite mass $m_1$ of the kernel crossing time zero
collapses the past contribution to the initial pairing below.

\begin{lemma}[Collapse of the past layer]\label{lem:Laplace-boundary}
With $u$ from Lemma~\ref{lem:nonlinear-renewal},
\begin{equation}\label{eq:Laplace-boundary-limit}
  \log\E\left[e^{-\ip{\cN_\epsilon^{\ell_\epsilon}}{f}}\right]
  \xrightarrow[\epsilon\downarrow0]{}
  -m_1\int_{\Rd}\ell(x)u(0,x)\dd x.
\end{equation}
\end{lemma}

\begin{proof}
The bound $0\leq J_\epsilon\leq U_f$ and Tonelli's theorem give
\[
 \int_{r<0}\int_{\Rd}J_\epsilon(r,y)\dd y\dd r
 \leq\epsilon^{-d}m_1
          \sup_{0\leq t\leq T}\|u_\epsilon(t,\cdot)\|_1\leq C_{f,T}\epsilon^{-d}.
\]
Thus, the Poisson formula \eqref{eq:full-Laplace} follows first for bounded parent sets and
then for the entire past by monotone exhaustion.
Since
$0\leq e^{-a}-1+a\leq a^2/2$ and $J_\epsilon^2\leq U_fJ_\epsilon$,
the error in linearizing \eqref{eq:full-Laplace} is at most
\begin{equation}\label{eq:Laplace-linearization-error}
 \frac{\epsilon^{d+2}}2\|\ell_\epsilon\|_\infty
       \int_{r<0}\int_{\Rd}\left(J_\epsilon(r,y)\right)^2\dd y\dd r
 \leq C\epsilon^2\|\ell_\epsilon\|_\infty.
\end{equation}
After setting $x=\epsilon y$ and $z=w-y$, the leading term is
\begin{equation}\label{eq:boundary-leading-varying}
 -\int_{\Rd}\ell_\epsilon(x)A_\epsilon(x)\dd x,\qquad
 A_\epsilon(x):=\int_{r<0}\int_{q>0}\int_{\Rd}
 h(q-r,z)u_\epsilon\left(\epsilon^2q,x+\epsilon z\right)\dd z\dd q\dd r.
\end{equation}
The crossing measure $h(q-r,z)\dd r\dd q\dd z$ on $r<0<q$ has
total mass $m_1$.  For each fixed $(q,z)$, Lemma
\ref{lem:spatial-localization} and continuity in $L^1$ imply
\[
\left\|u_\epsilon\left(\epsilon^2q,\cdot+\epsilon z\right)-u(0,\cdot)\right\|_1\to0.
\]
The uniform $L^1$ bound permits the dominated convergence theorem against the
finite crossing measure.  Hence
$\|A_\epsilon-m_1u(0,\cdot)\|_1\to0$.
Bounded locally uniform convergence of $\ell_\epsilon$ and
$u(0,\cdot)\in L^1$ now give the claimed limit.
\end{proof}

\subsubsection{Identification with super-Brownian occupation measure}

\begin{lemma}[Construction from bounded densities of infinite total mass]
\label{lem:locally-finite-SBM}
Let $\ell:\Rd\to\R_+$ be bounded and Borel measurable, and let
$\mu(\dd x)=m_1\ell(x)\dd x$.  There is a canonical process $X^\ell$
with continuous paths in $\cM_{\mathrm{loc}}(\Rd)$, initial measure $\mu$,
and martingale characteristics
\eqref{eq:SBM-martingale}--\eqref{eq:SBM-bracket}.  Its construction by
independent finite-measure components does not depend in law on the
chosen exhaustion.  For every nonnegative compactly supported Borel
function $\varphi$,
\begin{equation}\label{eq:SBM-local-mean}
  \E\left[\ip{X_t^\ell}{\varphi}\right]
  =\ip{\mu}{P_t^D\varphi}
  \leq m_1\norm{\ell}_\infty\|\varphi\|_{L^1(\Rd)}.
\end{equation}
Its occupation measure $\nu_X=m_1^{-1}X_t^\ell(\dd x)\dd t$ is locally
finite, and for $f\in C_c^\infty([0,T)\times\Rd)$, $f\geq0$,
\begin{equation}\label{eq:SBM-Laplace}
  \E\left[\exp\!\left\{-\frac1{m_1}\int_0^T
       \ip{X_t^\ell}{f(t,\cdot)}\dd t\right\}\right]
  =\exp\!\left\{-\ip{X_0^\ell}{u(0,\cdot)}\right\},
\end{equation}
where $u$ solves \eqref{eq:log-Laplace-PDE}.  In dimension one this
construction has a measurable density and the weak white-noise
representation \eqref{eq:candidate-spde}, with initial profile $\ell$.
\end{lemma}

\begin{proof}
For a finite initial measure $\mu_0$, apply
\cite[Theorem~II.5.1(a)--(b)]{perkins2002superprocess} with motion
generator $\cA$, branching variance $1/m_1$, and zero linear drift.
Brownian motion is a continuous strong Markov process whose semigroup
preserves bounded continuous functions, and both branching coefficients
are bounded and continuous, as required there.  For independent
solutions $Z^{(1)},\ldots,Z^{(n)}$ with initial measures
$\mu_1,\ldots,\mu_n$, write $Y=\sum_{j=1}^nZ^{(j)}$.
Their martingales are orthogonal, and hence
\begin{align*}
 M_t^Y(\psi)&=\sum_{j=1}^nM_t^{(j)}(\psi),\\
 \left\langle M^Y(\psi)\right\rangle_t
 &=\sum_{j=1}^n\left\langle M^{(j)}(\psi)\right\rangle_t
 =\frac1{m_1}\int_0^t\ip{Y_s}{\psi^2}\dd s.
\end{align*}
Thus $Y$ solves \eqref{eq:SBM-martingale}--\eqref{eq:SBM-bracket}
with initial measure $\sum_{j=1}^n\mu_j$, and uniqueness gives its law.
For $\varphi\in C_c^\infty(\Rd)$ and fixed $t>0$, let
$\varphi_s=P_{t-s}^D\varphi$, $0\leq s\leq t$.  This test function has bounded
space and time derivatives and satisfies
$\partial_s\varphi_s+\cA\varphi_s=0$.
The time-dependent formula
\cite[Proposition~II.5.7]{perkins2002superprocess} gives
\begin{equation}\label{eq:SBM-finite-mean-identity}
 \ip{X_t}{\varphi}
 =\ip{\mu_0}{P_t^D\varphi}+M_t(\varphi_\cdot),\qquad
 \E\left[\ip{X_t}{\varphi}\right]=\ip{\mu_0}{P_t^D\varphi}.
\end{equation}
Indeed, localization of the nonnegative total-mass martingale yields
$\E[X_s(\Rd)]\leq\mu_0(\Rd)$, and therefore,
\[
 \E\left[\langle M(\varphi_\cdot)\rangle_t\right]
 =\frac1{m_1}\int_0^t\E\left[\ip{X_s}{(P_{t-s}^D\varphi)^2}\right]\dd s
 \leq\frac{t}{m_1}\|\varphi\|_\infty^2\mu_0(\Rd)<\infty.
\]
Thus the martingale in \eqref{eq:SBM-finite-mean-identity} has zero
mean.  Equality on smooth compactly supported test functions identifies the
finite measures $\E[X_t]$ and $\mu_0P_t^D$; integration against any
nonnegative Borel function then gives the same identity, with
monotone convergence allowing unbounded test functions.

Let $A_1=B_1$ and
$A_j=B_j\setminus B_{j-1}$ for $j\geq2$.
Take independent finite-measure processes $Z^{(j)}$ started at
$\mu|_{A_j}$ and write
$X^{[n]}=\sum_{j=1}^n Z^{(j)}$.
These partial sums increase as measures and have the finite-measure
Dawson--Watanabe laws started at $\mu|_{B_n}$.
To control the sum simultaneously at all times, fix $p>d/2$ and let
$w(x)=(1+|x|^2)^{-p}$.  This weight is integrable and satisfies
$|\cA w|\leq Cw$, so $\ip{\mu}{w}<\infty$.
Choose $c\geq0$ with $\cA w\leq cw$.
For $Y^{n,m}=X^{[m]}-X^{[n]}$, $m>n$, the process
$e^{-ct}\ip{Y_t^{n,m}}{w}$ is a nonnegative supermartingale: its
finite-measure martingale decomposition has nonpositive drift.
The maximal inequality therefore gives
\begin{equation}\label{eq:SBM-weighted-tail}
  \mathbb{P}\!\left(\sup_{0\leq t\leq T}\left\langle Y_t^{n,m},w\right\rangle>\delta_0\right)
  \leq \frac{e^{cT}}{\delta_0}
       \ip{\mu}{w\mathbf 1_{\{|x|>n\}}}.
\end{equation}
The right-hand side of \eqref{eq:SBM-weighted-tail} tends to zero as $n\to\infty$, uniformly in $m>n$.
Passing $m\to\infty$ in the increasing events and choosing $n_k$ so that
$\ip{\mu}{w\mathbf 1_{\{|x|>n_k\}}}\leq2^{-3k}$, the
Borel--Cantelli lemma with $\delta_0=2^{-k}$ proves almost sure convergence
of the series uniformly in weighted mass on $[0,T]$.
Intersecting over integer $T$ gives a process
$X_t^\ell=\sum_{j\geq1}Z_t^{(j)}$ that is locally finite at every time.
For each $\psi\in C_c(\Rd)$, the bound $|\psi|\leq C_\psi w$ makes
$\ip{X_t^{[n]}}{\psi}$ converge uniformly on compact time intervals.
Thus the limit has continuous paths in the vague topology.

Applying \eqref{eq:SBM-finite-mean-identity} to $X^{[n]}$ and using
monotone convergence gives, for nonnegative compactly supported Borel
$\varphi$,
\[
 \E\left[\ip{X_t^\ell}{\varphi}\right]
 =\lim_{n\to\infty}\ip{\mu|_{B_n}}{P_t^D\varphi}
 =\ip{\mu}{P_t^D\varphi}
 \leq m_1\|\ell\|_\infty\|\varphi\|_1.
\]
This proves \eqref{eq:SBM-local-mean}, and the same identity with
$\varphi=w$, together with $P_t^Dw\leq e^{ct}w$, gives
$\E\left[\ip{X_t^\ell}{w}\right]\leq e^{ct}\ip{\mu}{w}$.
For $\psi\in C_c^2(\Rd)$, the martingale of a tail $Y^{n,m}$ has
second moment at time $T$ bounded by
\[
  \frac1{m_1}\int_0^T
       \ip{\mu|_{\{|x|>n\}}}{P_s^D\psi^2}\dd s
  \leq C_{\psi,T}\ip{\mu}{w\mathbf 1_{\{|x|>n\}}}.
\]
Here $\psi^2\leq C_\psi w$ and $P_s^Dw\leq e^{cs}w$.
The same weighted bound controls the drift tails with test function $\cA\psi$.
Consequently the martingales converge in $L^2$, uniformly in time by
Doob's inequality, and their brackets converge in $L^1$.
Passing to the limit gives precisely
\eqref{eq:SBM-martingale}--\eqref{eq:SBM-bracket}.

We next derive the occupation formula with its normalization.
For the finite process $X^{[n]}$, define
\[
 A_t^{[n]}:=\ip{X_t^{[n]}}{u(t,\cdot)}
       +\frac1{m_1}\int_0^t\ip{X_s^{[n]}}{f(s,\cdot)}\dd s,
 \qquad Z_t^{[n]}:=e^{-A_t^{[n]}}.
\]
The regularity in \eqref{eq:log-Laplace-spatial-derivatives} permits the
time-dependent test formula
\cite[Proposition~II.5.7]{perkins2002superprocess}.  Its martingale
part $L^{[n]}$ and \eqref{eq:log-Laplace-PDE} give
\begin{align}
 \dd A_t^{[n]}
 &=\ip{X_t^{[n]}}{\partial_tu+\cA u+f/m_1}\dd t
                         +\dd L_t^{[n]}
 =\frac1{2m_1}\ip{X_t^{[n]}}{(u(t,\cdot))^2}\dd t+\dd L_t^{[n]},
       \nonumber\\
 \dd\left\langle L^{[n]}\right\rangle_t
 &=\frac1{m_1}\ip{X_t^{[n]}}{(u(t,\cdot))^2}\dd t,\nonumber\\
 \dd Z_t^{[n]}
 &=-Z_t^{[n]}\dd L_t^{[n]}
   +Z_t^{[n]}\left[-\frac1{2m_1}\ip{X_t^{[n]}}{(u(t))^2}
                    +\frac1{2m_1}\ip{X_t^{[n]}}{(u(t))^2}\right]\dd t
 =-Z_t^{[n]}\dd L_t^{[n]}.
 \label{eq:SBM-occupation-exponential}
\end{align}
The last identity is It\^o's formula; it exhibits the cancellation of
the drift and half the quadratic variation.  Since $u,f\geq0$,
$0\leq Z_t^{[n]}\leq1$, so the local martingale is a true martingale.
Using $u(T,\cdot)=0$ in $\E\left[Z_T^{[n]}\right]=Z_0^{[n]}$ gives
\[
 \E\left[\exp\!\left\{-\frac1{m_1}\int_0^T
                    \ip{X_s^{[n]}}{f(s,\cdot)}\dd s\right\}\right]
 =\exp\!\left\{-\ip{\mu|_{B_n}}{u(0,\cdot)}\right\}.
\]

The linear heat equation supersolution gives
\begin{equation}\label{eq:SBM-Laplace-integrability}
  0\leq u(t,\cdot)\leq
     \frac1{m_1}\int_t^T P_{s-t}^D f(s,\cdot)\dd s,
  \qquad
  \|u(0,\cdot)\|_1\leq
     \frac1{m_1}\int_0^T\|f(s,\cdot)\|_1\dd s.
\end{equation}
In particular, $\ip{\mu}{u(0,\cdot)}<\infty$.
The occupation integrals increase with $n$, and bounded convergence of
their exponentials, together with monotone convergence of
$\ip{\mu|_{B_n}}{u(0,\cdot)}$, proves
\eqref{eq:SBM-Laplace} for $X^\ell$.
Moreover, for every bounded spatial Borel set $K$,
\begin{equation}\label{eq:SBM-occupation-mean}
  \E\left[\nu_X((0,T]\times K)\right]
  \leq T\norm{\ell}_\infty |K|.
\end{equation}
This proves local finiteness of the occupation measure, including control
of its mass as time decreases to zero.

To verify independence of the exhaustion for the whole process, let
$0<t_1<\cdots<t_k$ and $\varphi_i\in C_c^+(\Rd)$ for $1\leq i\leq k$.
Let $\mathcal U_t\varphi$ be the finite-measure transition log-Laplace
solution in \cite[Theorem~II.5.11(c)]{perkins2002superprocess} with zero
occupation source.  Its positivity and linear majorant give
$0\leq\mathcal U_t\varphi\leq P_t^D\varphi$.
The Markov property yields the recursion
\begin{equation}\label{eq:SBM-multitime-recursion}
 V_k:=\varphi_k,\qquad
 V_i:=\varphi_i+\mathcal U_{t_{i+1}-t_i}V_{i+1},\quad (1\leq i<k),
 \qquad V:=\mathcal U_{t_1}V_1.
\end{equation}
Iterating the transition formula therefore gives
\begin{equation}\label{eq:SBM-multitime-Laplace}
 \E\left[\exp\!\left\{-\sum_{i=1}^k\ip{X_{t_i}^{[n]}}{\varphi_i}\right\}\right]
 =e^{-\ip{\mu|_{B_n}}{V}},
 \qquad 0\leq V\leq\sum_{i=1}^kP_{t_i}^D\varphi_i.
\end{equation}
The last bound follows inductively from
\eqref{eq:SBM-multitime-recursion} and the semigroup property.  It is
integrable against $\mu$ because
\[
 \ip{\mu}{V}\leq m_1\|\ell\|_\infty
                        \sum_{i=1}^k\|\varphi_i\|_1<\infty.
\]
Let $\mu_n\uparrow\mu$ be any increasing finite-measure
restrictions, and construct $X^{\mu_n}$ by independent increments in
the initial measure.  The bound in \eqref{eq:SBM-multitime-Laplace}
and monotone convergence theorem give
\begin{align}
 \ip{\mu-\mu_n}{V}
 &\leq\sum_{i=1}^k\ip{\mu-\mu_n}{P_{t_i}^D\varphi_i}
 \longrightarrow0,\nonumber\\
 \E\left[\exp\left\{-\sum_{i=1}^k
                     \ip{X_{t_i}^{\mu}}{\varphi_i}\right\}\right]
 &=\lim_{n\to\infty}e^{-\ip{\mu_n}{V}}
 =e^{-\ip{\mu}{V}}.
 \label{eq:SBM-exhaustion-Laplace}
\end{align}
Here the process limit exists because \eqref{eq:SBM-weighted-tail}
applied to $\mu_m-\mu_n$ yields
\[
 \sup_{m>n}\mathbb P\left(\sup_{0\leq t\leq T}
        \ip{X_t^{\mu_m}-X_t^{\mu_n}}w>\delta_0\right)
 \leq\frac{e^{cT}}{\delta_0}\ip{\mu-\mu_n}w\longrightarrow0;
\]
a summable subsequence gives the same uniform convergence and path
continuity as above.  Equation \eqref{eq:SBM-exhaustion-Laplace}
does not involve the restrictions.  Replacing each $\varphi_i$ by
$a_i\varphi_i$, $a_i\geq0$, identifies all joint laws of these
pairings by uniqueness of multivariate Laplace transforms.
Choose a countable family $(\psi_j)_{j\geq1}\subset C_c^+(\Rd)$
that determines the vague topology.  On the continuous-path space,
\[
 \mathcal B\left(C\left([0,\infty),\cM_{\mathrm{loc}}(\Rd)\right)\right)
 =\sigma\left(X\mapsto\ip{X_q}{\psi_j}:
                    q\in\mathbb Q_+,\ j\geq1\right).
\]
To see this, continuity recovers every real-time evaluation from
rational times, and uniform distances on compact time intervals are
suprema over rational times.  Equality of the joint laws of these
coordinates proves independence of the exhaustion for the process
law.

Assume now $d=1$.  By
\cite[Theorem~III.4.2(a)--(b)]{perkins2002superprocess}, with the
Brownian and branching coefficients rescaled as in
\eqref{eq:SBM-martingale}--\eqref{eq:SBM-bracket}, each finite component
has a predictable measurable density $v_j$ and a white noise $W_j$.
Take these pairs independently.  For $j\geq1$ and $\psi\in C_c^2(\R)$,
the finite-component version of \eqref{eq:SBM-martingale} is
\begin{align}
 \ip{v_j(t,\cdot)}{\psi}
 &=\ip{\mu|_{A_j}}{\psi}
   +\frac{m_2}{2m_1}\int_0^t\ip{v_j(s,\cdot)}{\psi''}\dd s\nonumber\\
 &\qquad+\int_0^t\int_{\R}\psi(x)\sqrt{v_j(s,x)/m_1}\,
                                      W_j(\dd s,\dd x).
 \label{eq:SBM-finite-component-density}
\end{align}
Define $v:=\sum_{j=1}^\infty v_j$.  For each fixed $t>0$, monotone convergence
of the measures gives $X_t^\ell(\dd x)=v(t,x)\dd x$.  Moreover,
\begin{equation}\label{eq:SBM-density-local-integrability}
 \E\left[\int_0^T\int_K v(t,x)\dd x\dd t\right]
 =\int_0^T\E\left[X_t^\ell(K)\right]\dd t
 \leq m_1\|\ell\|_\infty T|K|<\infty
\end{equation}
for every bounded Borel set $K\subset\R$.  Thus $v$ is finite almost
everywhere in space--time, almost surely.  For a deterministic bounded
test function $a$ supported in $[0,T]\times K$, independence and the white-noise
isometry give
\begin{align}
 &\E\left[\left|\sum_{j=n+1}^m\int_0^T\int_{\R}
             a(t,x)\sqrt{v_j(t,x)/m_1}\,W_j(\dd t,\dd x)\right|^2\right]
       \nonumber\\
 &\qquad=\frac1{m_1}\E\left[\int_0^T\int_{\R}(a(t,x))^2
                       \sum_{j=n+1}^m v_j(t,x)\dd x\dd t\right]
 \longrightarrow0.
 \label{eq:SBM-noise-sum}
\end{align}
The limit is uniform in $m>n$, by
\eqref{eq:SBM-density-local-integrability} and dominated convergence theorem.

Let $W_0$ be an additional independent white noise.  For $j\geq1$,
let $\theta_j:=\sqrt{v_j/v}$ on $\{0<v<\infty\}$ and $\theta_j:=0$
elsewhere, and let $\theta_0:=\mathbf1_{\{v=0\}\cup\{v=\infty\}}$.
The set $\{v=\infty\}$ has zero space--time measure almost surely by
\eqref{eq:SBM-density-local-integrability}.  The predictable coefficients
satisfy $\sum_{j=0}^\infty\theta_j^2=1$ almost everywhere.  Define
\begin{equation}\label{eq:SBM-white-noise-construction}
 W(\dd t,\dd x):=\theta_0W_0(\dd t,\dd x)
                  +\sum_{j=1}^\infty\theta_jW_j(\dd t,\dd x).
\end{equation}
The series converges in $L^2$ against every deterministic compactly
supported square-integrable test function.  For spatial test functions $a,b\in L^2(\R)$,
its continuous martingales have covariations
\[
 \langle W(a),W(b)\rangle_t
 =\int_0^t\int_{\R}a(x)b(x)
                    \sum_{j=0}^\infty(\theta_j(s,x))^2\dd x\dd s
 =t\int_{\R}a(x)b(x)\dd x.
\]
The multidimensional L\'evy characterization applied to each finite
collection of test functions shows that $W$ is a space--time white noise.
Equations \eqref{eq:SBM-noise-sum} and
\eqref{eq:SBM-white-noise-construction} identify the summed martingale
measure with $\sqrt{v/m_1}\,W$.
For $\psi\in C_c^2(\R)$, choose $C_\psi<\infty$ so that
$|\psi|+|\psi''|\leq C_\psi w$.  Letting $m\to\infty$ in
\eqref{eq:SBM-weighted-tail} gives, for every $\delta_0>0$,
\begin{equation}\label{eq:SBM-infinite-weighted-tail}
 \mathbb P\!\left(\sup_{0\leq t\leq T}
       \ip{X_t^\ell-X_t^{[n]}}{w}>\delta_0\right)
 \leq\frac{e^{cT}}{\delta_0}
       \ip{\mu}{w\mathbf1_{\{|x|>n\}}}\longrightarrow0.
\end{equation}
Consequently,
\begin{align}
 &\sup_{0\leq t\leq T}\left|\ip{X_t^\ell-X_t^{[n]}}{\psi}\right|
 +\sup_{0\leq t\leq T}\left|\int_0^t
                 \ip{X_s^\ell-X_s^{[n]}}{\psi''}\dd s\right|\nonumber\\
 &\qquad\leq C_\psi(1+T)\sup_{0\leq s\leq T}
                       \ip{X_s^\ell-X_s^{[n]}}{w}
       \xrightarrow[n\to\infty]{\mathbb P}0.
 \label{eq:SBM-density-tail-limits}
\end{align}
For $0\leq t\leq T$, let
\begin{align*}
 M_t^{[n]}(\psi)
 &:=\sum_{j=1}^n\int_0^t\int_{\R}
       \psi(x)\sqrt{v_j(s,x)/m_1}\,W_j(\dd s,\dd x),\\
 M_t(\psi)
 &:=\int_0^t\int_{\R}\psi(x)\sqrt{v(s,x)/m_1}\,W(\dd s,\dd x).
\end{align*}
Summing \eqref{eq:SBM-finite-component-density} over $1\leq j\leq n$
gives
\begin{equation}\label{eq:SBM-summed-component-density}
 \ip{X_t^{[n]}}{\psi}
 =\ip{\mu|_{B_n}}{\psi}
  +\frac{m_2}{2m_1}\int_0^t\ip{X_s^{[n]}}{\psi''}\dd s
  +M_t^{[n]}(\psi).
\end{equation}
By \eqref{eq:SBM-density-tail-limits}, the left-hand side of
\eqref{eq:SBM-summed-component-density} converges to
$\ip{X_t^\ell}{\psi}$, and its drift term converges to
$\frac{m_2}{2m_1}\int_0^t\ip{X_s^\ell}{\psi''}\dd s$,
both in probability uniformly on $[0,T]$.
For the martingale term $M^{[n]}(\psi)$ in
\eqref{eq:SBM-summed-component-density},
\eqref{eq:SBM-noise-sum} and \eqref{eq:SBM-white-noise-construction},
together with Doob's maximal inequality, give
\begin{align*}
 \E\left[\sup_{0\leq t\leq T}
          \left|M_t(\psi)-M_t^{[n]}(\psi)\right|^2\right]
 &\leq4\E\left[\left|M_T(\psi)-M_T^{[n]}(\psi)\right|^2\right]\\
 &=\frac4{m_1}\E\left[\int_0^T\int_{\R}(\psi(x))^2
             \sum_{j=n+1}^\infty v_j(s,x)\dd x\dd s\right]
 \longrightarrow0.
\end{align*}
The limit follows from \eqref{eq:SBM-density-local-integrability}
and the dominated convergence theorem, since $\psi$ is bounded and
compactly supported and $\sum_{j=n+1}^\infty v_j\leq v$.
Since $\ip{\mu|_{B_n}}{\psi}\to m_1\ip{\ell}{\psi}$, passing to the
limit in \eqref{eq:SBM-summed-component-density} yields
\[
 \ip{v_t}{\psi}=m_1\ip{\ell}{\psi}
 +\frac{m_2}{2m_1}\int_0^t\ip{v_s}{\psi''}\dd s
 +\int_0^t\int_{\R}\psi(x)\sqrt{v(s,x)/m_1}\,W(\dd s,\dd x).
\]
Setting $\rho=v/m_1$ yields
\[
 m_1\ip{\rho_t}{\psi}=m_1\ip{\ell}{\psi}
 +\frac{m_2}{2}\int_0^t\ip{\rho_s}{\psi''}\dd s
 +\int_0^t\int_{\R}\psi(x)\sqrt{\rho(s,x)}\,W(\dd s,\dd x),
\]
which is \eqref{eq:candidate-spde}.  The initial value follows from
the vague continuity of $X^\ell$ at zero.
\end{proof}

Because $X_0^\ell=m_1\ell(x)\dd x$, the right-hand side of
\eqref{eq:SBM-Laplace} is exactly the limit in
\eqref{eq:Laplace-boundary-limit}.

\begin{proof}[Proof of Theorems~\ref{thm:critical} and
\ref{thm:critical-varying}]
Let $C_*:=\sup_{0<\epsilon\leq1}\|\ell_\epsilon\|_\infty$.
The construction and \eqref{eq:measure-mean-bound} give
\begin{equation}\label{eq:varying-mean-bound}
 \E\left[\ip{\cN_\epsilon^{\ell_\epsilon}}{a}\right]\leq C_*\|a\|_{L^1(E)}
 \qquad(a\geq0\text{ compactly supported Borel}).
\end{equation}
Proposition~\ref{prop:tightness} therefore applies uniformly.  Let
$\cN_{\epsilon_n}^{\ell_{\epsilon_n}}\Rightarrow\nu$ along a subsequence.
For $f\in C_c^\infty(E)$ with $f\geq0$, the map
$\xi\mapsto e^{-\ip{\xi}{f}}$ is bounded and vaguely continuous.
Equations \eqref{eq:Laplace-boundary-limit} and \eqref{eq:SBM-Laplace}
thus imply
\begin{equation}\label{eq:critical-subsequence-Laplace}
 \E\left[ e^{-\ip{\nu}{f}}\right]
 =\exp\!\left\{-m_1\int_{\Rd}\ell(x)u_f(0,x)\dd x\right\}
 =\E\left[e^{-\ip{\nu_X}{f}}\right],
 \qquad \nu_X:=m_1^{-1}X_t^\ell(\dd x)\dd t.
\end{equation}
The bound \eqref{eq:varying-mean-bound} passes to $\nu$ by truncating
the continuous random variable $\ip{\nu}{a}$.  For $\nu_X$ it follows
from \eqref{eq:SBM-local-mean}.  Hence, for either random measure
$\xi\in\{\nu,\nu_X\}$ and nonnegative test functions $f,h$,
\[
 \left|\E\left[e^{-\ip{\xi}{f}}\right]-\E\left[e^{-\ip{\xi}{h}}\right]\right|
 \leq\E\left[\ip{\xi}{|f-h|}\right]\leq C_*\|f-h\|_1.
\]
Approximating any $f\in C_c^+(E)$ by nonnegative smooth functions
supported in a common compact set extends
\eqref{eq:critical-subsequence-Laplace} to all such $f$.

The uniqueness assertion of Lemma~\ref{lem:Laplace-vague} now gives
$\operatorname{Law}(\nu)=\operatorname{Law}(\nu_X)$.
Thus the cluster Laplace limit identifies every subsequential law;
tightness yields \eqref{eq:varying-critical-convergence}.  Taking
$\ell_\epsilon=\ell=\lambda_0$ gives \eqref{eq:critical-convergence}.
\end{proof}

\subsection{Proof of the sparse limit}
\label{subsec:critical-sparse-proof}

\subsubsection{Reduction to a vanishing critical profile}

We now return to $0<\alpha<1$.  The elementary identity
\begin{equation}\label{eq:sparse-critical-identity}
  \ip{\cY_\epsilon}{g}
  =\epsilon^4\int_Eg(\epsilon^2t,\epsilon x)
       N^\epsilon(\dd t,\dd x)
\end{equation}
shows that $\cY_\epsilon$ is a critical-normalized family with the
triangular initial profile
\begin{equation}\label{eq:sparse-triangular-profile}
  \ell_\epsilon(x)=r_\epsilon\lambda_0(x),
  \qquad
  \epsilon^{d-2}\ell_\epsilon(\epsilon y)
  =\epsilon^{-(\alpha+1-d)}\lambda_0(\epsilon y).
\end{equation}
For a past parent $(s,y)$, let $\Xi_\epsilon^{s,y}$ assign mass $\epsilon^4$ to
each rescaled event in its parent clan, constructed by retaining only
direct positive-time children and then all their descendants as in
Proposition~\ref{prop:construction}.  Define its Poisson intensity measure
on $\cM_{\mathrm{loc}}(E)\setminus\{0\}$ by
\begin{equation}\label{eq:parent-clan-intensity}
  \mathfrak L_\epsilon(B)
  :=\epsilon^{-(\alpha+1-d)}\int_{s<0}\int_{\Rd}
       \lambda_0(\epsilon y)
       \mathbb{P}\left(\Xi_\epsilon^{s,y}\in B\right)\dd y\dd s.
\end{equation}

The limit profile is zero, but the first-order term in $r_\epsilon$ retains
the canonical cluster measure.  By the parent-clan construction in
\eqref{eq:parent-clan-intensity}, $\cY_\epsilon$ is the superposition of a
Poisson random measure of clans with intensity $\mathfrak L_\epsilon$.
Consequently,
\begin{equation}\label{eq:clan-Levy-exponent}
  -\log\E\left[e^{-\ip{\cY_\epsilon}{g}}\right]
  =\int_{\cM_{\mathrm{loc}}(E)\setminus\{0\}}
       \left(1-e^{-\ip{\nu}{g}}\right)
       \mathfrak L_\epsilon(\dd\nu).
\end{equation}

The Poisson clan representation isolates the first-order contribution
of the vanishing past intensity.  The next lemma identifies the
normalized clan Laplace exponent with the canonical occupation exponent.

\begin{lemma}[Canonical Laplace convergence]
\label{lem:canonical-parent-clan}
For every $g\in C_c^+(E)$,
\begin{equation}\label{eq:Levy-Laplace-convergence}
  \frac1{r_\epsilon}
  \int_{\cM_{\mathrm{loc}}(E)\setminus\{0\}}
       \left(1-e^{-\ip{\nu}{g}}\right)
       \mathfrak L_\epsilon(\dd\nu)
  \xrightarrow[\epsilon\downarrow0]{} \Psi_{\lambda_0}(g).
\end{equation}
\end{lemma}

\begin{proof}
First take $g$ smooth.  Use \eqref{eq:sparse-triangular-profile} in the
exact formula \eqref{eq:full-Laplace}.  In the proof of Lemma
\ref{lem:Laplace-boundary}, the leading term
\eqref{eq:boundary-leading-varying} is now
\[
  -r_\epsilon\int_{\Rd}\lambda_0(x)A_\epsilon(x)\dd x.
\]
The $L^1$ convergence of $A_\epsilon$ gives, after division by
$r_\epsilon$, the claimed limit.  Keeping the factor
$\ell_\epsilon=r_\epsilon\lambda_0$ in
\eqref{eq:Laplace-linearization-error} improves that error to
$\mathcal O(r_\epsilon\epsilon^2)$, so it remains negligible after the division.

To extend the statement to continuous test functions, note that for every
nonnegative compactly supported Borel function $a$,
\begin{equation}\label{eq:clan-first-moment-bound}
  \frac1{r_\epsilon}
  \int_{\cM_{\mathrm{loc}}(E)\setminus\{0\}}\ip{\nu}{a}\,
       \mathfrak L_\epsilon(\dd\nu)
  =\E\left[\ip{\cN_\epsilon}{a}\right]
  \leq\norm{\lambda_0}_\infty\int_E a(t,x)\dd t\dd x.
\end{equation}
The canonical first-moment formula
\cite[Theorems~II.7.2(iii) and II.7.3(b)]{perkins2002superprocess},
followed by the definition \eqref{eq:canonical-mixture}, similarly gives
\[
  \int_{\cM_{\mathrm{loc}}(E)\setminus\{0\}}
       \ip{\nu}{a}\,\mathbb Q_{\lambda_0}(\dd\nu)
  =\int_E a(t,x)\left(P_t^D\lambda_0\right)(x)\dd t\dd x
  \leq\norm{\lambda_0}_\infty\int_E a(t,x)\dd t\dd x.
\]
These bounds and $|e^{-b}-e^{-c}|\leq|b-c|$ for $b,c\geq0$ allow
uniform approximation by nonnegative smooth test functions with a common compact
support.  This proves \eqref{eq:Levy-Laplace-convergence} for all
$g\in C_c^+(E)$. The proof is complete.
\end{proof}

\subsubsection{Growing Laplace sources and the zero limit}

To recover the unmultiplied sparse observable, the Laplace source
must grow like $r_\epsilon^{-1}$.  We therefore need estimates in the next lemma
that remain uniform as the source amplitude increases.

\begin{lemma}[Bounds uniform in a growing source]
\label{lem:growing-source}
Fix $g\in C_c^\infty(E)$, $g\geq0$, with temporal support in $[0,T)$.
Let $u_{\epsilon,\theta}$ be the cluster solution for source $\theta g$.
There is $C_g<\infty$, independent of $\theta\geq1$ and sufficiently
small $\epsilon$, such that
\begin{align}
 \|u_{\epsilon,\theta}\|_\infty&\leq C_g\sqrt\theta,
 \qquad
 \sup_{0\leq t\leq T}\int_{\Rd}|x|^2u_{\epsilon,\theta}(t,x)\dd x
       \leq C_g\theta,
 \label{eq:growing-source-pointwise}\\
 \sup_{0\leq t\leq T}\|u_{\epsilon,\theta}(t,\cdot)\|_1
       &\leq C_g\theta^{(d+1)/(d+2)}.
 \label{eq:growing-source-L1}
\end{align}
\end{lemma}

\begin{proof}
Let $U_\theta:=\sqrt{2\theta\|g\|_\infty}+1$.
If $\epsilon^2U_\theta<1$, comparison as in
\eqref{eq:scheme-comparison} gives
$u_{\epsilon,\theta}\leq U_\theta\leq C_g\sqrt\theta$.
If $\epsilon^2U_\theta\geq1$, the definition \eqref{eq:u-epsilon},
with source $\theta g$ in \eqref{eq:cluster-Laplace}, gives
$0\leq u_{\epsilon,\theta}=(1-F_{\epsilon,\theta g})/\epsilon^2
\leq\epsilon^{-2}\leq U_\theta$, since $0\leq F_{\epsilon,\theta g}\leq1$.
Thus the same bound holds.
Let $b_{\epsilon,g}$ be the linear majorant
\eqref{eq:linear-cluster-majorant} with source $g$.
Linearity of that expression and \eqref{eq:u-weighted-tail} imply
\begin{equation}\label{eq:growing-source-weighted}
 0\leq u_{\epsilon,\theta}(t,x)\leq\theta b_{\epsilon,g}(t,x),
 \qquad
 \sup_{0\leq t\leq T}\int_{\Rd}\left(1+|x|^2\right)
                u_{\epsilon,\theta}(t,x)\dd x\leq C_g\theta.
\end{equation}
For $R\geq1$, let $B_R:=\{x\in\Rd:|x|\leq R\}$.
Using $\|u_{\epsilon,\theta}\|_\infty\leq C_g\sqrt\theta$ and
\eqref{eq:growing-source-weighted} gives
\begin{align*}
 \left\|u_{\epsilon,\theta}(t,\cdot)\right\|_1
 &=\int_{B_R}u_{\epsilon,\theta}(t,x)\dd x
       +\int_{\Rd\setminus B_R}u_{\epsilon,\theta}(t,x)\dd x\\
 &\leq |B_1|R^d\|u_{\epsilon,\theta}\|_\infty
       +R^{-2}\int_{\Rd\setminus B_R}|x|^2
                                  u_{\epsilon,\theta}(t,x)\dd x\\
 &\leq C_g\left(\sqrt\theta R^d+\theta R^{-2}\right),
       \qquad 0\leq t\leq T.
\end{align*}
For $R=\theta^{1/(2(d+2))}\geq1$, both powers of $\theta$ equal
$\theta^{(d+1)/(d+2)}$.  This proves \eqref{eq:growing-source-L1}.
\end{proof}


\begin{proof}[Proof of Theorem~\ref{thm:sparse-long}]
Equations \eqref{eq:clan-Levy-exponent} and
\eqref{eq:Levy-Laplace-convergence} prove
\eqref{eq:sparse-Laplace-first-order}.

To prove the zero limit, fix $g\in C_c^\infty(E)$ with $g\geq0$ and put
$\theta_\epsilon=r_\epsilon^{-1}$.  Then
$\theta_\epsilon\cY_\epsilon=\cN_\epsilon$, and Lemma
\ref{lem:growing-source} applies.  Let
$J_{\epsilon,\theta}$ be \eqref{eq:J-epsilon} with
$u_{\epsilon,\theta}$ in place of $u_\epsilon$.  Use the exact formula
\eqref{eq:full-Laplace} with profile $r_\epsilon\lambda_0$ and source
$\theta_\epsilon g$.  The inequality $1-e^{-a}\leq a$ and the substitution
$x=\epsilon y$ give
\begin{align*}
  0\leq-\log\E\left[e^{-\ip{\cN_\epsilon}{g}}\right]
  &\leq r_\epsilon\norm{\lambda_0}_\infty
    \int_{s<0}\int_{\Rd}
       J_{\epsilon,\theta_\epsilon}(s,x/\epsilon)\dd x\dd s\\
  &\leq r_\epsilon m_1\norm{\lambda_0}_\infty
      \sup_{t\geq0}\norm{u_{\epsilon,\theta_\epsilon}(t,\cdot)}_1.
\end{align*}
The second inequality follows by Tonelli's theorem, invariance of the
spatial integral under translations, and
\[
  \int_{s<0}\int_{q>0}\int_{\Rd} h(q-s,z)\dd z\dd q\dd s=m_1.
\]
Consequently,
\begin{equation}\label{eq:sparse-zero-Laplace-bound}
  -\log\E\left[e^{-\ip{\cN_\epsilon}{g}}\right]
  \leq C r_\epsilon^{1/(d+2)}\xrightarrow[\epsilon\downarrow0]{}0.
\end{equation}
For every $a>0$,
\[
  \mathbb{P}(\ip{\cN_\epsilon}{g}>a)
  \leq\frac{1-\E e^{-\ip{\cN_\epsilon}{g}}}{1-e^{-a}}
  \xrightarrow[\epsilon\downarrow0]{}0.
\]
Every nonnegative continuous compactly supported test function is
dominated by a nonnegative smooth compactly supported test function.
Consequently $L_{\cN_\epsilon}(f)\to1=L_0(f)$ for all $f\in C_c^+(E)$,
and the deterministic-limit assertion of Lemma~\ref{lem:Laplace-vague}
gives vague convergence in probability to zero.
The annealed assertion in \eqref{eq:sparse-zero}
is Proposition~\ref{prop:annealed-heat}.
\end{proof}

\section{Nearly Critical Case}
\label{sec:nearly-critical}

A perturbation of the branching ratio can survive under macroscopic
rescaling.  We treat finite-mean and infinite-mean reproduction delays
separately, retaining a prescribed Poisson past and a space--time event
observable in both cases.  The finite-mean theory reuses the critical
estimates with an added linear decay.  For a regularly varying temporal
tail, the space--time balance and past intensity change, and the
past leaves a temporal memory term.  The corresponding temporal scales
agree with \cite{jaisson2015limit,jaisson2016rough}.

\subsection{Finite-mean case}
\label{subsec:nc-light}

We retain the finite-mean assumptions
\eqref{eq:H1}--\eqref{eq:H4}; no exponential-moment assumption is
imposed.  We use the diffusive scaling and prescribed
Poisson past from Section~\ref{sec:model}.

\subsubsection{The perturbed kernel and its scale}

Let $h$ satisfy \eqref{eq:H1}--\eqref{eq:H4}, and suppose that
\begin{equation}\label{eq:nc-kernel}
  h_\epsilon(s,z):=a_\epsilon h(s,z),
  \qquad 0<a_\epsilon\leq1,
  \qquad
  \kappa_\epsilon:=\frac{1-a_\epsilon}{\epsilon^2}
       \xrightarrow[\epsilon\downarrow0]{}\kappa\in[0,\infty).
\end{equation}
Thus $a_\epsilon=1-\kappa\epsilon^2+o(\epsilon^2)$, and
$\|h_\epsilon\|_1=a_\epsilon$.  For $\kappa>0$ the kernels are
subcritical for all sufficiently small $\epsilon$.  The case $\kappa=0$
also includes deficits of order $o(\epsilon^2)$ and the exactly critical
choice $a_\epsilon\equiv1$.  The moments $m_1$ and $\Sigma$ below refer to
the fixed probability density $h$; conditional on an offspring event,
its reproduction delay and displacement have this density.

Use the same past field $\Pi^\epsilon$ as in
\eqref{eq:past-intensity}, and let $N^{\epsilon,a}$ have intensity
\begin{align}
  \lambda^{\epsilon,a}(t,x)
  ={}&a_\epsilon\int_{(-\infty,0)\times\Rd}
       h(t-s,x-y)\Pi^\epsilon(\dd s,\dd y)\nonumber\\
  &\qquad\qquad+a_\epsilon\int_{(0,t)\times\Rd}
       h(t-s,x-y)N^{\epsilon,a}(\dd s,\dd y).
  \label{eq:nc-intensity}
\end{align}
There is no additional positive-time immigration term.  Define
\begin{equation}\label{eq:nc-measure}
  \ip{\cN_\epsilon^a}{\varphi}
  :=\epsilon^{\alpha+3}
      \int_E\varphi(\epsilon^2t,\epsilon x)
           N^{\epsilon,a}(\dd t,\dd x),
  \qquad\varphi\in C_c(E).
\end{equation}
In particular, at the fluctuation threshold $\alpha=1$ the past
intensity is $\epsilon^{d-2}\lambda_0(\epsilon y)\dd s\dd y$ and each
rescaled event has mass $\epsilon^4$.

The forward offspring operator is $K_\epsilon^a=a_\epsilon K_\epsilon$.
For every $\varphi\in C_c^{1,2}([0,\infty)\times\Rd)$, Lemma
\ref{lem:generator} and \eqref{eq:nc-kernel} give
\begin{align}
  \frac{K_\epsilon^a\varphi-\varphi}{\epsilon^2}
  &=a_\epsilon\frac{K_\epsilon\varphi-\varphi}{\epsilon^2}
       -\kappa_\epsilon\varphi\nonumber\\
  &\xrightarrow[\epsilon\downarrow0]{}
       m_1\partial_t\varphi+\frac12\Sigma:D_x^2\varphi
       -\kappa\varphi
       =:(\cL-\kappa)\varphi
       \quad\text{in }L^1(E).
  \label{eq:nc-generator}
\end{align}
This identity determines the near-critical scale: a deficit of order
$\epsilon^2$ contributes on the same order as the first-order temporal
term and the centered second-order spatial term.  With observation horizon
$T=\epsilon^{-2}$, it is
exactly the condition $T(1-a_{T^{-1/2}})\to\kappa$ as $T\to\infty$,
corresponding to
\cite[equation~(3) and Theorem~2.2]{jaisson2015limit}.  Under the finite
moments used here, the physical spatial window has size $T^{1/2}$; the
branching deficit is of order $T^{-1}$ in every fixed dimension when
$\kappa>0$.  The
past intensity continues to account for the spatial volume factor.

\subsubsection{Limiting objects and convergence theorem}

Write
\begin{equation}\label{eq:nc-heat-profile}
  P_t^{D,\kappa}:=e^{-\kappa t/m_1}P_t^D,
  \qquad
  \rho^\kappa(t,\cdot):=P_t^{D,\kappa}\lambda_0.
\end{equation}
This is the bounded solution of
\begin{equation}\label{eq:nc-heat-equation}
  m_1\partial_t\rho^\kappa
       =\frac12\Sigma:D_x^2\rho^\kappa-\kappa\rho^\kappa,
  \qquad \rho^\kappa(0,\cdot)=\lambda_0.
\end{equation}
Let $X^\kappa$ be the canonical locally finite Dawson--Watanabe process
with initial measure $m_1\lambda_0(x)\dd x$ and characteristics
\begin{align}
  \cM_t^\kappa(\psi)
  &:={\ip{X_t^\kappa}{\psi}}-{\ip{X_0^\kappa}{\psi}}
       -\int_0^t\ip{X_s^\kappa}{
             \cA\psi-\frac{\kappa}{m_1}\psi}\dd s,
       \label{eq:nc-SBM-martingale}\\
  \left\langle\cM^\kappa(\psi)\right\rangle_t
  &=\frac1{m_1}\int_0^t\ip{X_s^\kappa}{\psi^2}\dd s,
       \qquad\psi\in C_c^2(\Rd).
       \label{eq:nc-SBM-bracket}
\end{align}
Here ``canonical'' refers to the finite-measure exhaustion construction,
as in Lemma~\ref{lem:locally-finite-SBM}; the construction with the added
linear decay is justified below.  Its occupation measure is
$\nu_{X^\kappa}(\dd t,\dd x)=m_1^{-1}X_t^\kappa(\dd x)\dd t$.

For $g\in C_c^+(E)$, choose $T$ beyond its temporal support and let
$u_g^\kappa$ denote the bounded nonnegative mild solution of
\begin{equation}\label{eq:nc-log-Laplace}
  \begin{cases}
  m_1\partial_tu+\frac12\Sigma:D_x^2u
        -\kappa u-\frac12u^2+g=0,&0\leq t<T,\\[0.3em]
  u(T,\cdot)=0.&
  \end{cases}
\end{equation}
Equivalently,
\begin{equation}\label{eq:nc-log-Laplace-mild}
  u_g^\kappa(t,\cdot)
  =\frac1{m_1}\int_t^T P_{s-t}^{D,\kappa}
     \left[g(s,\cdot)-\frac12\left(u_g^\kappa(s,\cdot)\right)^2\right]\dd s.
\end{equation}
The solution $u_g^\kappa$ defined by \eqref{eq:nc-log-Laplace-mild}
is independent of the terminal horizon beyond the temporal support of
$g$: if $T_2>T_1>\sup\{t:g(t,x)\ne0\text{ for some }x\}$, the
solution for $T_1$, extended by zero to $[T_1,T_2]$, satisfies the
same mild equation for $T_2$ and agrees with its unique solution.
Set
\begin{equation}\label{eq:nc-canonical-exponent}
  \Psi_{\lambda_0}^\kappa(g)
       :=m_1\int_{\Rd}\lambda_0(x)u_g^\kappa(0,x)\dd x.
\end{equation}
As in \eqref{eq:canonical-Laplace}, this is the canonical occupation
Laplace exponent, now for the superprocess with linear decay
$\kappa/m_1$.

The nearly critical perturbation adds linear decay to the limiting
heat flow and branching diffusion.  The fluctuation threshold remains
$\alpha=1$, as summarized by the following theorem.

\begin{theorem}[Nearly critical limits]\label{thm:nearly-critical}
Assume \eqref{eq:H1}--\eqref{eq:H4} and \eqref{eq:nc-kernel}, and let
$\lambda_0:\Rd\to\R_+$ be bounded and uniformly continuous.  For the
past intensity \eqref{eq:past-intensity} and process \eqref{eq:nc-intensity}:
\begin{enumerate}
  \item for every $\alpha>0$,
  \begin{equation}\label{eq:nc-annealed}
    \E\left[\cN_\epsilon^a\right]\xrightarrow[\epsilon\downarrow0]{v}
         \rho^\kappa(t,x)\dd t\dd x;
  \end{equation}
  \item if $\alpha>1$, then
  \begin{equation}\label{eq:nc-hydrodynamic}
    \cN_\epsilon^a\xrightarrow[\epsilon\downarrow0]{\mathbb P}
         \rho^\kappa(t,x)\dd t\dd x;
  \end{equation}
  \item if $\alpha=1$, then
  \begin{equation}\label{eq:nc-stochastic}
    \cN_\epsilon^a\xRightarrow[\epsilon\downarrow0]{}
         \frac1{m_1}X_t^\kappa(\dd x)\dd t;
  \end{equation}
  \item if $0<\alpha<1$, then
  \begin{equation}\label{eq:nc-sparse}
    \cN_\epsilon^a\xrightarrow[\epsilon\downarrow0]{\mathbb P}0,
    \qquad
    -\frac1{r_\epsilon}
       \log\E\left[e^{-\ip{r_\epsilon\cN_\epsilon^a}{g}}\right]
       \xrightarrow[\epsilon\downarrow0]{}\Psi_{\lambda_0}^\kappa(g)
       \quad(g\in C_c^+(E)),
  \end{equation}
  where $r_\epsilon=\epsilon^{1-\alpha}$.
\end{enumerate}
All random-measure convergences use the vague topology on
$\cM_{\mathrm{loc}}(E)$.  The stochastic assertion also holds for the
Poisson past with a varying profile at $\alpha=1$,
$\epsilon^{d-2}\ell_\epsilon(\epsilon y)\dd s\dd y$ under
\eqref{eq:varying-profile-assumption}; in that case the initial measure
of the limit is $m_1\ell(x)\dd x$.
\end{theorem}

\begin{proof}
Let $\mu_\epsilon^a:=\E\left[\cN_\epsilon^a\right]$ and
$C_0:=\|\lambda_0\|_\infty$.  For
$\varphi\in C_c^{1,2}([0,\infty)\times\Rd)$, define
$G_\epsilon^a\varphi:=\epsilon^{-2}(a_\epsilon K_\epsilon\varphi-\varphi)$.
Equations \eqref{eq:nc-master}, \eqref{eq:nc-mean-bound}, and
\eqref{eq:nc-generator} give
\begin{align*}
 \ip{\mu_\epsilon^a}{G_\epsilon^a\varphi}
       +a_\epsilon\E\left[\epsilon^{-2}\mathsf I_\epsilon(\varphi)\right]=0,
\end{align*}
and
\begin{align*}
 \E\left[\left|\ip{\cN_\epsilon^a}{G_\epsilon^a\varphi-(\cL-\kappa)\varphi}\right|\right]
       \leq C_0\left\|G_\epsilon^a\varphi-(\cL-\kappa)\varphi\right\|_1=o(1).
\end{align*}
The initial term converges by \eqref{eq:nc-boundary}, and
Lemma~\ref{lem:nc-linear-identification} then identifies every mean
subsequential limit as $\rho^\kappa(t,x)\dd t\dd x$.  This proves
\eqref{eq:nc-annealed}.
For $\alpha>1$, \eqref{eq:nc-noise-bound} also gives
\[
 \E\left[|\epsilon^{-2}M_\epsilon^a(\varphi)|^2\right]
       \leq C_\varphi\epsilon^{\alpha-1}\longrightarrow0.
\]
Tightness from Lemma~\ref{lem:nc-basic-estimates} and the random-limit
assertion of Lemma~\ref{lem:nc-linear-identification} therefore imply
that every random subsequential limit is the same deterministic
measure.  The metric argument in the proof of
Theorem~\ref{thm:main} yields \eqref{eq:nc-hydrodynamic}.

At $\alpha=1$, allow the varying profiles in
\eqref{eq:varying-profile-assumption}.  For smooth $f\geq0$, equations
\eqref{eq:nc-Laplace-limit} and \eqref{eq:nc-occupation-Laplace} give
\begin{align*}
 \lim_{\epsilon\downarrow0}\E\left[e^{-\ip{\cN_\epsilon^{a,\ell_\epsilon}}{f}}\right]
 &=\exp\!\left\{-m_1\int_{\Rd}\ell(x)u_f^\kappa(0,x)\dd x\right\}\\
 &=\E\left[\exp\!\left\{-\frac1{m_1}\int_0^T
                      \ip{X_t^\kappa}{f(t,\cdot)}\dd t\right\}\right].
\end{align*}
For $f,h\in C_c^+(E)$, the mean bound \eqref{eq:nc-mean-bound},
in its varying-profile version from Lemma~\ref{lem:nc-basic-estimates},
gives, with $C_*:=\sup_{0<\epsilon\leq1}\|\ell_\epsilon\|_\infty<\infty$
by \eqref{eq:varying-profile-assumption},
\begin{align*}
 \left|\E\left[e^{-\ip{\cN_\epsilon^{a,\ell_\epsilon}}{f}}\right]
       -\E\left[e^{-\ip{\cN_\epsilon^{a,\ell_\epsilon}}{h}}\right]\right|
 &\leq\E\left[\ip{\cN_\epsilon^{a,\ell_\epsilon}}{|f-h|}\right]\\
 &\leq C_*\|f-h\|_{L^1(E)},\qquad 0<\epsilon\leq1.
\end{align*}
The limiting occupation measure
$\nu_{X^\kappa}(\dd t,\dd x)=m_1^{-1}X_t^\kappa(\dd x)\dd t$
satisfies, for any compactly supported Borel function $a\geq0$,
\[
 \E[\ip{\nu_{X^\kappa}}a]
 =\int_0^\infty\ip\ell{P_t^{D,\kappa}a(t,\cdot)}\dd t
 \leq C_*\|a\|_{L^1(E)},
\]
by the first-moment formula \eqref{eq:nc-SBM-local-mean} and
Tonelli's theorem.

Smooth approximation therefore extends the Laplace convergence to
$C_c^+(E)$.  Lemma~\ref{lem:Laplace-vague}, with limiting random measure
$\nu_{X^\kappa}$, proves \eqref{eq:nc-stochastic} and its
varying-profile version.

For $0<\alpha<1$, let $r_\epsilon=\epsilon^{1-\alpha}$.
The measure $r_\epsilon\cN_\epsilon^a$ uses event mass $\epsilon^4$
and initial profile $r_\epsilon\lambda_0$ at the fluctuation threshold, by
\eqref{eq:sparse-triangular-profile}.  Equations
\eqref{eq:nc-Laplace-leading} and \eqref{eq:nc-preparation-L1} yield
\begin{align*}
 -\frac1{r_\epsilon}\log\E\left[e^{-\ip{r_\epsilon\cN_\epsilon^a}{f}}\right]
 =a_\epsilon\int_{\Rd}\lambda_0(x)A_\epsilon^a(x)\dd x
                        +\mathcal O(\epsilon^2),
\end{align*}
and
\begin{align*}
 \left|a_\epsilon\int_{\Rd}\lambda_0A_\epsilon^a\dd x
             -\Psi_{\lambda_0}^\kappa(f)\right|
 \leq\|\lambda_0\|_\infty\left(
 a_\epsilon\left\|A_\epsilon^a-m_1u_f^\kappa(0)\right\|_1
       +(1-a_\epsilon)m_1\left\|u_f^\kappa(0)\right\|_1\right)
 \longrightarrow0.
\end{align*}
To extend this to continuous $g\geq0$, let $\mathfrak L_\epsilon^a$ be
the clan intensity defined as in \eqref{eq:parent-clan-intensity}, with
the thinned clans.  The exact Poisson exponent gives
\begin{align*}
 &\left|\frac1{r_\epsilon}\log\E\left[e^{-\ip{r_\epsilon\cN_\epsilon^a}{f}}\right]
       -\frac1{r_\epsilon}\log\E\left[e^{-\ip{r_\epsilon\cN_\epsilon^a}{h}}\right]\right|\\
 &\quad\leq\frac1{r_\epsilon}
       \int_{\cM_{\mathrm{loc}}(E)\setminus\{0\}}
              \ip{\xi}{|f-h|}\,\mathfrak L_\epsilon^a(\dd\xi)
 =\E\left[\ip{\cN_\epsilon^a}{|f-h|}\right]
 \leq C_0\|f-h\|_1.
\end{align*}
Equation \eqref{eq:nc-canonical-identification} gives the same bound
for $|\Psi_{\lambda_0}^\kappa(f)-\Psi_{\lambda_0}^\kappa(h)|$.
Approximation in a common compact support proves the exponent limit
in \eqref{eq:nc-sparse}.  Finally, \eqref{eq:nc-domination} and
Theorem~\ref{thm:sparse-long} give, for $g\in C_c^+(E)$ and $a>0$,
\[
 \mathbb P(\ip{\cN_\epsilon^a}{g}>a)
 \leq\mathbb P(\ip{\cN_\epsilon}{g}>a)\longrightarrow0.
\]
A countable family of test functions determining the vague topology proves the
zero limit in \eqref{eq:nc-sparse}. This completes the proof.
\end{proof}

Thus the fluctuation threshold remains $\alpha=1$.  In the sparse regime,
\eqref{eq:nc-annealed} and \eqref{eq:nc-sparse} again imply loss of uniform
integrability whenever a nonnegative test function has positive limiting mean.

\subsubsection{Construction, moment bounds, and the deterministic limits}\label{sec:nearly:support:1}

Thinning transfers the construction and moment estimates of the
critical process to the perturbed kernel.  We collect the identities
needed to identify its deterministic limits.

\begin{lemma}[Basic estimates for nearly critical kernels]
\label{lem:nc-basic-estimates}
Assume \eqref{eq:H1}--\eqref{eq:H4} and \eqref{eq:nc-kernel}, and
let $\lambda_0$ be bounded and uniformly continuous.
The process is locally finite and unique in the mean-intensity class
\eqref{eq:Hawkes-mean-class}.  It admits a coupling with the
critical process such that
\begin{equation}\label{eq:nc-domination}
  N^{\epsilon,a}\leq N^\epsilon,
  \qquad
  \E\left[\lambda^{\epsilon,a}(t,x)\right]
       \leq\epsilon^{-(\alpha+1-d)}\|\lambda_0\|_\infty.
\end{equation}
Its rescaled event measures satisfy
\begin{equation}\label{eq:nc-mean-bound}
  \E\left[\ip{\cN_\epsilon^a}{f}\right]
       \leq\|\lambda_0\|_\infty
          \int_E f(t,x)\dd t\dd x,
       \qquad f\in C_c^+(E),
\end{equation}
and their laws are tight in $\cM_{\mathrm{loc}}(E)$.
Let $M_\epsilon^a$ denote the compensated event integral with
normalization \eqref{eq:scaled-measure}, and let
$\mathsf I_\epsilon$ be \eqref{eq:boundary-def} for the common
past and the base kernel $h$.  For
$\varphi\in C_c^{1,2}([0,\infty)\times\Rd)$,
\begin{align}
&\ip{\cN_\epsilon^a}{
       \frac{a_\epsilon K_\epsilon\varphi-\varphi}{\epsilon^2}}
  +a_\epsilon\frac{\mathsf I_\epsilon(\varphi)}{\epsilon^2}
  +\frac{M_\epsilon^a(\varphi)}{\epsilon^2}=0,\label{eq:nc-master}
\\
& a_\epsilon\frac{\mathsf I_\epsilon(\varphi)}{\epsilon^2}
       \xrightarrow[\epsilon\downarrow0]{}
       m_1\int_{\Rd}\varphi(0,x)\lambda_0(x)\dd x
       \quad\text{in }L^2,\label{eq:nc-boundary}
\\
&\E\left[\left|\epsilon^{-2}M_\epsilon^a(\varphi)\right|^2\right]
       \leq C_\varphi\epsilon^{\alpha-1}.\label{eq:nc-noise-bound}
\end{align}
The construction and bounds also apply to uniformly bounded varying
profiles, with their common bound replacing $\|\lambda_0\|_\infty$;
the boundary limit uses the locally uniform limiting profile.
\end{lemma}

\begin{proof}
Retain each offspring edge of a critical cluster independently with
probability $a_\epsilon$, retaining a vertex only when every edge on
its ancestral line from the past parent is retained.  Independent
Poisson thinning gives offspring intensity $a_\epsilon h$ and
$N^{\epsilon,a}\leq N^\epsilon$.  With
$c_\epsilon=\epsilon^{-(\alpha+1-d)}\|\lambda_0\|_\infty$,
\eqref{eq:construction-generation-bound} yields
\[
 \E[\lambda^{\epsilon,a}(t,x)]
 \leq a_\epsilon c_\epsilon\overline H(t)
       +a_\epsilon c_\epsilon\int_0^tH(s)\dd s
 \leq c_\epsilon.
\]
This proves \eqref{eq:nc-domination} and local finiteness.  For two
solutions driven by the same Poisson measure, their expected intensity
difference satisfies
\[
 q(t,x)\leq a_\epsilon\int_0^t\int_{\Rd}
                           h(t-s,x-y)q(s,y)\dd y\dd s.
\]
With $q_T^*$ and $k_b<1$ as in the proof of
Proposition~\ref{prop:construction}, this gives
$q_T^*\leq a_\epsilon k_bq_T^*$ and hence $q=0$.
The Poisson embedding and Picard construction in Proposition~\ref{prop:construction}
therefore give pathwise uniqueness and uniqueness in law.

For nonnegative compactly supported Borel $f$, changing
$t=\epsilon^2s$, $x=\epsilon y$ gives
\begin{align*}
 \E[\ip{\cN_\epsilon^a}f]
 &=\epsilon^{3+\alpha}\int_Ef(\epsilon^2s,\epsilon y)
                       \E[\lambda^{\epsilon,a}(s,y)]\dd s\dd y\\
 &\leq\epsilon^{3+\alpha}c_\epsilon\epsilon^{-2-d}
                         \|f\|_{L^1(E)}
 =\|\lambda_0\|_\infty\|f\|_{L^1(E)}.
\end{align*}
This is \eqref{eq:nc-mean-bound}.  For a compact exhaustion $K_j$ of
$E$, choose $A_j>0$ so that
$\|\lambda_0\|_\infty|K_j|/A_j\leq\zeta2^{-j}$.  Then
\[
 \sup_{0<\epsilon\leq1}\mathbb P
   (\cN_\epsilon^a(K_j)>A_j\text{ for some }j\geq1)
 \leq\sum_{j=1}^\infty
             \frac{\|\lambda_0\|_\infty|K_j|}{A_j}\leq\zeta.
\]
The compactness argument of Proposition~\ref{prop:tightness} proves
tightness.

Define
$\lambda_\epsilon^a(t,x):=\epsilon^{\alpha+1-d}
\lambda^{\epsilon,a}(t/\epsilon^2,x/\epsilon)$.
Substituting the past-parent and positive-time-parent terms of
\eqref{eq:nc-intensity} into this rescaled intensity and applying
Fubini's theorem gives the analogue of \eqref{eq:compensator}:
\begin{equation}\label{eq:nc-compensator}
 \int_E\varphi\lambda_\epsilon^a\dd t\dd x
 =a_\epsilon\ip{\cN_\epsilon^a}{K_\epsilon\varphi}
                         +a_\epsilon\mathsf I_\epsilon(\varphi).
\end{equation}
The compensated integral is
\begin{equation}\label{eq:nc-compensated-integral}
 M_\epsilon^a(\varphi)
 =\ip{\cN_\epsilon^a}\varphi-
                              \int_E\varphi\lambda_\epsilon^a\dd t\dd x.
\end{equation}
Substituting \eqref{eq:nc-compensator} into
\eqref{eq:nc-compensated-integral} and dividing by $\epsilon^2$ proves
\eqref{eq:nc-master}.  If
$B_\varphi=m_1\int_{\Rd}\varphi(0,x)\lambda_0(x)\dd x$, then
\eqref{eq:boundary-convergence} implies
\[
 \left\|a_\epsilon\epsilon^{-2}\mathsf I_\epsilon(\varphi)-B_\varphi\right\|_{L^2}
 \leq a_\epsilon\left\|\epsilon^{-2}\mathsf I_\epsilon(\varphi)-B_\varphi\right\|_{L^2}
          +(1-a_\epsilon)|B_\varphi|\longrightarrow0.
\]
This proves \eqref{eq:nc-boundary}.
Here $\langle M_\epsilon^a(\varphi)\rangle_{\mathrm{pred}}$ denotes
the terminal predictable quadratic variation of the compensated
event martingale in \eqref{eq:nc-compensated-integral}.
The nearly critical version of \eqref{eq:martingale-bracket} is
\begin{equation}\label{eq:nc-martingale-bracket}
 \left\langle M_\epsilon^a(\varphi)\right\rangle_{\mathrm{pred}}
 =\epsilon^{3+\alpha}\int_E\left(\varphi(t,x)\right)^2
                                    \lambda_\epsilon^a(t,x)\dd t\dd x.
\end{equation}
Taking expectations in \eqref{eq:nc-martingale-bracket} and using
\eqref{eq:nc-domination} gives
\begin{align*}
 \E\left[\left|\epsilon^{-2}M_\epsilon^a(\varphi)\right|^2\right]
 &=\epsilon^{\alpha-1}\int_E\left(\varphi(t,x)\right)^2
                                      \E\left[\lambda_\epsilon^a(t,x)\right]\dd t\dd x\\
 &\leq\epsilon^{\alpha-1}\|\lambda_0\|_\infty\|\varphi\|_2^2,
\end{align*}
which proves \eqref{eq:nc-noise-bound}.

For profiles $\ell_\epsilon$ satisfying
\eqref{eq:varying-profile-assumption}, all these inequalities hold
with $C_*:=\sup_{0<\epsilon\leq1}\|\ell_\epsilon\|_\infty<\infty$.
Writing $\mathsf I_\epsilon^{\ell_\epsilon}$ for the corresponding
boundary term, \eqref{eq:boundary-mean} becomes
\begin{align*}
 \E\left[\epsilon^{-2}\mathsf I_\epsilon^{\ell_\epsilon}(\varphi)\right]
 &=\int_0^\infty\int_q^\infty\int_{\Rd}\int_{\Rd}
   \varphi(\epsilon^2q,x)h(r,z)\ell_\epsilon(x-\epsilon z)
                                          \dd x\dd z\dd r\dd q\\
 &\longrightarrow m_1\int_{\Rd}\varphi(0,x)\ell(x)\dd x.
\end{align*}
For fixed $(q,z)$ the spatial integral converges by local uniform
convergence, and it is bounded by
$C_*\sup_{t\geq0}\|\varphi(t,\cdot)\|_1$.
The crossing measure has mass $m_1$, so the dominated convergence
theorem applies.  Equation \eqref{eq:boundary-variance} also gives
\[
 \operatorname{Var}\left(\epsilon^{-2}
           \mathsf I_\epsilon^{\ell_\epsilon}(\varphi)\right)
 \leq C_\varphi C_*m_1\epsilon^{\alpha+1}\longrightarrow0.
\]
Multiplication by $a_\epsilon\to1$ proves the varying-profile
boundary limit.
\end{proof}

The next lemma identifies subsequential limits through the linear
weak equation.  Its growth condition is satisfied by the moment
bound above.

\begin{lemma}[Identification of the linear limits]
\label{lem:nc-linear-identification}
Under the assumptions of Lemma~\ref{lem:nc-basic-estimates}, the
equation
\begin{equation}\label{eq:nc-weak-equation}
  \ip{\nu}{(\cL-\kappa)\varphi}
       +m_1\int_{\Rd}\varphi(0,x)\lambda_0(x)\dd x=0
\end{equation}
for all $\varphi\in C_c^{1,2}([0,\infty)\times\Rd)$ has the unique
nonnegative Radon solution
$\nu(\dd t,\dd x)=\rho^\kappa(t,x)\dd t\dd x$ in the
polynomial-growth class of Proposition~\ref{prop:uniqueness}.
Every vague subsequential limit of the mean measures belongs to
this class and satisfies the equation.  If $\alpha>1$, the same
holds almost surely for every subsequential random-measure limit.
\end{lemma}

\begin{proof}
Let $G_\epsilon^a\varphi:=\epsilon^{-2}(a_\epsilon K_\epsilon\varphi-\varphi)$
and $G\varphi:=(\cL-\kappa)\varphi$.
Equations \eqref{eq:nc-generator} and \eqref{eq:nc-mean-bound} imply
\begin{equation}\label{eq:nc-generator-pairing-error}
 \E\left[\left|\ip{\cN_\epsilon^a}{G_\epsilon^a\varphi-G\varphi}\right|\right]
 \leq C_0\left\|G_\epsilon^a\varphi-G\varphi\right\|_1=o(1),
 \qquad C_0:=\|\lambda_0\|_\infty.
\end{equation}
Consider a vague subsequential limit $\mu$ of the mean measures.
The same mean bound gives $\mu(A)\leq C_0|A|$ for Borel $A\subset E$.
If $\zeta_\tau$ vanishes on $[0,\tau]$ and equals one on
$[2\tau,\infty)$, then $\zeta_\tau G\varphi\in C_c(E)$.
For $B_R$ containing the spatial support of $\varphi$ and
$\widetilde\mu\in\{\E[\cN_\epsilon^a],\mu\}$,
\begin{equation}\label{eq:nc-initial-strip-error}
 |\ip{\widetilde\mu}{(1-\zeta_\tau)G\varphi}|
 \leq2C_0\tau|B_R|\|G\varphi\|_\infty.
\end{equation}
Taking expectations in \eqref{eq:nc-master}, passing to the vague
limit with this cutoff, and then letting $\tau\downarrow0$, equations
\eqref{eq:nc-generator-pairing-error} and \eqref{eq:nc-boundary} give
\eqref{eq:nc-weak-equation} for $\mu$.  Also
$\mu((0,T]\times B_R)\leq C_0T|B_R|$, so $\mu$ has polynomial growth.

For uniqueness, let $g\in C_c^\infty(E)$ and define
\begin{equation}\label{eq:nc-backward-test}
 \varphi(t,\cdot):=\frac1{m_1}\int_t^T
                 P_{s-t}^{D,\kappa}g(s,\cdot)\dd s,
 \qquad (\cL-\kappa)\varphi=-g,
\end{equation}
where $T$ exceeds the temporal support of $g$.
The factor $e^{-\kappa(s-t)/m_1}\leq1$ preserves the Gaussian bound
\eqref{eq:backward-test-tail}.  For the cutoff $\chi_R$ used in
\eqref{eq:backward-cutoff-error}, the product rule yields
\[
 (\cL-\kappa)(\chi_R\varphi)
 =-\chi_Rg+E_R,\qquad
 E_R=\varphi\,\frac12\Sigma:D_x^2\chi_R
           +(\Sigma\nabla\chi_R)\cdot\nabla_x\varphi.
\]
Thus, for a solution $\nu$ satisfying the polynomial growth condition
\eqref{eq:limiting-measure-growth},
\begin{align*}
 |\ip{\nu}{E_R}|
 &\leq C\left(R^{-1}+R^{-2}\right)e^{-cR^2}\nu((0,T]\times B_{2R})\\
 &\leq C C_T\left(R^{-1}+R^{-2}\right)(1+2R)^{p_T}e^{-cR^2}
       \longrightarrow0.
\end{align*}
Insert $\chi_R\varphi$ in \eqref{eq:nc-weak-equation} and let
$R\to\infty$.  The integrable initial pairing and symmetry of
$P^{D,\kappa}$ give
\[
 \ip{\nu}{g}
 =m_1\ip{\lambda_0}{\varphi(0)}
 =\int_0^T\int_{\Rd}g(t,x)\left(P_t^{D,\kappa}\lambda_0\right)(x)\dd x\dd t.
\]
This proves uniqueness and identifies the mean limit.

If $\alpha>1$, take a subsequential random limit $\nu$ using tightness.
The first-moment inequality passes to it by truncation and gives
$\E[\nu(A)]\leq C_0|A|$.  In particular
\eqref{eq:nc-initial-strip-error} holds in expectation for $\nu$.
For each fixed cutoff the pairing converges in law; letting
$\tau\downarrow0$ therefore gives
$\ip{\cN_\epsilon^a}{G\varphi}\Rightarrow\ip{\nu}{G\varphi}$ along
that subsequence.  The error
\eqref{eq:nc-generator-pairing-error}, the boundary limit
\eqref{eq:nc-boundary}, and
\[
 \E\left[\left|\epsilon^{-2}M_\epsilon^a(\varphi)\right|^2\right]
       \leq C_\varphi\epsilon^{\alpha-1}\longrightarrow0
\]
then imply \eqref{eq:nc-weak-equation} almost surely, by
\eqref{eq:nc-master}.  A countable family dense in the $C^{1,2}$ norm
on each compact support gives one probability-one event for all test functions.
Finally,
\[
 \E\left[\nu((0,n]\times B_n)\right]\leq Cn^{d+1},\qquad
 \mathbb P\left(\nu((0,n]\times B_n)>n^{d+3}\right)\leq Cn^{-2}.
\]
The Borel--Cantelli lemma gives a finite $C(\omega)$ such that
$\nu((0,n]\times B_n)\leq C(\omega)n^{d+3}$ for all $n\geq1$.
For arbitrary $T<\infty$ and $R\geq1$, choosing
$n=\lceil\max\{T,R,1\}\rceil$ gives
\[
 \nu((0,T]\times B_R)\leq C_T(\omega)(1+R)^{d+3}.
\]
The uniqueness calculation with \eqref{eq:nc-backward-test} thus
applies pathwise and identifies every random limit.
This completes the proof.
\end{proof}

\subsubsection{The stochastic occupation limit and sparse expansion}\label{sec:nearly:support:2}

The single-cluster estimates also control initial profiles whose
amplitude tends to zero.  The error below retains that amplitude,
which is needed in the sparse expansion.

For $f\in C_c^\infty([0,T)\times\Rd)$ with $f\geq0$, define
$F_\epsilon^a$ from a full cluster with offspring kernel $a_\epsilon h$
as in \eqref{eq:cluster-Laplace}, and put
$u_\epsilon^a=(1-F_\epsilon^a)/\epsilon^2$.  The branching equation is
\begin{equation}\label{eq:nc-renewal}
  \log\left(1-\epsilon^2u_\epsilon^a\right)
       =-\epsilon^4f-a_\epsilon\epsilon^2K_\epsilon u_\epsilon^a,
       \qquad u_\epsilon^a(T,\cdot)=0.
\end{equation}

To pass from the single-cluster equation to the full Hawkes field, we
integrate the cluster cost over the negative-time Poisson parents.
The next lemma identifies this contribution of the Poisson past and retains an
error bound uniform over bounded initial profiles, as required for
the sparse expansion.

\begin{lemma}[Nearly critical contribution of the Poisson past]
\label{lem:nc-preparation}
Assume \eqref{eq:H1}--\eqref{eq:H4} and \eqref{eq:nc-kernel}, and
let $f$ and $u_\epsilon^a$ be as above.
Define
\[
  A_\epsilon^a(x)
       :=\int_{r<0}\int_{q>0}\int_{\Rd}
          h(q-r,z)u_\epsilon^a\left(\epsilon^2q,x+\epsilon z\right)
          \dd z\dd q\dd r.
\]
For every uniformly bounded family of nonnegative initial profiles
$\ell_\epsilon$, the event measure with past intensity
$\epsilon^{d-2}\ell_\epsilon(\epsilon y)\dd s\dd y$ and event mass
$\epsilon^4$ satisfies
\begin{equation}\label{eq:nc-Laplace-leading}
  \log\E\left[e^{-\ip{\cN_\epsilon^{a,\ell_\epsilon}}{f}}\right]
       =-a_\epsilon\int_{\Rd}\ell_\epsilon(x)A_\epsilon^a(x)\dd x
          +\mathcal O\left(\epsilon^2\|\ell_\epsilon\|_\infty\right),
\end{equation}
where the implicit constant is independent of the profiles, and
\begin{equation}\label{eq:nc-preparation-L1}
 \left\|A_\epsilon^a-m_1u_f^\kappa(0,\cdot)\right\|_1\xrightarrow[\epsilon\downarrow0]{}0.
\end{equation}
If also $\ell_\epsilon\to\ell$ locally uniformly, then
\begin{equation}\label{eq:nc-Laplace-limit}
  \log\E\left[ e^{-\ip{\cN_\epsilon^{a,\ell_\epsilon}}{f}}\right]
       \xrightarrow[\epsilon\downarrow0]{}
       -m_1\int_{\Rd}\ell(x)u_f^\kappa(0,x)\dd x.
\end{equation}
\end{lemma}

\begin{proof}
Lemmas~\ref{lem:nonlinear-renewal} and
\ref{lem:spatial-localization}, applied to $a_\epsilon h$, give
\begin{equation}\label{eq:nc-cluster-convergence}
  \left\|u_\epsilon^a-u_f^\kappa\right\|_\infty\xrightarrow[\epsilon\downarrow0]{}0,
  \qquad
  \sup_{0\leq t\leq T}
       \left\|u_\epsilon^a(t,\cdot)-u_f^\kappa(t,\cdot)\right\|_1\xrightarrow[\epsilon\downarrow0]{}0.
\end{equation}
This includes the initial time and the case $\kappa=0$ under
\eqref{eq:H1}--\eqref{eq:H4} alone.

For the Poisson past with a varying profile at the fluctuation threshold, define
\[
  J_\epsilon^a(r,y)
       :=\int_{q>0}\int_{\Rd}
          h(q-r,w-y)u_\epsilon^a\left(\epsilon^2q,\epsilon w\right)\dd w\dd q.
\]
The positive-time children of a past parent have intensity
$a_\epsilon h$.  Their Laplace contribution is therefore
$\exp\left\{-a_\epsilon\epsilon^2J_\epsilon^a(r,y)\right\}$, and the exact
Poisson formula is
\begin{equation}\label{eq:nc-full-Laplace}
  \log\E\left[e^{-\ip{\cN_\epsilon^{a,\ell_\epsilon}}{f}}\right]
       =\epsilon^{d-2}\int_{r<0}\int_{\Rd}
          \ell_\epsilon(\epsilon y)
          \left(e^{-a_\epsilon\epsilon^2J_\epsilon^a(r,y)}-1\right)
          \dd y\dd r.
\end{equation}
Here $\cN_\epsilon^{a,\ell_\epsilon}$ uses event mass $\epsilon^4$.
Linearization and the change of variables $x=\epsilon y$ give
\eqref{eq:nc-Laplace-leading}.
Indeed, the quadratic error is bounded by
$a_\epsilon^2\epsilon^{d+2}\|\ell_\epsilon\|_\infty/2$ times the
integral of $(J_\epsilon^a)^2$ in $(r,y)$.
By \eqref{eq:u-uniform-bound} and \eqref{eq:u-weighted-tail},
$0\leq J_\epsilon^a\leq\|u_\epsilon^a\|_\infty\leq U_f$ and
$\sup_{0\leq t\leq T}\|u_\epsilon^a(t,\cdot)\|_1\leq C_{f,T}$.
Tonelli's theorem and the change of variables $x=\epsilon y$ therefore give
\begin{align*}
  \int_{-\infty}^0\int_{\Rd}(J_\epsilon^a(r,y))^2\dd y\dd r
  &\leq U_f\int_{-\infty}^0\int_{\Rd}J_\epsilon^a(r,y)\dd y\dd r\\
  &\leq U_f\epsilon^{-d}m_1
                  \sup_{0\leq t\leq T}\|u_\epsilon^a(t,\cdot)\|_1
  \leq C_{f,T}\epsilon^{-d}m_1.
\end{align*}
This also justifies \eqref{eq:nc-full-Laplace} for the infinite past:
first apply \eqref{eq:Poisson-exponential-formula} to bounded sets of
parents with
$\Lambda(\dd r,\dd y)=\epsilon^{d-2}\ell_\epsilon(\epsilon y)\dd r\dd y$
and $L(r,y)=e^{-a_\epsilon\epsilon^2J_\epsilon^a(r,y)}$, and then
exhaust the past, using $1-e^{-z}\leq z$ and the finite integral above.
The measure $h(q-r,z)\dd r\dd q\dd z$ on $r<0<q$ has total mass
$m_1$.  For each fixed $(q,z)$, \eqref{eq:nc-cluster-convergence}
and continuity of $t\mapsto u_f^\kappa(t,\cdot)$ in $L^1$ give
\[
  \left\|u_\epsilon^a\left(\epsilon^2q,\cdot+\epsilon z\right)
                           -u_f^\kappa(0,\cdot)\right\|_1\to0.
\]
The common $L^1$ bound permits the dominated convergence theorem against that
finite measure, proving \eqref{eq:nc-preparation-L1}.
Since $a_\epsilon\to1$, bounded locally uniform convergence of the
profiles and spatial localization now yield
\eqref{eq:nc-Laplace-limit}.
In particular, the extra factor in the contribution of the prescribed past tends to one:
the initial measure is still $m_1\ell(x)\dd x$.
\end{proof}

The limiting Laplace exponent of the Poisson past is the occupation exponent of a
branching diffusion with linear decay.  We record its construction
and the canonical measure used in the sparse regime.

\begin{lemma}[Occupation law and canonical exponent]
\label{lem:nc-occupation-law}
For every bounded nonnegative initial profile $\ell$, there is a
canonical locally finite process $X^\kappa$, with initial measure
$m_1\ell(x)\dd x$, satisfying
\eqref{eq:nc-SBM-martingale}--\eqref{eq:nc-SBM-bracket}.
Its law is independent of the finite-measure exhaustion, and its
locally finite occupation measure satisfies, for every $T<\infty$
and $f\in C_c^\infty([0,T)\times\Rd)$ with $f\geq0$,
\begin{equation}\label{eq:nc-occupation-Laplace}
  \E\left[\exp\!\left\{-\frac1{m_1}\int_0^T
                   \ip{X_t^\kappa}{f(t,\cdot)}\dd t\right\}\right]
       =\exp\!\left\{-\ip{X_0^\kappa}{u_f^\kappa(0,\cdot)}\right\}.
\end{equation}
For $\ell=\lambda_0$, there is a sigma-finite canonical occupation
measure $\mathbb Q_{\lambda_0}^\kappa$ on nonzero Radon measures
such that, for every $g\in C_c^+(E)$,
\begin{equation}\label{eq:nc-canonical-identification}
 \begin{aligned}
  &\int_{\cM_{\mathrm{loc}}(E)\setminus\{0\}}
       \left(1-e^{-\ip\nu g}\right)\mathbb Q_{\lambda_0}^\kappa(\dd\nu)
       =\Psi_{\lambda_0}^\kappa(g),\\
  &\int_{\cM_{\mathrm{loc}}(E)\setminus\{0\}}
       \ip\nu g\,\mathbb Q_{\lambda_0}^\kappa(\dd\nu)
       =\int_E g(t,x)\rho^\kappa(t,x)\dd t\dd x
       \leq\|\lambda_0\|_\infty\|g\|_1.
 \end{aligned}
\end{equation}
\end{lemma}

\begin{proof}
For a finite initial measure $\mu_0$, apply
\cite[Theorem~II.5.1(a)--(b)]{perkins2002superprocess} with Brownian
motion generator $\cA$, branching variance $1/m_1$, and the bounded
continuous linear drift $-\kappa/m_1$.  With
$\varphi_s=P_{t-s}^{D,\kappa}\varphi$, one has
\[
 \partial_s\varphi_s+(\cA-\kappa/m_1)\varphi_s=0.
\]
The time-dependent martingale formula, as in
\eqref{eq:SBM-finite-mean-identity}, therefore gives
\begin{equation}\label{eq:nc-SBM-first-moment}
 \E_{\mu_0}[\ip{X_t^\kappa}\varphi]
       =\ip{\mu_0}{P_t^{D,\kappa}\varphi},\qquad \varphi\geq0.
\end{equation}
The identity extends from smooth test functions to nonnegative Borel functions
by equality of mean measures and monotone convergence theorem.

Let $\mu(\dd x)=m_1\ell(x)\dd x$ and construct independent components
from $\mu|_{A_j}$, using the annuli of
Lemma~\ref{lem:locally-finite-SBM}.  Write $X^{\kappa,[n]}$ for their
first $n$ components and
$Y^{n,m}=X^{\kappa,[m]}-X^{\kappa,[n]}$.
For $w(x)=(1+|x|^2)^{-p}$, $p>d/2$, choose $c\geq0$ with
$(\cA-\kappa/m_1)w\leq cw$.  The finite-measure decomposition gives
\[
 \dd\left(e^{-ct}\ip{Y_t^{n,m}}w\right)
 =e^{-ct}\dd M_t^{n,m}(w)
   +e^{-ct}\ip{Y_t^{n,m}}{(\cA-\kappa/m_1-c)w}\dd t.
\]
The drift is nonpositive.  Thus the nonnegative supermartingale
maximal inequality used in \eqref{eq:SBM-weighted-tail} yields
\begin{equation}\label{eq:nc-SBM-weighted-tail}
 \mathbb P\left(\sup_{0\leq t\leq T}\ip{Y_t^{n,m}}w>\zeta\right)
 \leq\frac{e^{cT}}\zeta\ip\mu{w\mathbf1_{B_n^c}},\qquad m>n.
\end{equation}

Choose $n_k$ with $\ip\mu{w\mathbf1_{B_{n_k}^c}}\leq2^{-3k}$.
After $m\to\infty$, the sum of the probabilities in
\eqref{eq:nc-SBM-weighted-tail} for $(n,\zeta)=(n_k,2^{-k})$
is at most $e^{cT}\sum_{k=1}^\infty2^{-2k}<\infty$.
The Borel--Cantelli lemma proves uniform convergence in weighted mass
on $[0,T]$, almost surely, and hence continuity of the limiting
process $X^\kappa$ in the vague topology.
Monotone convergence in \eqref{eq:nc-SBM-first-moment}, applied to
$X^{\kappa,[n]}$, gives, for every nonnegative compactly supported
Borel function $\varphi$,
\begin{equation}\label{eq:nc-SBM-local-mean}
 \begin{aligned}
 \E\left[\ip{X_t^\kappa}{\varphi}\right]
 &=\lim_{n\to\infty}\ip{\mu|_{B_n}}{P_t^{D,\kappa}\varphi}
 =\ip{\mu}{P_t^{D,\kappa}\varphi}\\
 &\leq m_1\|\ell\|_\infty e^{-\kappa t/m_1}\|\varphi\|_1.
 \end{aligned}
\end{equation}
For $\psi\in C_c^2(\Rd)$, Doob's $L^2$ maximal inequality and
\eqref{eq:nc-SBM-first-moment} give
\begin{align}
 \E\left[\sup_{0\leq t\leq T}\left|M_t^{n,m}(\psi)\right|^2\right]
 &\leq\frac4{m_1}\int_0^T
             \ip{\mu|_{B_n^c}}{P_s^{D,\kappa}\psi^2}\dd s
 \leq C_{\psi,T}\ip\mu{w\mathbf1_{B_n^c}}\longrightarrow0,\nonumber\\
 \E\left[\int_0^T
        \ip{Y_s^{n,m}}{\left|(\cA-\kappa/m_1)\psi\right|}\dd s\right]
 &\leq C_{\psi,T}\ip\mu{w\mathbf1_{B_n^c}}\longrightarrow0.
 \label{eq:nc-SBM-martingale-tail}
\end{align}
Here $\psi^2+|(\cA-\kappa/m_1)\psi|\leq C_\psi w$ and
$P_s^{D,\kappa}w\leq e^{cs}w$.
Let $\cM^{\kappa,[n]}(\psi)$ be the martingale in
\eqref{eq:nc-SBM-martingale} with $X^{\kappa,[n]}$ and initial measure
$\mu|_{B_n}$.  Its predictable quadratic variation, given by the
finite-measure version of \eqref{eq:nc-SBM-bracket}, is
\begin{equation}\label{eq:nc-SBM-finite-bracket}
 B_t^{[n]}(\psi):=\left\langle\cM^{\kappa,[n]}(\psi)\right\rangle_t
       =\frac1{m_1}\int_0^t\ip{X_s^{\kappa,[n]}}{\psi^2}\dd s.
\end{equation}
Set $B_t(\psi):=m_1^{-1}\int_0^t\ip{X_s^\kappa}{\psi^2}\dd s$.
The increasing construction, \eqref{eq:nc-SBM-first-moment}, and
the same weighted estimate used in \eqref{eq:nc-SBM-martingale-tail}
give
\begin{align}
 \E\left[\sup_{0\leq t\leq T}|B_t(\psi)-B_t^{[n]}(\psi)|\right]
 &=\frac1{m_1}\int_0^T
       \E\left[\ip{X_s^\kappa-X_s^{\kappa,[n]}}{\psi^2}\right]\dd s
       \nonumber\\
 &=\frac1{m_1}\int_0^T
                    \ip{\mu|_{B_n^c}}{P_s^{D,\kappa}\psi^2}\dd s
 \leq C_{\psi,T}\ip\mu{w\mathbf1_{B_n^c}}\longrightarrow0.
 \label{eq:nc-SBM-bracket-convergence}
\end{align}
Thus these predictable quadratic variations converge in $L^1$,
uniformly on $[0,T]$.  Together with the $L^2$ martingale convergence
in \eqref{eq:nc-SBM-martingale-tail}, this permits passage to the
limit in the martingale and squared-martingale identities and proves
\eqref{eq:nc-SBM-martingale}--\eqref{eq:nc-SBM-bracket}.

Independence of the exhaustion follows explicitly from the finite-time
transition log-Laplace operators $\mathcal U_t^\kappa$:
for $0<t_1<\cdots<t_k$ and $\varphi_i\in C_c^+(\Rd)$, let
\[
 V_k:=\varphi_k,\qquad
 V_i:=\varphi_i+\mathcal U_{t_{i+1}-t_i}^\kappa V_{i+1},\quad (1\leq i<k),
 \qquad V:=\mathcal U_{t_1}^\kappa V_1.
\]
By \cite[Theorem~II.5.11(c)]{perkins2002superprocess},
\[
 \E\left[e^{-\sum_{i=1}^k\ip{X_{t_i}^{\kappa,[n]}}{\varphi_i}}\right]
 =e^{-\ip{\mu|_{B_n}}V}\longrightarrow e^{-\ip\mu V},\qquad
 0\leq V\leq\sum_{i=1}^kP_{t_i}^{D,\kappa}\varphi_i.
\]
The majorant has $\mu$-integral at most
$m_1\|\ell\|_\infty\sum_{i=1}^k\|\varphi_i\|_1$.
The argument of \eqref{eq:SBM-exhaustion-Laplace}, including the
rational-time coordinates of continuous paths, now gives the same
process law for every finite-measure exhaustion.

For the occupation formula let $u:=u_f^\kappa$ and, for the finite
process $X^{\kappa,[n]}$, define
\[
 A_t^{[n]}:=\ip{X_t^{\kappa,[n]}}{u(t)}
       +\frac1{m_1}\int_0^t\ip{X_s^{\kappa,[n]}}{f(s)}\dd s.
\]
The regularity proved in \eqref{eq:log-Laplace-spatial-derivatives}
and the time-dependent formula
\cite[Proposition~II.5.7]{perkins2002superprocess} apply here.
Equation \eqref{eq:nc-log-Laplace} gives
\begin{align*}
 \dd A_t^{[n]}
 &=\ip{X_t^{\kappa,[n]}}{\partial_tu+\cA u-\kappa u/m_1+f/m_1}\dd t
                      +\dd L_t^{[n]}\\
 &=\frac1{2m_1}\ip{X_t^{\kappa,[n]}}{(u(t))^2}\dd t+\dd L_t^{[n]},
 \qquad
 \dd\left\langle L^{[n]}\right\rangle_t
 =\frac1{m_1}\ip{X_t^{\kappa,[n]}}{(u(t))^2}\dd t.
\end{align*}
It\^o's formula therefore yields
$\dd e^{-A_t^{[n]}}=-e^{-A_t^{[n]}}\dd L_t^{[n]}$.
Since $0\leq e^{-A_t^{[n]}}\leq1$, taking expectations at $0$ and $T$
proves
\[
 \E\left[e^{-m_1^{-1}\int_0^T\ip{X_t^{\kappa,[n]}}{f(t)}\dd t}\right]
       =e^{-\ip{\mu|_{B_n}}{u(0)}}.
\]
The linear majorant in the mild equation gives
\[
 0\leq u(t)\leq\frac1{m_1}\int_t^TP_{s-t}^{D,\kappa}f(s)\dd s,
 \qquad
 \ip\mu{u(0)}\leq\|\ell\|_\infty\int_0^T\|f(s)\|_1\dd s<\infty.
\]
Monotone convergence of the occupation integrals and bounded
convergence of their exponentials prove
\eqref{eq:nc-occupation-Laplace}.  Also, for bounded Borel $K$,
\[
 \E\left[\frac1{m_1}\int_0^TX_t^\kappa(K)\dd t\right]
 =\int_0^T\ip\ell{P_t^{D,\kappa}\mathbf1_K}\dd t
 \leq T\|\ell\|_\infty|K|,
\]
which proves local finiteness.

For the canonical measure, adjoin an isolated absorbing state $\Delta$
and kill Brownian motion at rate $\kappa/m_1$.  Its conservative
semigroup on $\Rd\cup\{\Delta\}$ is
\[
 \widetilde P_tf(x)=e^{-\kappa t/m_1}P_t^D(f|_{\Rd})(x)
        +\left(1-e^{-\kappa t/m_1}\right)f(\Delta),\quad x\in\Rd,
 \qquad \widetilde P_tf(\Delta)=f(\Delta).
\]
It is a Borel strong Markov process preserving bounded continuous
functions.  Apply \cite[Theorem~II.7.3]{perkins2002superprocess}
with branching variance $1/m_1$ and zero linear branching drift,
and project its canonical measures onto paths restricted to $\Rd$.
Denote these measures by $\mathbb N_x^\kappa$.  For a test function extended
by zero at $\Delta$, the motion generator is
$\cA-\kappa/m_1$, so the projected processes have exactly
\eqref{eq:nc-SBM-martingale}--\eqref{eq:nc-SBM-bracket}.
Let $\mathcal O(X)(\dd t,\dd y):=m_1^{-1}X_t(\dd y)\dd t$ and set
\[
 \mathbb Q_{\lambda_0}^\kappa(B)
 :=m_1\int_{\Rd}\lambda_0(x)
               \mathbb N_x^\kappa(\mathcal O(X)\in B)\dd x,
 \qquad B\subset\cM_{\mathrm{loc}}(E)\setminus\{0\}.
\]
The canonical Poisson decomposition and the finite-initial-measure
occupation formula yield
\begin{equation}\label{eq:nc-canonical-pointwise}
 \mathbb N_x^\kappa\left(1-e^{-\ip{\mathcal O(X)}g}\right)=u_g^\kappa(0,x).
\end{equation}
Indeed both sides are the negative logarithm of the occupation
Laplace functional of the process started at $\delta_x$.
Integrating \eqref{eq:nc-canonical-pointwise} against
$m_1\lambda_0(x)\dd x$ proves the first identity in
\eqref{eq:nc-canonical-identification}.

For the first moment let
$v_g(t):=m_1^{-1}\int_t^TP_{s-t}^{D,\kappa}g(s)\dd s$.
The mild equation \eqref{eq:nc-log-Laplace-mild} gives, for
$0<\theta\leq1$, $ 0\leq u_{\theta g}^\kappa\leq\theta v_g$, and
\begin{align*}
 0\leq v_g(t)-\theta^{-1}u_{\theta g}^\kappa(t)
 =\frac1{2m_1\theta}\int_t^TP_{s-t}^{D,\kappa}
                     \left(u_{\theta g}^\kappa(s)\right)^2\dd s
 \leq\frac{\theta}{2m_1}\int_t^TP_{s-t}^{D,\kappa}\left(v_g(s)\right)^2\dd s.
\end{align*}
The last bound tends to zero uniformly.  Monotone convergence theorem
applied to $\frac{1-e^{-\theta\ip{\mathcal O(X)}g}}{\theta}$ in
\eqref{eq:nc-canonical-pointwise} thus gives
$\mathbb N_x^\kappa[\ip{\mathcal O(X)}g]=v_g(0,x)$.
The semigroup $P_t^{D,\kappa}$ in \eqref{eq:nc-heat-profile} has the
symmetric Gaussian kernel in \eqref{eq:heat-kernel} multiplied by
$e^{-\kappa t/m_1}$.  Tonelli's theorem and this symmetry imply
\begin{align*}
 \int_{\cM_{\mathrm{loc}}(E)\setminus\{0\}}\ip\nu g
                                  \mathbb Q_{\lambda_0}^\kappa(\dd\nu)
 &=m_1\int_{\Rd}\lambda_0(x)v_g(0,x)\dd x\\
 &=\int_0^T\int_{\Rd}g(t,x)\left(P_t^{D,\kappa}\lambda_0\right)(x)\dd x\dd t
 \leq\|\lambda_0\|_\infty\|g\|_1.
\end{align*}
This is the second identity in \eqref{eq:nc-canonical-identification}.
The identities extend from smooth to continuous compactly supported
test functions by approximation, using this mean bound and
$|e^{-a}-e^{-b}|\leq|a-b|$ for $a,b\geq0$.
Finally, choose $\chi_j\in C_c^+(E)$ equal to one on a compact
exhaustion.  The sets
$B_{j,k}=\left\{\nu:\ip\nu{\chi_j}>1/k\right\}$, $j,k\geq1$, satisfy
\[
 \cM_{\mathrm{loc}}(E)\setminus\{0\}
       =\bigcup_{j,k\geq1}B_{j,k},\qquad
 \mathbb Q_{\lambda_0}^\kappa(B_{j,k})
       \leq k\|\lambda_0\|_\infty\|\chi_j\|_1<\infty.
\]
This proves sigma-finiteness.
The proof is complete.
\end{proof}

\subsubsection{One-dimensional representation and shorter observation windows}

When $d=1$, the finite-measure density result used in Lemma
\ref{lem:locally-finite-SBM} remains valid with constant linear decay.
Indeed, after stopping when the total mass exceeds a fixed level,
Dawson's Girsanov theorem
\cite[Theorem~IV.1.6]{perkins2002superprocess} transfers the density
property from zero drift to the bounded drift $-\kappa/m_1$ on each
finite time interval.  In fact, if
$\tau_L=\inf\{t:X_t(1)\geq L\}$, the quadratic variation of the
Girsanov exponent up to $T\wedge\tau_L$ is at most
$\kappa^2LT/m_1$, so the stopped exponential satisfies Novikov's
condition.  Nonexplosion of the finite-measure total mass lets
$L\to\infty$.  Removing the stopping and summing the independent
finite-measure components gives a measurable density
$v^\kappa=m_1\rho$ for $X^\kappa$.  The same orthogonal-martingale-measure
construction as in Lemma~\ref{lem:locally-finite-SBM} provides a white noise
$W$.  The martingale problem
\eqref{eq:nc-SBM-martingale}--\eqref{eq:nc-SBM-bracket} gives the
spatially weak formulation: for each $\psi\in C_c^2(\R)$, almost surely for all
$t\geq0$,
\begin{align*}
 m_1\ip{\rho_t}{\psi}
 &=m_1\ip{\lambda_0}{\psi}
   +\int_0^t\ip{\rho_s}{\frac{m_2}{2}\psi''-\kappa\psi}\dd s\\
 &\qquad+\int_0^t\int_{\R}
       \psi(x)\sqrt{\rho(s,x)}\,W(\dd s,\dd x).
\end{align*}
The first-moment bound \eqref{eq:nc-SBM-local-mean} gives
$\E[\int_0^t\int_{\R}\psi(x)^2\rho(s,x)\dd x\dd s]
\leq t\|\lambda_0\|_\infty\|\psi\|_{L^2(\R)}^2$, so the white-noise
integral is square integrable.  The initial condition is again the
vague limit $\rho_t(x)\dd x\to\lambda_0(x)\dd x$ as $t\downarrow0$,
with the white-noise and initial-value conventions used for
\eqref{eq:candidate-spde}.
In stochastic heat equation notation, the test-function identity reads
\begin{equation}\label{eq:nc-spde}
  m_1\dd\rho_t
       =\left(\frac{m_2}{2}\partial_{xx}\rho_t
                       -\kappa\rho_t\right)\dd t
          +\sqrt{\rho_t}\,W(\dd t,\dd x),
       \qquad\rho_0=\lambda_0.
\end{equation}
The measure-valued formulation remains the assertion in general dimension.

If the initial profile is also integrable, the total mass
$Z_t=\ip{X_t^\kappa}{1}$ satisfies
\begin{equation}\label{eq:nc-total-mass}
  \dd Z_t=-\frac{\kappa}{m_1}Z_t\dd t
             +\sqrt{\frac{Z_t}{m_1}}\dd B_t,
  \qquad Z_0=m_1\int_{\Rd}\lambda_0(x)\dd x.
\end{equation}
This follows from the finite-measure martingale problem with constant
test function one.  Thus the total mass is a Feller branching diffusion
with linear decay.  The positive immigration term in the CIR limit of
\cite[Theorem~2.2]{jaisson2015limit} comes from its positive baseline
intensity; it is absent in the prescribed-past model studied here.
The heavy-tailed setting of \cite{jaisson2016rough}, in which the mean
reproduction delay is infinite, is treated under separate assumptions in
Section~\ref{subsec:nc-heavy}.

On a spatial observation scale $\epsilon^{-\gamma}$, the relevant decay
parameter is $(1-a_\epsilon)/\epsilon^{2\gamma}$.  Thus the same
microscopic perturbation becomes negligible on a shorter parabolic
window.  In particular, the one-dimensional local critical limit from
Section~\ref{subsec:critical-results} persists.

\begin{corollary}[Local limit for nearly critical kernels]
\label{cor:nc-local}
Assume $d=1$, \eqref{eq:H1}--\eqref{eq:H4}, \eqref{eq:nc-kernel}, and
$0<\alpha<1$.  Define $\widehat{\cN}_\epsilon^{a,\alpha,x_\star}$ by
\eqref{eq:local-alpha-measure} with $N^{\epsilon,a}$ in place of
$N^\epsilon$.  Then, for every $x_\star\in\R$,
\begin{equation}\label{eq:nc-local-limit}
  \widehat{\cN}_\epsilon^{a,\alpha,x_\star}
       \xRightarrow[\epsilon\downarrow0]{}\frac1{m_1}X_t^{x_\star}(\dd x)\dd t,
\end{equation}
where $X^{x_\star}$ is the critical process of Theorem
\ref{thm:local-critical}, with initial measure
$m_1\lambda_0(x_\star)\dd x$ and zero linear decay.
\end{corollary}

\begin{proof}
Let $\delta=\epsilon^\alpha$ and translate by $x_\star/\epsilon$ as in
Theorem~\ref{thm:local-critical}.  The profile in $\delta$ coordinates is
\eqref{eq:frozen-profile}, and
\[
  \frac{1-a_\epsilon}{\delta^2}
       =\kappa_\epsilon\epsilon^{2(1-\alpha)}\xrightarrow[\epsilon\downarrow0]{}0.
\]
Applying the varying-profile assertion of Theorem~\ref{thm:nearly-critical}
with small parameter $\delta$, scalar coefficients
$a_{\delta^{1/\alpha}}$, and limiting decay parameter zero completes the proof.
\end{proof}

\subsection{Heavy-tail case}
\label{subsec:nc-heavy}

An infinite mean reproduction delay changes both the clock and the
contribution of the prescribed past.  We prove a space--time random-measure limit for a
factorized kernel with a regularly varying temporal tail.  The tail
index below is denoted by $\eta$ to distinguish it from the density
parameter $\alpha$.

\subsubsection{Assumptions, scaling, and the surviving past}

In this subsection replace \eqref{eq:H3} by the factorization and tail assumption
\begin{equation}\label{eq:ht-kernel}
  h(s,z)=H(s)p(z),\qquad
  \overline H(s):=\int_s^\infty H(r)\dd r
       \sim c_Hs^{-\eta},\qquad 0<\eta<1,\quad c_H>0.
\end{equation}
Here $H$ and $p$ are probability densities on $\R_+$ and $\Rd$,
respectively, and
\begin{equation}\label{eq:ht-spatial-moments}
  \int_{\Rd}zp(z)\dd z=0,\qquad
  \Sigma:=\int_{\Rd}zz^{\top}p(z)\dd z
       \quad\text{ is finite and positive definite}.
\end{equation}
The product structure separates spatial averaging from the temporal
tail in the functional of the Poisson past below.  The renewal limit itself
also holds for fixed dependent increment laws with finite spatial
second moment: \eqref{eq:clock-mixed-remainder} makes the mixed term
$o(\delta_\epsilon)$ and leaves the limiting symbol additive.
Here $m_1=\infty$, and the scaling constants are
\begin{equation}\label{eq:ht-constants}
  c_\eta:=c_H\Gamma(1-\eta),\qquad
  \delta_\epsilon:=\epsilon^{2\eta},\qquad
  \frac{1-a_\epsilon}{\delta_\epsilon}\xrightarrow[\epsilon\downarrow0]{}\kappa\geq0,
  \qquad 0<a_\epsilon\leq1.
\end{equation}
This includes exact criticality $a_\epsilon=1$.  For $\kappa>0$,
writing the observation horizon as $T=\epsilon^{-2}$ gives
$T^\eta(1-a_{T^{-1/2}})\to\kappa$ as $T\to\infty$, the positive balance in
\cite[Assumption~2]{jaisson2016rough}.

For $\alpha>0$ and bounded uniformly continuous
$\lambda_0:\Rd\to\R_+$, prescribe a Poisson past with intensity
\begin{equation}\label{eq:ht-past}
  q_\epsilon\lambda_0(\epsilon^\eta y)\dd s\dd y,
  \qquad s<0,\qquad
  q_\epsilon:=\epsilon^{\,2-\eta(\alpha+3-d)}.
\end{equation}
Let $N^{\epsilon,\eta}$ be its positive-time descendant process with
kernel $a_\epsilon H(s)p(z)$, as in \eqref{eq:nc-intensity}, and define
\begin{equation}\label{eq:ht-measure}
  \ip{\cN_\epsilon^{\eta,\alpha}}{g}
       :=\epsilon^{\eta(\alpha+3)}
          \int_E g(\epsilon^2t,\epsilon^\eta x)
                   N^{\epsilon,\eta}(\dd t,\dd x).
\end{equation}
The physical spatial window is $T^{\eta/2}$.  At $\alpha=1$, each
rescaled event has mass $T^{-2\eta}$ and the microscopic past
intensity is $T^{\eta(4-d)/2-1}\lambda_0(x/T^{\eta/2})$.
The exponents in \eqref{eq:ht-past}--\eqref{eq:ht-measure} reduce to
those of the finite-mean normalization when formally setting $\eta=1$.
For $0<\eta<1$ they are a different family of prescribed pasts and observation scales;
keeping \eqref{eq:past-intensity} unchanged would not preserve the same
fluctuation threshold under this new clock.

The generation-by-generation construction in Proposition
\ref{prop:construction} only uses that $H$ and $p$ have mass one.
It therefore still gives a locally finite process, uniqueness in the
bounded mean-intensity class, and
\begin{equation}\label{eq:ht-mean-bound}
  \E\left[N^{\epsilon,\eta}(\dd t,\dd x)\right]
       \leq q_\epsilon\|\lambda_0\|_\infty\dd t\dd x,
  \qquad
  \E\left[\ip{\cN_\epsilon^{\eta,\alpha}}{g}\right]
       \leq\|\lambda_0\|_\infty
          \int_E g(t,x)\dd t\dd x
       \quad(g\in C_c^+(E)).
\end{equation}
If $M_\epsilon^{\eta,\alpha}(g)$ is the compensated integral with the
normalization in \eqref{eq:ht-measure}, its variance satisfies
\begin{equation}\label{eq:ht-bracket}
  \E\left[\left|\delta_\epsilon^{-1}
                   M_\epsilon^{\eta,\alpha}(g)\right|^2\right]
       \leq\epsilon^{\eta(\alpha-1)}
          \|\lambda_0\|_\infty
          \int_E (g(t,x))^2\dd t\dd x.
\end{equation}
Thus near-criticality does not merge the three density regimes.
The branching deficit determines the drift scale, while $\alpha$
still determines whether the normalized noise disappears, survives,
or requires a sparse-limit analysis.

\subsubsection{The fractional operator and the limiting law}

Write $\mathcal B\psi:=\frac12\Sigma:D_x^2\psi$, and denote its
heat semigroup by $P_r^\Sigma$.  For $r>0$, its Gaussian density and
transition measure are
\begin{equation}\label{eq:ht-brownian-kernel}
 p_r^\Sigma(y):=
 \frac{\exp\!\left[-\dfrac{y^\top\Sigma^{-1}y}{2r}\right]}
      {(2\pi r)^{d/2}\sqrt{\det\Sigma}},
 \qquad p_r^\Sigma(\dd y):=p_r^\Sigma(y)\dd y.
\end{equation}
Thus $P_r^\Sigma f(x)=\int_{\Rd}f(x+y)p_r^\Sigma(\dd y)$; at
$r=0$, let $p_0^\Sigma(\dd y)=\delta_0(\dd y)$.
Let $S_r$ be an $\eta$-stable
subordinator and $B_r^\Sigma$ an independent Brownian motion, normalized
by
\[
  \E\left[e^{-zS_r}\right]=e^{-r c_\eta z^\eta},\qquad
  \E\left[B_r^\Sigma\left(B_r^\Sigma\right)^{\top}\right]=r\Sigma.
\]
For bounded functions on $[0,\infty)\times\Rd$, define
\begin{equation}\label{eq:ht-resolvent}
  \mathsf P_r^\eta f(t,x):=\E\left[f\left(t+S_r,x+B_r^\Sigma\right)\right],\qquad
  R_b^\eta f:=\int_0^\infty e^{-br}\mathsf P_r^\eta f\dd r,
       \qquad b>0.
\end{equation}
For $b=0$, the same integral defines $R_0^\eta$ on the finite-horizon
space $\mathcal X_T$ of Lemma~\ref{lem:renewal-framework}: for $q>0$
and $p=1,\infty$, its norm in $\|\cdot\|_{q,p}$ is at most
$(c_\eta q^\eta)^{-1}$.  Here $r$ is a genealogical clock.  The generator of $\mathsf P^\eta$ is
\begin{equation}\label{eq:ht-generator}
  \mathcal G_\eta\varphi
       :=\mathcal B\varphi-c_\eta D_+^\eta\varphi,
  \qquad
  D_+^\eta\varphi(t,x)
       :=\frac{\eta}{\Gamma(1-\eta)}
         \int_0^\infty
          \frac{\varphi(t,x)-\varphi(t+s,x)}{s^{1+\eta}}\dd s.
\end{equation}
For terminal-value problems the test function is extended by zero beyond the
terminal time; this convention includes the whole future in $D_+^\eta$.

For $g\in C_c^+(E)$, choose $T$ beyond its temporal support and let
$u_g^{\eta,\kappa}$ be the unique bounded nonnegative solution in
$\mathcal X_T$ of
\begin{equation}\label{eq:ht-nonlinear-mild}
  u=R_\kappa^\eta\left(g-\frac12u^2\right).
\end{equation}
Its existence and uniqueness are proved below.  The solution vanishes
after the temporal support of $g$ and is independent of the choice of $T$.  This is the mild interpretation of
\begin{equation}\label{eq:ht-nonlinear-equation}
  c_\eta D_+^\eta u-\mathcal B u+\kappa u+\frac12u^2=g.
\end{equation}
Define
\begin{equation}\label{eq:ht-exponent}
  \Psi_{\lambda_0}^{\eta,\kappa}(g)
       :=c_H\int_0^\infty t^{-\eta}
          \int_{\Rd}\lambda_0(x)u_g^{\eta,\kappa}(t,x)\dd x\dd t.
\end{equation}
The singularity at zero is integrable since $\eta<1$.  Spatial
integrability will follow from the bounds used in the proof.
In contrast to \eqref{eq:nc-canonical-exponent}, the entire function
$u_g^{\eta,\kappa}(t,\cdot)$ enters this exponent.  The infinite-mean
past contributes the temporal density $c_Ht^{-\eta}$ instead of a
finite mass concentrated at $t=0$.

The mean profile has a separate inverse-clock representation.  With
$E_t:=\inf\{r:S_r>t\}$, let
\begin{equation}\label{eq:ht-mean-profile}
  \rho^{\eta,\kappa}(t,\cdot)
       :=\E\left[e^{-\kappa E_t}P_{E_t}^\Sigma\lambda_0\right].
\end{equation}
It is the bounded mild solution of the Caputo equation
\begin{equation}\label{eq:ht-caputo}
  c_\eta\,{}^{\mathrm C}D_t^\eta\rho^{\eta,\kappa}
       =\mathcal B\rho^{\eta,\kappa}
           -\kappa\rho^{\eta,\kappa},
       \qquad \rho^{\eta,\kappa}(0,\cdot)=\lambda_0,
\end{equation}
where for an absolutely continuous scalar function
\[
  {}^{\mathrm C}D_t^\eta v(t)
       :=\frac1{\Gamma(1-\eta)}
            \int_0^t(t-s)^{-\eta}v'(s)\dd s.
\]
Equation \eqref{eq:ht-caputo} is interpreted through
\eqref{eq:ht-mean-profile} for a general bounded uniformly continuous
initial profile.

With heavy-tailed delays, the prescribed past contributes throughout
positive physical time. The following theorem gives the
three density regimes under this new clock, including exact
criticality.

\begin{theorem}[Heavy-tail critical and nearly critical limits]\label{thm:ht-limits}
Assume \eqref{eq:ht-kernel}--\eqref{eq:ht-constants}, with
$\kappa\geq0$, and use \eqref{eq:ht-past}--\eqref{eq:ht-measure}.
For bounded uniformly continuous $\lambda_0\geq0$:
\begin{enumerate}
  \item for every $\alpha>0$,
  \begin{equation}\label{eq:ht-annealed}
    \E\left[\cN_\epsilon^{\eta,\alpha}\right]\xrightarrow[\epsilon\downarrow0]{v}
        \rho^{\eta,\kappa}(t,x)\dd t\dd x;
  \end{equation}
  \item if $\alpha>1$, then
  \begin{equation}\label{eq:ht-hydrodynamic}
    \cN_\epsilon^{\eta,\alpha}\xrightarrow[\epsilon\downarrow0]{\mathbb P}
        \rho^{\eta,\kappa}(t,x)\dd t\dd x;
  \end{equation}
  \item if $\alpha=1$, there is a locally finite random measure
  $\nu_{\lambda_0}^{\eta,\kappa}$ such that
  \begin{equation}\label{eq:ht-stochastic}
    \cN_\epsilon^{\eta,1}\xRightarrow[\epsilon\downarrow0]{}
          \nu_{\lambda_0}^{\eta,\kappa},\qquad
    \E\left[e^{-\ip{\nu_{\lambda_0}^{\eta,\kappa}}{g}}\right]
          =e^{-\Psi_{\lambda_0}^{\eta,\kappa}(g)}
          \quad(g\in C_c^+(E));
  \end{equation}
  this Laplace functional determines its law;
  \item if $0<\alpha<1$, let $r_\epsilon=\epsilon^{\eta(1-\alpha)}$.
  Then
  \begin{equation}\label{eq:ht-sparse}
    \cN_\epsilon^{\eta,\alpha}\xrightarrow[\epsilon\downarrow0]{\mathbb P}0,
    \qquad
    -\frac1{r_\epsilon}
      \log\E\left[e^{-\ip{r_\epsilon\cN_\epsilon^{\eta,\alpha}}{g}}\right]
        \xrightarrow[\epsilon\downarrow0]{}\Psi_{\lambda_0}^{\eta,\kappa}(g)
        \quad(g\in C_c^+(E)).
  \end{equation}
\end{enumerate}
All random-measure convergences use the vague topology.  The result
holds in every fixed $d\geq1$ and for every $0<\eta<1$ at the level
of event measures.  Proposition~\ref{prop:ht-lifted} gives a
superprocess occupation representation of the stochastic limit.
\end{theorem}

\begin{proof}
Fix $g\in C_c^+(E)$ and choose $T<\infty$ beyond its temporal
support.  The mean bound \eqref{eq:ht-mean-bound} gives tightness
by the compact-exhaustion argument of Proposition~\ref{prop:tightness}.
Lemma~\ref{lem:ht-linear-mean} proves \eqref{eq:ht-annealed}.

For $\alpha=1$, Lemmas~\ref{lem:ht-cluster-costs},
\ref{lem:ht-preparation-functional}, and
\ref{lem:ht-critical-exponent} give
\[
 \log\E\left[e^{-\ip{\cN_\epsilon^{\eta,1}}{g}}\right]
       \xrightarrow[\epsilon\downarrow0]{}-\Psi_{\lambda_0}^{\eta,\kappa}(g).
\]
Tightness provides a subsequential limit, whose Laplace functional is
$\exp\{-\Psi_{\lambda_0}^{\eta,\kappa}(g)\}$ for every $g\in C_c^+(E)$.
Lemma~\ref{lem:Laplace-vague} shows that all subsequential limits have
this same law.  This proves both existence of
$\nu_{\lambda_0}^{\eta,\kappa}$ and \eqref{eq:ht-stochastic}.
For $Y=\ip{\nu_{\lambda_0}^{\eta,\kappa}}g$, monotone convergence theorem
and \eqref{eq:ht-exponent-derivative} give the explicit first moment
\begin{align}
 \E[Y]
 &=\lim_{\theta\downarrow0}\frac{1-\E[e^{-\theta Y}]}\theta
 =\lim_{\theta\downarrow0}
       \frac{1-e^{-\Psi_{\lambda_0}^{\eta,\kappa}(\theta g)}}\theta
       \nonumber\\
 &=\lim_{\theta\downarrow0}
       \frac{\Psi_{\lambda_0}^{\eta,\kappa}(\theta g)}\theta
 =\int_Eg(t,x)\rho^{\eta,\kappa}(t,x)\dd t\dd x.
 \label{eq:ht-limit-first-moment}
\end{align}
The third equality uses
$\Psi_{\lambda_0}^{\eta,\kappa}(\theta g)=\mathcal O(\theta)$ and
$(1-e^{-z})/z\to1$ as $z\downarrow0$.

For $\alpha>1$, let $r_\epsilon=\epsilon^{\eta(1-\alpha)}\to\infty$.
The exact Laplace formula for the Poisson past uses profile $r_\epsilon\lambda_0$
and source $g/r_\epsilon$.
Lemma~\ref{lem:ht-varying-sources} replaces
$r_\epsilon u_{\epsilon,g/r_\epsilon}$ by
$R_\epsilon^{(0)}g$ with an error tending to zero, including the
parent linearization error.
Lemmas~\ref{lem:ht-preparation-functional} and
\ref{lem:ht-linear-mean} then give
\[
 \log\E\left[e^{-\ip{\cN_\epsilon^{\eta,\alpha}}{g}}\right]
    \xrightarrow[\epsilon\downarrow0]{}
        -\int_E g(t,x)\rho^{\eta,\kappa}(t,x)\dd t\dd x.
\]
Tightness and the deterministic limiting law prove
\eqref{eq:ht-hydrodynamic}.

For $0<\alpha<1$, the observable
$r_\epsilon\cN_\epsilon^{\eta,\alpha}$ has source $g$, with the
initial profile multiplied by $r_\epsilon$.
The exact Poisson formula \eqref{eq:ht-past-Laplace} therefore gives
\[
 -\frac1{r_\epsilon}
  \log\E\left[e^{-\ip{r_\epsilon\cN_\epsilon^{\eta,\alpha}}{g}}\right]
 =-\log\E\left[e^{-\ip{\cN_\epsilon^{\eta,1}}{g}}\right],
\]
which converges to $\Psi_{\lambda_0}^{\eta,\kappa}(g)$ by
Lemma~\ref{lem:ht-critical-exponent}.
For the unmultiplied observable, the bound for the Poisson past and
Lemma~\ref{lem:ht-varying-sources}, with $\theta=r_\epsilon^{-1}$,
yield
\begin{equation}\label{eq:ht-sparse-zero-bound}
 0\leq-\log\E\left[e^{-\ip{\cN_\epsilon^{\eta,\alpha}}{g}}\right]
        \leq C_{g,T}r_\epsilon^{1/(d+2)}\to0.
\end{equation}
Thus the Laplace functionals converge to $L_0(g)=1$ for every
$g\in C_c^+(E)$.  The deterministic-limit assertion of
Lemma~\ref{lem:Laplace-vague} proves the remaining assertion
of \eqref{eq:ht-sparse}.
The proof is complete.
\end{proof}

\subsubsection{Resolvent convergence and the nonlinear cluster equation}\label{sec:heavy:support:1}

The tail assumption fixes the clock exponent in the common
renewal lemma.  We first record the resulting linear estimates.

The forward offspring operator in this scaling is
\begin{align}\label{offspring:operator}
  \mathsf K_\epsilon\varphi(t,x)
       :=\int_0^\infty\int_{\Rd}H(s)p(z)
            \varphi\left(t+\epsilon^2s,x+\epsilon^\eta z\right)\dd z\dd s.
\end{align}

For $b\geq0$, define
\begin{equation}\label{eq:ht-discrete-resolvent}
 R_\epsilon^{(b)}:=\mathcal R_{\epsilon,b}
  =\frac{\delta_\epsilon}{1+b\delta_\epsilon}
    \sum_{n\geq0}
       \left(\frac{a_\epsilon}{1+b\delta_\epsilon}\right)^n
       \mathsf K_\epsilon^n,\qquad b\geq0.
\end{equation}
\begin{lemma}[Heavy-tail operator and resolvent estimates]
\label{lem:ht-resolvent-estimates}
Assume \eqref{eq:ht-kernel}--\eqref{eq:ht-constants}.
The temporal tail and Laplace transform satisfy
\begin{equation}\label{eq:ht-tail-bounds}
  \overline H(s)\leq C\min\left(1,s^{-\eta}\right),\qquad
  1-\widehat H(z)\sim c_\eta z^\eta\quad(z\downarrow0),
\end{equation}
where $\widehat H(z):=\int_0^\infty e^{-zs}H(s)\dd s$.
The offspring operator \eqref{offspring:operator} has the expansion
\begin{equation}\label{eq:ht-operator-limit}
  \frac{a_\epsilon\mathsf K_\epsilon\varphi-\varphi}
       {\delta_\epsilon}
       \xrightarrow[\epsilon\downarrow0]{}(\mathcal G_\eta-\kappa)\varphi
       \quad\text{in }L^1(E),
       \qquad\varphi\in C_c^{1,2}([0,\infty)\times\Rd).
\end{equation}
For each $T<\infty$ and $b\geq0$, the resolvents satisfy
\begin{equation}\label{eq:ht-resolvent-strong}
 \left\|R_\epsilon^{(b)}v-R_{\kappa+b}^\eta v\right\|_\infty\to0
       \qquad(v\in\mathcal X_T).
\end{equation}
In the weighted norms of Lemma~\ref{lem:renewal-framework},
for $q>0$ and $p=1,\infty$,
\begin{equation}\label{eq:ht-finite-horizon-bound}
 \left\|R_\epsilon^{(b)}\right\|_{q,p}
 \leq\frac{\delta_\epsilon}
       {1+b\delta_\epsilon-a_\epsilon\widehat H(q\epsilon^2)}
       \xrightarrow[\epsilon\downarrow0]{}\frac1{\kappa+b+c_\eta q^\eta},
 \qquad
 \left\|R_{\kappa+b}^\eta\right\|_{q,p}
       \leq\frac1{\kappa+b+c_\eta q^\eta}.
\end{equation}
These statements include $b=\kappa=0$.
\end{lemma}

\begin{proof}
The first bound in \eqref{eq:ht-tail-bounds} follows from the tail
asymptotic and $\overline H\leq1$.  For the second statement in \eqref{eq:ht-tail-bounds}, we use
$1-\widehat H(z)=z\int_0^\infty e^{-zs}\overline H(s)\dd s$,
change variables, and apply the dominated convergence theorem using the first bound.
The same tail computation and the centered spatial Taylor expansion give
\eqref{eq:ht-operator-limit}.
Indeed, the temporal difference quotient is
\begin{align}
 T_\epsilon\varphi(t,x)
 &:=\frac1{\delta_\epsilon}\int_0^\infty H(r)
             \left[\varphi\left(t+\epsilon^2r,x\right)-\varphi(t,x)\right]\dd r\nonumber\\
 &=\frac1{\delta_\epsilon}\int_0^\infty H(r)
                   \int_0^{\epsilon^2r}\partial_t\varphi(t+s,x)\dd s\dd r
       \nonumber\\
 &=\int_0^\infty\frac{\overline H(s/\epsilon^2)}{\delta_\epsilon}
                              \partial_t\varphi(t+s,x)\dd s.
 \label{eq:ht-temporal-quotient}
\end{align}
Fubini's theorem applies because $\partial_t\varphi$ is bounded
and compactly supported and the tail bound in
\eqref{eq:ht-tail-bounds} is integrable near zero after rescaling.
Choose $T<\infty$ beyond the temporal support of $\varphi$.
By \eqref{eq:ht-kernel} and \eqref{eq:ht-tail-bounds},
$\overline H(s/\epsilon^2)/\delta_\epsilon\to c_Hs^{-\eta}$ for
$s>0$, and this ratio is at most $Cs^{-\eta}$.
For $0<a<b<\infty$, integration by parts gives
\begin{align}
 c_H\int_a^b s^{-\eta}\partial_t\varphi(t+s,x)\dd s
 &=c_H\left[s^{-\eta}
       \left(\varphi(t+s,x)-\varphi(t,x)\right)\right]_{s=a}^{s=b}
       \nonumber\\
 &\qquad\qquad+\eta c_H\int_a^b
       \frac{\varphi(t+s,x)-\varphi(t,x)}{s^{1+\eta}}\dd s.
 \label{eq:ht-temporal-integration-by-parts}
\end{align}
The two endpoint terms in \eqref{eq:ht-temporal-integration-by-parts}
satisfy
\begin{align*}
 \left|a^{-\eta}\left(\varphi(t+a,x)-\varphi(t,x)\right)\right|
 &\leq a^{1-\eta}\|\partial_t\varphi\|_\infty
       \xrightarrow[a\downarrow0]{}0,\\
 \left|b^{-\eta}\left(\varphi(t+b,x)-\varphi(t,x)\right)\right|
 &\leq2b^{-\eta}\|\varphi\|_\infty
       \xrightarrow[b\to\infty]{}0.
\end{align*}
Letting $a\downarrow0$ and $b\to\infty$ in
\eqref{eq:ht-temporal-integration-by-parts}, and using
\eqref{eq:ht-constants} and \eqref{eq:ht-generator}, gives
\begin{equation}\label{eq:ht-temporal-fractional-identity}
 c_H\int_0^\infty s^{-\eta}\partial_t\varphi(t+s,x)\dd s
 =-c_\eta D_+^\eta\varphi(t,x).
\end{equation}
Equations \eqref{eq:ht-temporal-quotient} and
\eqref{eq:ht-temporal-fractional-identity} therefore yield
\begin{align}
 \left\|T_\epsilon\varphi+c_\eta D_+^\eta\varphi\right\|_{L^1(E)}
 &\leq\int_0^T
   \left|\frac{\overline H(s/\epsilon^2)}{\delta_\epsilon}
                         -c_Hs^{-\eta}\right|
   \left\|\partial_t\varphi(\cdot+s,\cdot)\right\|_{L^1(E)}\dd s
 \longrightarrow0.
 \label{eq:ht-temporal-L1}
\end{align}
The dominated convergence theorem applies because the integrand in
\eqref{eq:ht-temporal-L1} is bounded by
\[
 (C+c_H)s^{-\eta}\|\partial_t\varphi\|_{L^1(E)},\qquad 0<s<T,
\]
which is integrable since $0<\eta<1$.

To separate the temporal and spatial terms in
\eqref{eq:ht-operator-limit}, define the averaging operators
\begin{equation}\label{eq:ht-separate-averages}
 \begin{aligned}
 (\mathsf H_\epsilon w)(t,x)
   &:=\int_0^\infty H(s)w(t+\epsilon^2s,x)\dd s,\\
 (\mathsf P_\epsilon w)(t,x)
   &:=\int_{\Rd}p(z)w(t,x+\epsilon^\eta z)\dd z.
 \end{aligned}
\end{equation}
The factorization \eqref{eq:ht-kernel} gives
$\mathsf K_\epsilon=\mathsf H_\epsilon\mathsf P_\epsilon$.
Let $V_\epsilon:=(\mathsf P_\epsilon\varphi-\varphi)/\delta_\epsilon$.
Equations \eqref{eq:ht-temporal-quotient} and
\eqref{eq:ht-separate-averages} give the exact decomposition
\begin{equation}\label{eq:ht-operator-decomposition}
 \frac{a_\epsilon\mathsf K_\epsilon\varphi-\varphi}{\delta_\epsilon}
 =a_\epsilon T_\epsilon\varphi
   +a_\epsilon\mathsf H_\epsilon V_\epsilon
   -\frac{1-a_\epsilon}{\delta_\epsilon}\varphi.
\end{equation}
We next estimate the spatial term $a_\epsilon\mathsf H_\epsilon V_\epsilon$
in \eqref{eq:ht-operator-decomposition}.
Taylor's formula with integral remainder and the centering condition
\eqref{eq:ht-spatial-moments} imply
\begin{align}
 V_\epsilon(t,x)
 &=\frac1{\delta_\epsilon}\int_{\Rd}p(z)
   \left[\varphi(t,x+\epsilon^\eta z)-\varphi(t,x)
                   -\epsilon^\eta z\cdot\nabla_x\varphi(t,x)\right]\dd z
       \nonumber\\
 &=\int_{\Rd}p(z)\int_0^1(1-\theta)
      z^{\top}D_x^2\varphi(t,x+\theta\epsilon^\eta z)z\dd\theta\dd z.
 \label{eq:ht-spatial-quotient}
\end{align}
Since $\delta_\epsilon=\epsilon^{2\eta}$ and
$\int_0^1(1-\theta)\dd\theta=1/2$, the covariance identity in
\eqref{eq:ht-spatial-moments} gives
\begin{align}
 \left\|V_\epsilon-\mathcal B\varphi\right\|_{L^1(E)}
 \leq\int_{\Rd}|z|^2p(z)\int_0^1(1-\theta)\cdot\left\|D_x^2\varphi(\cdot,\cdot+\theta\epsilon^\eta z)
                         -D_x^2\varphi\right\|_{L^1(E)}\dd\theta\dd z
 \longrightarrow0.
 \label{eq:ht-spatial-L1}
\end{align}
The norm of the Hessian uses the Euclidean matrix norm.
For every fixed $(z,\theta)$ the norm inside the integral tends to
zero by spatial translation continuity in $L^1(E)$.  It is bounded
by $2\|D_x^2\varphi\|_{L^1(E)}$, so the finite second moment in
\eqref{eq:ht-spatial-moments} justifies the dominated convergence theorem.

For every $w\in L^1(E)$, Tonelli's theorem applied to
\eqref{eq:ht-separate-averages} and the unit mass of $H$ give
\begin{equation}\label{eq:ht-time-average-contraction}
 \|\mathsf H_\epsilon w\|_{L^1(E)}
 \leq\int_0^\infty H(s)\int_{\epsilon^2s}^\infty
               \int_{\Rd}|w(r,x)|\dd x\dd r\dd s
 \leq\|w\|_{L^1(E)}.
\end{equation}
Moreover, \eqref{eq:ht-separate-averages} implies
\begin{equation}\label{eq:ht-time-average-strong}
 \left\|\mathsf H_\epsilon w-w\right\|_{L^1(E)}
 \leq\int_0^\infty H(s)
       \left\|w\left(\cdot+\epsilon^2s,\cdot\right)-w\right\|_{L^1(E)}\dd s
 \longrightarrow0.
\end{equation}
Indeed, extending $w$ by zero to negative times proves translation
continuity, and the integrand is bounded by
$2H(s)\|w\|_{L^1(E)}$, which is integrable because $H$ has mass one.
Combining \eqref{eq:ht-spatial-L1} and
\eqref{eq:ht-time-average-strong} with
\eqref{eq:ht-time-average-contraction} yields
\begin{align}
 \left\|\mathsf H_\epsilon V_\epsilon-\mathcal B\varphi\right\|_{L^1(E)}
 &\leq\left\|V_\epsilon-\mathcal B\varphi\right\|_{L^1(E)}
       +\left\|\mathsf H_\epsilon\mathcal B\varphi
                         -\mathcal B\varphi\right\|_{L^1(E)}
 \longrightarrow0.
 \label{eq:ht-averaged-spatial-L1}
\end{align}
Finally, inserting \eqref{eq:ht-temporal-L1},
\eqref{eq:ht-averaged-spatial-L1}, and \eqref{eq:ht-constants} into
\eqref{eq:ht-operator-decomposition} shows that its left-hand side
converges in $L^1(E)$ to
$-c_\eta D_+^\eta\varphi+\mathcal B\varphi-\kappa\varphi
=(\mathcal G_\eta-\kappa)\varphi$, proving
\eqref{eq:ht-operator-limit}.

The remaining renewal estimates are an application of Lemma
\ref{lem:renewal-framework} with
\[
 \tau_\epsilon=\epsilon^2,\qquad
 \delta_\epsilon=\epsilon^{2\eta},\qquad
 \chi(q)=c_\eta q^\eta,\qquad D_r^\chi=S_r.
\]
The discrete resolvent is \eqref{eq:ht-discrete-resolvent}.
Fix $T<\infty$.
The exponentially tilted convolution measures have
Laplace--Fourier transforms
\begin{equation}\label{eq:ht-resolvent-transform}
 \frac{\delta_\epsilon}
 {1+b\delta_\epsilon-a_\epsilon
    \widehat H(\epsilon^2(q+z))\widehat p(\epsilon^\eta\xi)}
 \xrightarrow[\epsilon\downarrow0]{}
 \frac1{\kappa+b+c_\eta(q+z)^\eta+\xi^{\top}\Sigma\xi/2},
 \quad q>0,\ z\geq0.
\end{equation}
Here $\widehat p$ is the spatial characteristic function.
Part (ii) of that lemma gives \eqref{eq:ht-resolvent-strong},
including when $b=\kappa=0$.
Part (i) gives \eqref{eq:ht-finite-horizon-bound}.
Thus these operators are bounded uniformly in $\epsilon$ in the
unweighted finite-horizon norms as well, with an additional factor
$e^{qT}$.  The finite horizon provides the bound at zero killing.
This completes the proof.
\end{proof}

Fix $g\in C_c^+(E)$ and $T<\infty$ beyond its temporal support.
For a full descendant cluster, give each event mass
$\delta_\epsilon^2$ and evaluate $g$ at the rescaled coordinates.
Let $F_{\epsilon,g}(t,x)$ be the resulting Laplace transform for a
cluster rooted at $(t/\epsilon^2,x/\epsilon^\eta)$, including the root.
Critical and subcritical Poisson Galton--Watson clusters are finite
almost surely.  Define their costs by
\begin{equation}\label{eq:ht-cluster-cost}
 u_{\epsilon,g}:=\frac{1-F_{\epsilon,g}}{\delta_\epsilon},
 \qquad 0\leq F_{\epsilon,g}\leq1,
 \qquad 0\leq u_{\epsilon,g}\leq\delta_\epsilon^{-1}.
\end{equation}
The branching equation is
\begin{equation}\label{eq:ht-cluster-renewal}
 \log\left(1-\delta_\epsilon u_{\epsilon,g}\right)
    =-\delta_\epsilon^2g
       -a_\epsilon\delta_\epsilon\mathsf K_\epsilon u_{\epsilon,g}.
\end{equation}

Equation~\eqref{eq:ht-cluster-renewal} is the heavy-tail specialization
of the common nonlinear renewal framework.  The next lemma records the
resulting convergence of the cluster costs and the spatial integrability
estimates needed to average them over the Poisson past.

\begin{lemma}[Convergence of heavy-tail cluster costs]
\label{lem:ht-cluster-costs}
Under the assumptions of Lemma~\ref{lem:ht-resolvent-estimates},
\eqref{eq:ht-nonlinear-mild} has a unique bounded nonnegative solution
in $\mathcal X_T$, and the cluster costs satisfy
\begin{equation}\label{eq:ht-cluster-uniform}
 \left\|u_{\epsilon,g}-u_g^{\eta,\kappa}\right\|_\infty\xrightarrow[\epsilon\downarrow0]{}0.
\end{equation}
\begin{equation}\label{eq:ht-linear-bound}
 0\leq u_{\epsilon,g}\leq R_\epsilon^{(0)}g,
 \qquad 0\leq u_g^{\eta,\kappa}\leq R_\kappa^\eta g.
\end{equation}
For all sufficiently small $\epsilon$,
\begin{equation}\label{eq:ht-spatial-bound}
 \sup_{0\leq t\leq T}\int_{\Rd}\left(1+|x|^2\right)
               R_\epsilon^{(0)}g(t,x)\dd x\leq C_{g,T}.
\end{equation}
The spatial bound also holds for the limiting linear and nonlinear
costs, and
\begin{equation}\label{eq:ht-cluster-L1}
 \sup_{0\leq t\leq T}\left\|u_{\epsilon,g}(t,\cdot)
                    -u_g^{\eta,\kappa}(t,\cdot)\right\|_1\to0.
\end{equation}
The linear resolvents applied to $g$ satisfy
\begin{equation}\label{eq:ht-linear-L1}
 \sup_{0\leq t\leq T}
       \left\|R_\epsilon^{(0)}g(t,\cdot)
                    -R_\kappa^\eta g(t,\cdot)\right\|_1\to0.
\end{equation}
\end{lemma}

\begin{proof}
Apply Lemma~\ref{lem:renewal-framework} with
$\tau_\epsilon=\epsilon^2$, $\delta_\epsilon=\epsilon^{2\eta}$,
and $\chi(q)=c_\eta q^\eta$.
Its nonlinear statement (iii) gives existence, uniqueness, the linear
majorant, and the uniform and $L^1$ convergence.
The spatial estimate follows from (i) and
\eqref{eq:clock-geometric-bound}, and passes to the limits by
Fatou's lemma.  The linear $L^1$ assertion is part (ii) of the
same lemma.
\end{proof}

\subsubsection{The Poisson past and the stochastic limit}\label{sec:heavy:support:2}

The functional of the Poisson past integrates cluster costs over the entire
prescribed past.  Its bound will also be used when the source amplitude
varies.

For a past parent $(r,y)$, let
\[
  J_{\epsilon,g}(r,y)
       :=\int_{q>0}\int_{\Rd}H(q-r)p(w-y)
              u_{\epsilon,g}(\epsilon^2q,\epsilon^\eta w)\dd w\dd q.
\]
Define the spatial averaging operator
$\mathsf P_\epsilon v(t,x)
=\int_{\Rd}p(z)v(t,x+\epsilon^\eta z)\dd z$ and
\begin{equation}\label{eq:ht-boundary-functional}
  b_\epsilon(t):=\frac{\overline H(t/\epsilon^2)}{\delta_\epsilon},
  \qquad
  \mathfrak B_\epsilon(v)
       :=\int_0^\infty b_\epsilon(t)
             \ip{\lambda_0}{\mathsf P_\epsilon v(t,\cdot)}\dd t.
\end{equation}

The functional $\mathfrak B_\epsilon$ collects the contribution of the
entire negative-time Poisson field after spatial averaging.  The next
lemma gives the exact Poisson Laplace formula, identifies the limiting
memory functional, and establishes the uniform bound needed for
convergence.

\begin{lemma}[Poisson past functional]
\label{lem:ht-preparation-functional}
Assume \eqref{eq:ht-kernel}--\eqref{eq:ht-constants}, and let
$\lambda_0$ be bounded, nonnegative, and uniformly continuous.
For $g\in C_c^+(E)$, at $\alpha=1$ the exact Laplace formula is
\begin{equation}\label{eq:ht-past-Laplace}
  \log\E\left[e^{-\ip{\cN_\epsilon^{\eta,1}}{g}}\right]
       =q_\epsilon\int_{r<0}\int_{\Rd}
         \lambda_0(\epsilon^\eta y)
          \left(e^{-a_\epsilon\delta_\epsilon
                                J_{\epsilon,g}(r,y)}-1\right)
                       \dd y\dd r,
\end{equation}
where $q_\epsilon$ is evaluated at $\alpha=1$.
The boundary densities satisfy
\begin{equation}\label{eq:ht-boundary-limit}
  b_\epsilon(t)\xrightarrow[\epsilon\downarrow0]{} c_Ht^{-\eta}\quad(t>0),
  \qquad b_\epsilon(t)\leq Ct^{-\eta}.
\end{equation}
For every $T<\infty$ and measurable $v$ vanishing for $t\geq T$ with finite
$L^\infty([0,T];L^1(\Rd))$ norm,
\begin{equation}\label{eq:ht-boundary-norm}
 |\mathfrak B_\epsilon(v)|
 \leq\frac{C\|\lambda_0\|_\infty T^{1-\eta}}{1-\eta}
                         \sup_{0\leq t\leq T}\|v(t,\cdot)\|_1.
\end{equation}
For the cluster costs of Lemma~\ref{lem:ht-cluster-costs},
\begin{equation}\label{eq:ht-boundary-cost-limit}
  \mathfrak B_\epsilon(u_{\epsilon,g})
       \xrightarrow[\epsilon\downarrow0]{}\Psi_{\lambda_0}^{\eta,\kappa}(g).
\end{equation}
\end{lemma}

\begin{proof}
For a bounded set $A$ of past parents, let $\nu_{\epsilon,A}$ be the
superposition of their rescaled positive-time clans, each event having
mass $\delta_\epsilon^2$.  In the Poisson exponential formula
\eqref{eq:Poisson-exponential-formula}, the independent child clusters give
$L(r,y)=e^{-a_\epsilon\delta_\epsilon J_{\epsilon,g}(r,y)}$.
Applying \eqref{eq:Poisson-exponential-formula} with
$\Lambda(\dd r,\dd y)=q_\epsilon\lambda_0(\epsilon^\eta y)\dd r\dd y$
yields
\begin{equation}\label{eq:ht-bounded-parent-Laplace}
 \log\E\left[e^{-\ip{\nu_{\epsilon,A}}g}\right]
 =q_\epsilon\int_A\lambda_0(\epsilon^\eta y)
    \left(e^{-a_\epsilon\delta_\epsilon J_{\epsilon,g}(r,y)}-1\right)
                                                \dd r\dd y.
\end{equation}
We compute an integrable bound before exhausting the whole past.
For a nonnegative $v$ supported before $T$, let $J_{\epsilon,v}$ denote
$J_{\epsilon,g}$ with $u_{\epsilon,g}$ replaced by $v$.
Tonelli's theorem and the identity
$\int_{-\infty}^0 H(q-r)\dd r=\overline H(q)$ give
\begin{align}
 &q_\epsilon\delta_\epsilon\int_{-\infty}^0\int_{\Rd}
                 \lambda_0(\epsilon^\eta y)J_{\epsilon,v}(r,y)\dd y\dd r
       \nonumber\\
 &\quad=q_\epsilon\delta_\epsilon\int_0^\infty\overline H(q)
           \int_{\Rd}\int_{\Rd}\lambda_0(\epsilon^\eta y)p(z)
                 v\left(\epsilon^2q,\epsilon^\eta(y+z)\right)\dd z\dd y\dd q
       \nonumber\\
 &\quad=\int_0^T\frac{\overline H(t/\epsilon^2)}{\delta_\epsilon}
           \int_{\Rd}\lambda_0(x)\int_{\Rd}p(z)
                 v\left(t,x+\epsilon^\eta z\right)\dd z\dd x\dd t
       =\mathfrak B_\epsilon(v).
 \label{eq:ht-parent-linear-identity}
\end{align}
Here $q_\epsilon=\epsilon^{2-\eta(4-d)}$ at $\alpha=1$, and the
changes $t=\epsilon^2q$ and $x=\epsilon^\eta y$ use
$q_\epsilon\delta_\epsilon\epsilon^{-2-\eta d}=\delta_\epsilon^{-1}$.
This proves that the linear term in \eqref{eq:ht-past-Laplace} is
$-a_\epsilon\mathfrak B_\epsilon(u_{\epsilon,g})$.

By \eqref{eq:ht-tail-bounds}, for $t>0$,
\begin{align*}
 b_\epsilon(t)
 &=\epsilon^{-2\eta}\overline H\left(t/\epsilon^2\right)
        \longrightarrow c_Ht^{-\eta},\\
 b_\epsilon(t)&\leq C\epsilon^{-2\eta}
                   \min\left\{1,\left(t/\epsilon^2\right)^{-\eta}\right\}
                   \leq Ct^{-\eta}.
\end{align*}
These imply \eqref{eq:ht-boundary-limit}.  Since $\mathsf P_\epsilon$
is an $L^1$ contraction, for signed $v$ we have
\begin{align*}
 |\mathfrak B_\epsilon(v)|
 \leq\|\lambda_0\|_\infty\int_0^T b_\epsilon(t)
                              \|\mathsf P_\epsilon v(t)\|_1\dd t
 \leq\frac{C\|\lambda_0\|_\infty T^{1-\eta}}{1-\eta}
                   \sup_{0\leq t\leq T}\|v(t)\|_1.
\end{align*}
This proves \eqref{eq:ht-boundary-norm}.  The linear majorant and spatial bound
\eqref{eq:ht-linear-bound}--\eqref{eq:ht-spatial-bound} imply
\[
 \sup_{0\leq t\leq T}\left\|u_{\epsilon,g}(t)\right\|_1\leq C_{g,T}.
\]
Thus, the right-hand side of \eqref{eq:ht-boundary-norm} is finite for
$v=u_{\epsilon,g}$;
 hence \eqref{eq:ht-parent-linear-identity} and
$1-e^{-a_\epsilon\delta_\epsilon J}\leq a_\epsilon\delta_\epsilon J$
give
\[
 q_\epsilon\int_{-\infty}^0\int_{\Rd}\lambda_0(\epsilon^\eta y)
       \left(1-e^{-a_\epsilon\delta_\epsilon J_{\epsilon,g}(r,y)}\right)\dd y\dd r
 \leq a_\epsilon\mathfrak B_\epsilon(u_{\epsilon,g})<\infty.
\]
Let $A_n:=[-n,0)\times B_n$.  In the common cluster construction,
\[
 A_n\uparrow(-\infty,0)\times\Rd,\qquad
 \nu_{\epsilon,A_n}\uparrow\cN_\epsilon^{\eta,1},\qquad
 Y_n:=\ip{\nu_{\epsilon,A_n}}g\uparrow
                         Y:=\ip{\cN_\epsilon^{\eta,1}}g.
\]
Define
\[
 I_n:=q_\epsilon\int_{A_n}\lambda_0(\epsilon^\eta y)
       \left(1-e^{-a_\epsilon\delta_\epsilon J_{\epsilon,g}(r,y)}\right)\dd r\dd y.
\]
By \eqref{eq:ht-parent-linear-identity} and
\eqref{eq:ht-boundary-norm}, monotone convergence theorem gives
\begin{align*}
 I_n\uparrow I
 :=q_\epsilon\int_{-\infty}^0\int_{\Rd}
        \lambda_0(\epsilon^\eta y)
       \left(1-e^{-a_\epsilon\delta_\epsilon J_{\epsilon,g}(r,y)}\right)\dd y\dd r
       \leq a_\epsilon\mathfrak B_\epsilon(u_{\epsilon,g})<\infty.
\end{align*}
Since $0\leq e^{-Y_n}\leq1$, bounded convergence theorem and
\eqref{eq:ht-bounded-parent-Laplace} imply
\[
 \E\left[e^{-Y}\right]=\lim_{n\to\infty}\E\left[e^{-Y_n}\right]
          =\lim_{n\to\infty}e^{-I_n}=e^{-I}>0.
\]
Taking logarithms proves \eqref{eq:ht-past-Laplace}.

Write $u=u_g^{\eta,\kappa}$.  For each $0\leq t\leq T$, the
$L^1$ convergence in \eqref{eq:ht-cluster-L1} and translation
continuity give
\begin{align*}
 \left\|\mathsf P_\epsilon u_{\epsilon,g}(t)-u(t)\right\|_1
 \leq\left\|u_{\epsilon,g}(t)-u(t)\right\|_1
 +\int_{\Rd}p(z)
       \left\|u(t,\cdot+\epsilon^\eta z)-u(t,\cdot)\right\|_1\dd z
       \longrightarrow0.
\end{align*}
The integrand in the second term is bounded by $2\|u(t)\|_1$,
so the dominated convergence theorem applies to $p(z)\dd z$.
In particular, the difference of the functionals of the Poisson past satisfies
\begin{align}
 \left|\mathfrak B_\epsilon(u_{\epsilon,g})
                  -c_H\int_0^Tt^{-\eta}\ip{\lambda_0}{u(t)}\dd t\right|
 &\leq\|\lambda_0\|_\infty\int_0^T b_\epsilon(t)
               \left\|\mathsf P_\epsilon u_{\epsilon,g}(t)-u(t)\right\|_1\dd t
       \nonumber\\
 &\qquad+\|\lambda_0\|_\infty\int_0^T
             \left|b_\epsilon(t)-c_Ht^{-\eta}\right|\,\|u(t)\|_1\dd t
       \longrightarrow0.
 \label{eq:ht-boundary-cost-error}
\end{align}
Both integrands are bounded by a constant times $t^{-\eta}$, using
the common $L^1$ bound and \eqref{eq:ht-boundary-limit}.
Since $\int_0^Tt^{-\eta}\dd t=T^{1-\eta}/(1-\eta)<\infty$, the
dominated convergence theorem proves the limit in
\eqref{eq:ht-boundary-cost-error}.  Equivalently, each discarded
interval $[0,\tau]$ has cost at most $C\tau^{1-\eta}/(1-\eta)$.
By \eqref{eq:ht-exponent}, the limiting integral is
$\Psi_{\lambda_0}^{\eta,\kappa}(g)$, proving
\eqref{eq:ht-boundary-cost-limit}.
\end{proof}

The quadratic error in linearizing each past parent's contribution
vanishes when the past is scaled at the fluctuation threshold. This yields the
limiting Laplace exponent.

\begin{lemma}[Convergence of the critical Laplace exponent]
\label{lem:ht-critical-exponent}
Under the assumptions of Lemma~\ref{lem:ht-preparation-functional},
for every $g\in C_c^+(E)$,
\begin{equation}\label{eq:ht-critical-Laplace-limit}
 \log\E\left[e^{-\ip{\cN_\epsilon^{\eta,1}}{g}}\right]
       \xrightarrow[\epsilon\downarrow0]{}-\Psi_{\lambda_0}^{\eta,\kappa}(g).
\end{equation}
More precisely, replacing the exponential in
\eqref{eq:ht-past-Laplace} by its linear term incurs a nonnegative
error bounded by
\begin{equation}\label{eq:ht-parent-error}
  \frac{a_\epsilon\delta_\epsilon}{2}
       \|u_{\epsilon,g}\|_\infty\,
          a_\epsilon\mathfrak B_\epsilon(u_{\epsilon,g})
       =\mathcal O(\delta_\epsilon).
\end{equation}
\end{lemma}

\begin{proof}
Since $0\leq J_{\epsilon,g}\leq\|u_{\epsilon,g}\|_\infty$,
the inequality $0\leq e^{-z}-1+z\leq z^2/2$ gives
\eqref{eq:ht-parent-error}.  The uniform cluster bound from
Lemma~\ref{lem:ht-cluster-costs} and the boundary estimate from
Lemma~\ref{lem:ht-preparation-functional} make this $\mathcal O(\delta_\epsilon)$.
The linear term is
$-a_\epsilon\mathfrak B_\epsilon(u_{\epsilon,g})$.
Use \eqref{eq:ht-boundary-cost-limit} and $a_\epsilon\to1$ to obtain
\eqref{eq:ht-critical-Laplace-limit}.
\end{proof}

\subsubsection{Mean, high-density, and sparse limits}\label{sec:heavy:support:3}

The linear resolvent determines both the limiting mean of the event
measures and the derivative of the limiting Laplace exponent at zero.
The latter will identify the first moment of the random limit.

\begin{lemma}[Linear limit and fractional mean profile]
\label{lem:ht-linear-mean}
Under the assumptions of Lemma~\ref{lem:ht-preparation-functional},
for every $\alpha>0$ and $g\in C_c^+(E)$,
\begin{equation}\label{eq:ht-mean-exact}
  \E\left[\ip{\cN_\epsilon^{\eta,\alpha}}{g}\right]
       =a_\epsilon\mathfrak B_\epsilon\left(R_\epsilon^{(0)}g\right)
       \xrightarrow[\epsilon\downarrow0]{}
       c_H\int_0^\infty t^{-\eta}
                     \left\langle\lambda_0,R_\kappa^\eta g(t,\cdot)\right\rangle\dd t.
\end{equation}
The limit on the right equals
$\int_E g(t,x)\rho^{\eta,\kappa}(t,x)\dd t\dd x$, where
\eqref{eq:ht-mean-profile} is the unique bounded continuous mild
solution of \eqref{eq:ht-caputo} with initial trace $\lambda_0$.
Moreover,
\begin{equation}\label{eq:ht-exponent-derivative}
 \lim_{\theta\downarrow0}
   \frac{\Psi_{\lambda_0}^{\eta,\kappa}(\theta g)}{\theta}
     =\int_E g(t,x)\rho^{\eta,\kappa}(t,x)\dd t\dd x.
\end{equation}
\end{lemma}

\begin{proof}
Fix $g\in C_c^+(E)$ and choose $T$ beyond its temporal support.
For a full microscopic cluster rooted at
$(t/\epsilon^2,x/\epsilon^\eta)$, the expected sum of the values of $g$
over generation $n$ is $a_\epsilon^n\mathsf K_\epsilon^ng(t,x)$.
Consequently, with event mass
$w_\epsilon:=\epsilon^{\eta(\alpha+3)}$,
\begin{equation}\label{eq:ht-full-cluster-mean}
 \E_{t/\epsilon^2,x/\epsilon^\eta}
       \left[w_\epsilon\int_{[0,\infty)\times\Rd}
         g\left(\epsilon^2s,\epsilon^\eta z\right)C(\dd s,\dd z)\right]
 =w_\epsilon\sum_{n=0}^\infty a_\epsilon^n\mathsf K_\epsilon^ng(t,x)
 =\frac{w_\epsilon}{\delta_\epsilon}R_\epsilon^{(0)}g(t,x).
\end{equation}
The series is finite by \eqref{eq:ht-finite-horizon-bound}.
Integrating \eqref{eq:ht-full-cluster-mean} over the direct
positive-time children of every past parent gives, by Tonelli's theorem,
\begin{align}
 \E\left[\ip{\cN_\epsilon^{\eta,\alpha}}{g}\right]
 &=\frac{a_\epsilon q_\epsilon w_\epsilon}{\delta_\epsilon}
       \int_{-\infty}^0\int_{\Rd}\lambda_0(\epsilon^\eta y)
       \int_0^\infty\int_{\Rd}H(q-r)p(z)
       R_\epsilon^{(0)}g\left(\epsilon^2q,\epsilon^\eta(y+z)\right)
                    \dd z\dd q\dd y\dd r\nonumber\\
 &=a_\epsilon\mathfrak B_\epsilon\left(R_\epsilon^{(0)}g\right).
\label{eq:ht-mean-parent-integration}
\end{align}
For the second equality in \eqref{eq:ht-mean-parent-integration},
 integrate $r$ to obtain $\overline H(q)$
and change $t=\epsilon^2q$, $x=\epsilon^\eta y$ as in
\eqref{eq:ht-parent-linear-identity}.  For every $\alpha>0$, the
coefficient is
\[
 \frac{q_\epsilon w_\epsilon}{\delta_\epsilon}
                \epsilon^{-2-\eta d}=\delta_\epsilon^{-1}.
\]
By \eqref{eq:ht-mean-parent-integration}, \eqref{eq:ht-boundary-norm},
\eqref{eq:ht-spatial-bound}, and $a_\epsilon\leq1$,
\begin{align}
 \E\left[\ip{\cN_\epsilon^{\eta,\alpha}}{g}\right]
 &\leq\frac{C\|\lambda_0\|_\infty T^{1-\eta}}{1-\eta}
       \sup_{0\leq t\leq T}\left\|R_\epsilon^{(0)}g(t,\cdot)\right\|_1
       \nonumber\\
 &\leq\frac{C C_{g,T}\|\lambda_0\|_\infty T^{1-\eta}}{1-\eta}<\infty.
 \label{eq:ht-event-mean-bound}
\end{align}
This bound is uniform over $\alpha>0$ and sufficiently small $\epsilon$.
The linear costs $R_\epsilon^{(0)}g$ and $R_\kappa^\eta g$ are defined
by \eqref{eq:ht-discrete-resolvent} and \eqref{eq:ht-resolvent},
respectively.  Apply the estimate \eqref{eq:ht-boundary-cost-error}
to this pair in place of $u_{\epsilon,g}$ and $u_g^{\eta,\kappa}$,
defined by \eqref{eq:ht-cluster-cost} and \eqref{eq:ht-nonlinear-mild}.
Their convergence in \eqref{eq:ht-linear-L1} and the spatial bounds
in \eqref{eq:ht-spatial-bound} then give
\[
 a_\epsilon\mathfrak B_\epsilon\left(R_\epsilon^{(0)}g\right)
 \longrightarrow c_H\int_0^T t^{-\eta}
                    \ip{\lambda_0}{R_\kappa^\eta g(t)}\dd t.
\]
This proves \eqref{eq:ht-mean-exact}.

For the derivative of the limiting exponent, let
$v:=R_\kappa^\eta g$ and $u_\theta:=u_{\theta g}^{\eta,\kappa}$.
Equations \eqref{eq:ht-nonlinear-mild} and \eqref{eq:ht-linear-bound}
give, for $0<\theta\leq1$,
\[
 0\leq v-u_\theta/\theta
       =\frac1{2\theta}R_\kappa^\eta(u_\theta^2),
 \qquad 0\leq u_\theta\leq\theta v.
\]
The finite-horizon bounds \eqref{eq:ht-finite-horizon-bound}, with
$p=1,\infty$, imply
\begin{align}
 \sup_{0\leq t\leq T}\|v(t)-u_\theta(t)/\theta\|_1
 &\leq\frac{e^{qT}}{2\theta(\kappa+c_\eta q^\eta)}
               \sup_{0\leq t\leq T}\left\|(u_\theta(t))^2\right\|_1\nonumber\\
 &\leq\frac{\theta e^{qT}}{2(\kappa+c_\eta q^\eta)}
                \|v\|_\infty\sup_{0\leq t\leq T}\|v(t)\|_1
 \leq C_{g,T}\theta,
 \label{eq:ht-linearization-L1}
\end{align}
where $q>0$ is fixed.  Thus \eqref{eq:ht-exponent} gives
\begin{equation}\label{eq:ht-exponent-linearization}
 \left|\frac{\Psi_{\lambda_0}^{\eta,\kappa}(\theta g)}{\theta}
        -c_H\int_0^T t^{-\eta}\ip{\lambda_0}{v(t)}\dd t\right|
 \leq\frac{c_H\|\lambda_0\|_\infty T^{1-\eta}}{1-\eta}
                          C_{g,T}\theta\longrightarrow0.
\end{equation}

We identify the common linear limit by its space--time kernel.
Let
\begin{equation}\label{eq:ht-potential-measure}
 U_\kappa(\dd s,\dd z)
 :=\int_0^\infty e^{-\kappa r}
       \mathbb P\left(S_r\in\dd s,B_r^\Sigma\in\dd z\right)\dd r
\end{equation}
be the potential measure of \eqref{eq:ht-resolvent}.  Independence of
$S$ and $B^\Sigma$ gives, for $z>0$ and $\xi\in\Rd$,
\begin{equation}\label{eq:ht-potential-transform}
 \int_{[0,\infty)\times\Rd}e^{-zs+\mathrm i\xi\cdot y}
                      U_\kappa(\dd s,\dd y)
 =\int_0^\infty e^{-r(\kappa+c_\eta z^\eta+\xi^\top\Sigma\xi/2)}\dd r
 =\frac1{\kappa+c_\eta z^\eta+\xi^\top\Sigma\xi/2}.
\end{equation}
In particular its temporal tilt is finite even when $\kappa=0$.
Define the physical-time kernel, for almost every $t>0$, by
\begin{equation}\label{eq:ht-mean-kernel}
 K_\kappa(t,\dd y)
 :=c_H\int_{[0,t)}(t-s)^{-\eta}U_\kappa(\dd s,\dd y).
\end{equation}
Fubini's theorem applies after the exponential temporal tilt, and
\eqref{eq:ht-potential-transform} implies
\begin{align}
 \int_0^\infty e^{-zt}\int_{\Rd}e^{\mathrm i\xi\cdot y}
                         K_\kappa(t,\dd y)\dd t
 =\frac{c_H\int_0^\infty e^{-za}a^{-\eta}\dd a}
             {\kappa+c_\eta z^\eta+\xi^\top\Sigma\xi/2}
 =\frac{c_\eta z^{\eta-1}}
             {c_\eta z^\eta+\kappa+\xi^\top\Sigma\xi/2}.
 \label{eq:ht-mean-transform}
\end{align}
Equations \eqref{eq:ht-resolvent} and \eqref{eq:ht-potential-measure}
give the forward representation
\begin{equation}\label{eq:ht-forward-potential-representation}
 (R_\kappa^\eta g)(a,z)
 =\int_{[0,\infty)\times\Rd}g(a+s,z+y)U_\kappa(\dd s,\dd y).
\end{equation}
Substitute \eqref{eq:ht-forward-potential-representation} into the
limiting pairing in \eqref{eq:ht-mean-exact}.  Tonelli's theorem and
the change of variables $t=a+s$, $x=z+y$, together with
\eqref{eq:ht-mean-kernel}, give
\begin{equation}\label{eq:ht-linear-kernel-pairing}
 c_H\int_0^T a^{-\eta}\ip{\lambda_0}{R_\kappa^\eta g(a)}\dd a
 =\int_0^T\int_{\Rd}g(t,x)
                \int_{\Rd}\lambda_0(x-y)K_\kappa(t,\dd y)\dd x\dd t.
\end{equation}

For $E_t=\inf\{r:S_r>t\}$, the relation
$\mathbb P(E_t>r)=\mathbb P(S_r\leq t)$ gives
\begin{align}
 \int_0^\infty e^{-zt}\mathbb P(E_t>r)\dd t
 &=\E\left[\int_{S_r}^\infty e^{-zt}\dd t\right]
   =z^{-1}e^{-rc_\eta z^\eta},\label{to:take:derivative}\\
 \int_0^\infty e^{-zt}\mathbb P(E_t\in\dd r)\dd t
 &=c_\eta z^{\eta-1}e^{-rc_\eta z^\eta}\dd r,
 \label{eq:ht-inverse-clock-transform}
\end{align}
where the equality in \eqref{eq:ht-inverse-clock-transform} follows by taking the negative distributional
derivative in $r$ with respect to \eqref{to:take:derivative}.
Using the Brownian transition measure $p_r^\Sigma(\dd y)$ defined in
\eqref{eq:ht-brownian-kernel}, let
\[
 K'_\kappa(t,\dd y)
 :=\int_0^\infty e^{-\kappa r}p_r^\Sigma(\dd y)
                              \mathbb P(E_t\in\dd r).
\]
This measure has the transform \eqref{eq:ht-mean-transform}, by
\eqref{eq:ht-inverse-clock-transform}.  Uniqueness of the
Laplace--Fourier transform for the finite tilted measures gives
$K_\kappa(t,\dd y)\dd t=K'_\kappa(t,\dd y)\dd t$.
Substitution in \eqref{eq:ht-linear-kernel-pairing} identifies its
inner integral as
\[
 \int_{\Rd}\lambda_0(x-y)K'_\kappa(t,\dd y)
 =\E\left[e^{-\kappa E_t}P_{E_t}^\Sigma\lambda_0(x)\right]
 =\rho^{\eta,\kappa}(t,x).
\]
This proves the mean identification and, by
\eqref{eq:ht-exponent-linearization}, \eqref{eq:ht-exponent-derivative}.

Finally, define the time Laplace transform as a bounded-function-valued
integral.  Equation \eqref{eq:ht-inverse-clock-transform} gives
\begin{align}
 \widetilde\rho(z)
 &:=\int_0^\infty e^{-zt}\rho^{\eta,\kappa}(t,\cdot)\dd t\nonumber\\
 &=c_\eta z^{\eta-1}\int_0^\infty
             e^{-r(c_\eta z^\eta+\kappa)}P_r^\Sigma\lambda_0\dd r
 =c_\eta z^{\eta-1}\left(c_\eta z^\eta+\kappa-\mathcal B\right)^{-1}\lambda_0,
 \label{eq:ht-caputo-resolvent}
\end{align}
so that
\[
 c_\eta\left(z^\eta\widetilde\rho(z)-z^{\eta-1}\lambda_0\right)
       =(\mathcal B-\kappa)\widetilde\rho(z).
\]
This is the Laplace form of the mild Caputo equation
\eqref{eq:ht-caputo}; it applies to bounded uniformly continuous
$\lambda_0$ without taking its spatial Fourier transform.
The stable subordinator is strictly increasing, and $E_t\to0$
almost surely as $t\downarrow0$.  Strong continuity of $P_r^\Sigma$
on bounded uniformly continuous functions and
\[
 \left\|\rho^{\eta,\kappa}(t)-\lambda_0\right\|_\infty
 \leq\E\left[\left\|P_{E_t}^\Sigma\lambda_0-\lambda_0\right\|_\infty\right]
       +\|\lambda_0\|_\infty\E\left[1-e^{-\kappa E_t}\right]
       \longrightarrow0
\]
give the initial trace.  The inverse clock is continuous, and the same
semigroup estimate gives continuity in $t$.  Every bounded continuous
mild solution has the transform \eqref{eq:ht-caputo-resolvent};
uniqueness of the time Laplace transform and continuity therefore give
uniqueness of the solution.
\end{proof}

The other density regimes require control of small and growing
Laplace sources.  A fixed physical-time horizon makes these estimates
uniform even when $\kappa=0$.

\begin{lemma}[Bounds for varying Laplace sources]
\label{lem:ht-varying-sources}
Under the assumptions of Lemma~\ref{lem:ht-preparation-functional},
fix $g\in C_c^+(E)$ and $T<\infty$ such that $g(t,x)=0$ for
$t\geq T$.  There is $C_{g,T}<\infty$, independent of the amplitudes
below and of sufficiently small $\epsilon$, such that:
\begin{enumerate}[label=(\roman*)]
\item For every $r\geq1$,
\begin{equation}\label{eq:ht-small-source}
  \sup_{0\leq t\leq T}\left\|
     r u_{\epsilon,g/r}(t,\cdot)
                 -R_\epsilon^{(0)}g(t,\cdot)\right\|_1
       \leq\frac{C_{g,T}}{r}.
\end{equation}
In the Poisson formula with profile $r\lambda_0$ and source $g/r$,
the error from linearizing the parent exponential is at most
$C_{g,T}\delta_\epsilon/r$.
\item For every $\theta\geq1$,
\begin{equation}\label{eq:ht-growing-source}
  \sup_{0\leq t\leq T}\left\|u_{\epsilon,\theta g}(t,\cdot)\right\|_1
       \leq C_{g,T}\theta^{(d+1)/(d+2)}.
\end{equation}
\end{enumerate}
\end{lemma}

\begin{proof}
\textit{(i).}
Use
$u_{\epsilon,g/r}
 \leq r^{-1}R_\epsilon^{(0)}g$ and the exact identity
\[
  R_\epsilon^{(0)}(g/r)
          -u_{\epsilon,g/r}
       =R_\epsilon^{(0)}
                    \Phi_\epsilon(u_{\epsilon,g/r}).
\]
On this bounded range $\Phi_\epsilon(v)\leq Cv^2$.
The uniform finite-horizon resolvent bounds
\eqref{eq:ht-finite-horizon-bound}, for $p=1,\infty$, then prove \eqref{eq:ht-small-source}.  The boundary functional
of its error is also $\mathcal O(r^{-1})$ since
$\int_0^T b_\epsilon(t)\dd t$ is bounded.
The parent linearization error, including the rescaled intensity of the past,
is $\mathcal O(\delta_\epsilon/r)$ by the same estimate as
\eqref{eq:ht-parent-error}.

\textit{(ii).}
For $\theta\geq1$, comparison with a constant in
\eqref{eq:ht-cluster-renewal} gives
$\|u_{\epsilon,\theta g}\|_\infty\leq C_{g,T}\sqrt\theta$;
if this proposed bound exceeds $\delta_\epsilon^{-1}$ the assertion
follows from $u_{\epsilon,\theta g}\leq\delta_\epsilon^{-1}$ in
\eqref{eq:ht-cluster-cost}, applied with source $\theta g$.
The linear bound \eqref{eq:ht-linear-bound} and linearity of
$R_\epsilon^{(0)}$ give
$u_{\epsilon,\theta g}\leq\theta R_\epsilon^{(0)}g$.
Consequently, \eqref{eq:ht-spatial-bound} implies
\begin{align*}
 \sup_{0\leq t\leq T}\int_{\Rd}|x|^2u_{\epsilon,\theta g}(t,x)\dd x
 &\leq\theta\sup_{0\leq t\leq T}\int_{\Rd}|x|^2
                            R_\epsilon^{(0)}g(t,x)\dd x\\
 &\leq C_{g,T}\theta.
\end{align*}
Splitting space into $B_R$ and its complement gives
\[
  \sup_{0\leq t\leq T}\left\|u_{\epsilon,\theta g}(t,\cdot)\right\|_1
       \leq C_{g,T}\left(\sqrt\theta R^d+\theta R^{-2}\right).
\]
Take $R=\theta^{1/(2(d+2))}$ to obtain
\eqref{eq:ht-growing-source}.
\end{proof}

\subsubsection{A superprocess on space--time}

The fractional law has a Markovian realization on $\overline E$,
with physical time included in the state.  Define
\begin{equation}\label{eq:ht-lifted-initial}
 \mu_H(\dd t,\dd x):=c_Ht^{-\eta}\lambda_0(x)\dd t\dd x.
\end{equation}
This is a locally finite measure, with no atom at $t=0$.
The semigroup $\mathsf P^\eta$ in \eqref{eq:ht-resolvent} describes motion
$(t,x)\mapsto(t+S_r,x+B_r^\Sigma)$ in genealogical time $r$.

The construction in the following proposition realizes the limiting Laplace law as a
superprocess occupation measure.  Increasing finite-measure
restrictions accommodate the possibly infinite mass of $\mu_H$.

\begin{proposition}[Space--time superprocess representation]
\label{prop:ht-lifted}
For every $\kappa\geq0$, let $\mu_0=0$ and
$\mu_n=\mu_H|_{[0,n]\times B_n}$ for $n\geq1$.
For each $n$, let
$\mathsf Z^{[n]}$ be the finite-measure superprocess on $\overline E$
with motion generator $\mathcal G_\eta$, branching mechanism
$\psi(v)=\kappa v+v^2/2$, and initial measure $\mu_n$.
It is characterized by the continuous martingales
\begin{align}
 \mathsf M_r^{[n]}(\varphi)
  &:=\ip{\mathsf Z_r^{[n]}}{\varphi}-\ip{\mu_n}{\varphi}
     -\int_0^r\ip{\mathsf Z_s^{[n]}}{
              \mathcal G_\eta\varphi-\kappa\varphi}\dd s,
       \label{eq:ht-lifted-martingale}\\
 \left\langle\mathsf M^{[n]}(\varphi)\right\rangle_r
  &=\int_0^r\ip{\mathsf Z_s^{[n]}}{\varphi^2}\dd s,
       \label{eq:ht-lifted-bracket}
\end{align}
for bounded test functions in the generator domain.
Couple these processes increasingly by independent superprocesses
started from the successive differences $\mu_n-\mu_{n-1}$.
Then the occupation measures have an increasing, locally finite limit,
and
\begin{equation}\label{eq:ht-lifted-occupation}
 \nu_{\lambda_0}^{\eta,\kappa}
   \ \stackrel{\mathrm{law}}{=}\quad
   \lim_{n\to\infty}\int_0^\infty\mathsf Z_r^{[n]}\dd r
       \quad\text{on }E.
\end{equation}
Its law is independent of the finite-measure exhaustion.
\end{proposition}

\begin{proof}
The process $L_r:=\left(S_r,B_r^\Sigma\right)$ is a L\'evy process in
$\R^{d+1}$ with c\`adl\`ag paths.  The strong Markov property and
quasi-left-continuity of L\'evy processes are given in
\cite[Theorem~3.1 and Lemma~3.2]{kyprianou2006levy}; the stopping-time
arguments there apply to vector-valued increments as well.
In particular, for stopping times $\tau_n\uparrow\tau$,
\[
 L_{\tau_n}\longrightarrow L_\tau
 \quad\text{almost surely on }\{\tau<\infty\}.
\]
Thus $L$, translated by its initial point, is a Hunt process, i.e.
a c\`adl\`ag strong Markov process with this quasi-left-continuity
property.  Since $S$ is nondecreasing, a starting point in
$\overline E$ stays in $\overline E$.
For $f\in C_0(\overline E)$, extend $f$ to a uniformly continuous
function on $\R^{d+1}$ by $f(t,x)=f(0,x)$ for $t<0$ when estimating
translations.  Stochastic continuity gives, for $a>0$,
\[
 \left\|\mathsf P_r^\eta f-f\right\|_\infty
 \leq\omega_f(a)+2\|f\|_\infty\mathbb P(|L_r|>a)
 \longrightarrow\omega_f(a)\qquad(r\downarrow0).
\]
Letting $a\downarrow0$ proves strong continuity.
For $f\in C_c(\overline E)$ with support in $[0,T]\times B_R$,
monotonicity of $S_r$ and the triangle inequality give
\begin{align}
 \mathsf P_r^\eta f(t,x)&=0\qquad(t>T),\nonumber\\
 \sup_{0\leq t\leq T}\left|\mathsf P_r^\eta f(t,x)\right|
 &\leq\|f\|_\infty\mathbb P\left(x+B_r^\Sigma\in B_R\right)
 \leq\|f\|_\infty\mathbb P\left(\left|B_r^\Sigma\right|\geq|x|-R\right)
 \xrightarrow[|x|\to\infty]{}0.
 \label{eq:ht-motion-vanishing}
\end{align}
For any $f\in C_0(\overline E)$ choose $f_n\in C_c(\overline E)$
with $\|f_n-f\|_\infty\to0$.  The contraction estimate
\[
 \left\|\mathsf P_r^\eta f_n-\mathsf P_r^\eta f\right\|_\infty
 \leq\|f_n-f\|_\infty,
\]
and \eqref{eq:ht-motion-vanishing} prove preservation of
$C_0(\overline E)$ once continuity is checked.  For
$f\in C_b(\overline E)$ and $(t_n,x_n)\to(t,x)$, the dominated
convergence theorem gives
\[
 \left|\mathsf P_r^\eta f(t_n,x_n)-\mathsf P_r^\eta f(t,x)\right|
 \leq\E\left[\left|f\left(t_n+S_r,x_n+B_r^\Sigma\right)
                         -f\left(t+S_r,x+B_r^\Sigma\right)\right|\right]\longrightarrow0.
\]
Thus the motion $(t,x)\mapsto\left(t+S_r,x+B_r^\Sigma\right)$, whose semigroup
is $\mathsf P_r^\eta$ in \eqref{eq:ht-resolvent}, preserves both
$C_0(\overline E)$ and bounded continuous functions.

The finite-measure construction and martingale characterization follow
from \cite[Theorem~II.5.1]{perkins2002superprocess}, with branching
variance one and linear drift $-\kappa$.
For finite initial measure $\mu$, the occupation Laplace formula
\cite[Theorem~II.5.11(a),(c)]{perkins2002superprocess} gives, for $R<\infty$,
\begin{equation}\label{eq:ht-finite-genealogical-Laplace}
 \E_\mu\left[e^{-\int_0^R\ip{\mathsf Z_r}g\dd r}\right]
       =e^{-\ip\mu{u_R}},\qquad
 u_R=\int_0^R e^{-\kappa s}\mathsf P_s^\eta
                           \left(g-u_{R-s}^2/2\right)\dd s.
\end{equation}
Take $T$ beyond the physical-time support of $g$.  The first-moment
semigroup is $e^{-\kappa r}\mathsf P_r^\eta$.  For a single initial
unit $\delta_{(t,x)}$, this gives the explicit tail estimate
\begin{align}
 \E_{\delta_{(t,x)}}\left[\int_R^\infty\ip{\mathsf Z_r}g\dd r\right]
 &=\int_R^\infty e^{-\kappa r}
                   \E\left[g\left(t+S_r,x+B_r^\Sigma\right)\right]\dd r\nonumber\\
 &\leq\|g\|_\infty\int_R^\infty e^{-\kappa r}
                                 \mathbb P(S_r\leq T)\dd r\nonumber\\
 &\leq\frac{\|g\|_\infty e^{qT}}{\kappa+c_\eta q^\eta}
                  e^{-(\kappa+c_\eta q^\eta)R},\qquad q>0.
 \label{eq:ht-genealogical-tail}
\end{align}
The last inequality uses
$\mathbf1_{\{S_r\leq T\}}\leq e^{qT-qS_r}$ and
$\E[e^{-qS_r}]=e^{-c_\eta q^\eta r}$.  The bound is uniform in
$(t,x)\in\overline E$.

For a unit initial measure, set
$Y_R:=\int_0^R\ip{\mathsf Z_r}g\dd r$ and $Y_\infty:=\lim_{R\to\infty}Y_R$.
Jensen's inequality gives
$0\leq u_R\uparrow u_\infty\leq R_\kappa^\eta g$.
Let $C:=\|R_\kappa^\eta g\|_\infty<\infty$, by
\eqref{eq:ht-finite-horizon-bound}.  Since
$\E[e^{-Y_\infty}]\geq e^{-C}$,
\begin{align*}
 0\leq u_\infty(t,x)-u_R(t,x)
 &=\log\frac{\E_{\delta_{(t,x)}}[e^{-Y_R}]}
                 {\E_{\delta_{(t,x)}}[e^{-Y_\infty}]}\\
 &\leq e^C\E_{\delta_{(t,x)}}[Y_\infty-Y_R]
 \leq\frac{e^C\|g\|_\infty e^{qT}}{\kappa+c_\eta q^\eta}
                    e^{-(\kappa+c_\eta q^\eta)R}.
\end{align*}
Thus, \eqref{eq:ht-genealogical-tail} proves uniform convergence.
We next verify that $u_R\in\mathcal X_T$.  Let
$Q_s:=e^{-\kappa s}\mathsf P_s^\eta$ and recall
$\mathcal X_T=\left\{v\in C_0(\overline E):v(t,x)=0\text{ for }t\geq T\right\}$.
This is a closed subspace of $C_0(\overline E)$ in the uniform norm.
The motion properties proved in \eqref{eq:ht-motion-vanishing} and
the surrounding argument imply that $Q_s$ is a strongly continuous
contraction semigroup on this space: since $S_s\geq0$,
\[
 Q_sv(t,x)=e^{-\kappa s}\E\left[v\left(t+S_s,x+B_s^\Sigma\right)\right]=0,
 \qquad(t\geq T,\ v\in\mathcal X_T).
\]
To make the iteration globally Lipschitz, define
\[
 \Phi_C(z):=\tfrac12\left((z\vee0)\wedge C\right)^2,
 \qquad |\Phi_C(z)-\Phi_C(z')|\leq C|z-z'|,
 \qquad \Phi_C(0)=0.
\]
Thus $\Phi_C(v)\in\mathcal X_T$ for every $v\in\mathcal X_T$.
Choose $a>0$ with $aC<1$ and let
$\mathcal Y_a:=C([0,a];\mathcal X_T)$, with norm
$\|v\|_{\mathcal Y_a}:=\sup_{0\leq r\leq a}\|v(r)\|_\infty$.
For an initial function $v_*\in\mathcal X_T$, define
\begin{equation}\label{eq:ht-Picard-map}
 (\mathcal T_{v_*}v)(r)
 :=Q_rv_*+\int_0^r Q_s\left[g-\Phi_C(v(r-s))\right]\dd s,
 \qquad 0\leq r\leq a.
\end{equation}
Strong continuity of $Q_s$ and the Lipschitz bound for $\Phi_C$
show that $\mathcal T_{v_*}$ maps $\mathcal Y_a$ into itself.
Moreover,
\begin{equation}\label{eq:ht-Picard-contraction}
 \left\|\mathcal T_{v_*}v-\mathcal T_{v_*}\widetilde v\right\|_{\mathcal Y_a}
 \leq aC\left\|v-\widetilde v\right\|_{\mathcal Y_a}.
\end{equation}
Starting from $v^{(0)}(r)=Q_rv_*$, the Picard iterates
$v^{(k+1)}=\mathcal T_{v_*}v^{(k)}$, $k\geq0$, remain in
$\mathcal Y_a$, and \eqref{eq:ht-Picard-contraction} gives
\[
 \left\|v^{(k+1)}-v^{(k)}\right\|_{\mathcal Y_a}
 \leq(aC)^k\left\|v^{(1)}-v^{(0)}\right\|_{\mathcal Y_a}.
\]
Their uniform limit $v\in\mathcal Y_a$ is the unique solution of
$v=\mathcal T_{v_*}v$, with $\mathcal T_{v_*}$ defined in
\eqref{eq:ht-Picard-map}; explicitly,
\begin{equation}\label{eq:ht-Picard-fixed-point}
 v(r)=Q_rv_*+\int_0^r Q_s\left[g-\Phi_C(v(r-s))\right]\dd s,
 \qquad 0\leq r\leq a.
\end{equation}
The contraction estimate \eqref{eq:ht-Picard-contraction} also gives
uniqueness for \eqref{eq:ht-Picard-fixed-point} among bounded measurable solutions on
$[0,a]\times\overline E$.

Since $0\leq u_R\leq C$, the nonlinearity $u_R^2/2$ in
\eqref{eq:ht-finite-genealogical-Laplace} equals $\Phi_C(u_R)$.
Taking $v_*=u_0=0$ therefore identifies the first fixed point with
$r\mapsto u_r$ for $0\leq r\leq a$.
For a genealogical time $R_0\geq0$, splitting the integral in
\eqref{eq:ht-finite-genealogical-Laplace} at $s=r$ and using the
semigroup property yields
\begin{equation}\label{eq:ht-Picard-restart}
 u_{R_0+r}=Q_ru_{R_0}
       +\int_0^r Q_s\left[g-\Phi_C(u_{R_0+r-s})\right]\dd s,
 \qquad 0\leq r\leq a.
\end{equation}
Once $u_{R_0}\in\mathcal X_T$, this is \eqref{eq:ht-Picard-map}
with $v_*=u_{R_0}$.  Applying \eqref{eq:ht-Picard-contraction}
successively at $R_0=ja$, $j=0,1,\ldots$, proves
$u_R\in\mathcal X_T$ for every finite $R$.
The uniform limit $u_\infty$ also belongs to this closed space.
In \eqref{eq:ht-finite-genealogical-Laplace} the integrand
is bounded in absolute value by
\[
 \left(\|g\|_\infty+C^2/2\right)e^{-\kappa s}\mathbb P(S_s\leq T)
 \leq\left(\|g\|_\infty+C^2/2\right)e^{qT-(\kappa+c_\eta q^\eta)s}.
\]
The dominated convergence theorem therefore gives
$u_\infty=R_\kappa^\eta(g-u_\infty^2/2)$.
Uniqueness in \eqref{eq:ht-nonlinear-mild} identifies
$u_\infty=u_g^{\eta,\kappa}$.  Bounded convergence in
\eqref{eq:ht-finite-genealogical-Laplace} now yields
\[
 \E_\mu\left[e^{-\int_0^\infty\ip{\mathsf Z_r}g\dd r}\right]
       =e^{-\ip\mu{u_g^{\eta,\kappa}}}.
\]
All the estimates remain valid for $\kappa=0$.

For $g\in C_c^+(E)$ supported before $T$, let
$g_*(x):=\sup_{0\leq t\leq T} g(t,x)$.  Applying
\eqref{eq:ht-finite-horizon-bound} with $b=0$ and $p=1$ gives,
for every $q>0$,
\begin{align*}
 \sup_{0\leq t\leq T}\int_{\Rd}R_\kappa^\eta g(t,x)\dd x
 &\leq\left\|R_\kappa^\eta g\right\|_{q,1}
 \leq\frac{\|g\|_{q,1}}{\kappa+c_\eta q^\eta}
 \leq\frac{e^{qT}}{\kappa+c_\eta q^\eta}\left\|g_*\right\|_1.
\end{align*}
Since $S_r\geq0$ and $g=0$ for $t\geq T$, the definition
\eqref{eq:ht-resolvent} also gives
$R_\kappa^\eta g(t,\cdot)=0$ for $t\geq T$.
Consequently,
\begin{equation}\label{eq:ht-lifted-local-bound}
 \ip{\mu_H}{R_\kappa^\eta g}
 \leq\frac{c_H\|\lambda_0\|_\infty T^{1-\eta}e^{qT}}
               {(1-\eta)(\kappa+c_\eta q^\eta)}\left\|g_*\right\|_1<\infty.
\end{equation}
This bounds the expected occupation of every compact set, by choosing
$g$ equal to one there.  The increasing occupation limit is therefore
locally finite almost surely.  Monotone convergence of its Laplace
functional yields
\[
\exp\left\{-\ip{\mu_H}{u_g^{\eta,\kappa}}\right\}
=\exp\left\{-\Psi_{\lambda_0}^{\eta,\kappa}(g)\right\}.
\]
This proves \eqref{eq:ht-lifted-occupation} and independence of the
exhaustion.
\end{proof}

Branches share their ancestral history, and after branching their
stable time increments evolve independently.  The occupation integral
is taken in genealogical time and is a measure on physical space--time.
The surviving past provides the initial measure $\mu_H$ for this
enlarged-state process.

\begin{corollary}[Two temporal tails, one occupation construction]
\label{cor:clock-unification}
Let $\kappa\geq0$.  Under the respective finite-mean and heavy-tail
assumptions, the stochastic limits at $\alpha=1$ admit the following
occupation representations.  In both cases the space--time
superprocess has branching mechanism $\psi(v)=\kappa v+v^2/2$.
\begin{enumerate}[label=(\roman*)]
  \item \textit{Finite-mean delays.}
  The motion from $(t,x)$ is $(t+m_1r,x+B_r^\Sigma)$, and the
  initial measure is
  \[
    \mu_F(\dd t,\dd x)=m_1\delta_0(\dd t)\lambda_0(x)\dd x.
  \]
  The occupation measure is $m_1^{-1}X_t^\kappa(\dd x)\dd t$,
  where $X^\kappa$ has motion generator $\cA=\mathcal B/m_1$,
  branching variance $1/m_1$, linear drift $-\kappa/m_1$, and
  initial measure $m_1\lambda_0(x)\dd x$.

  \item \textit{Heavy-tailed delays.}
  The motion from $(t,x)$ is $(t+S_r,x+B_r^\Sigma)$, and the
  initial measure is
  \[
    \mu_H(\dd t,\dd x)=c_Ht^{-\eta}\lambda_0(x)\dd t\dd x.
  \]
  The occupation measure has law $\nu_{\lambda_0}^{\eta,\kappa}$,
  as constructed in Proposition~\ref{prop:ht-lifted}.
\end{enumerate}

Moreover, on the closed space--time $\overline E$ the source measures
generated by the prescribed past have the vague limits
\begin{equation}\label{eq:preparation-clock-limits}
 \frac{\overline H(t/\epsilon^2)}{\delta_\epsilon}
       \lambda_0(x)\dd t\dd x
 \xrightarrow[\epsilon\downarrow0]{v}
 \begin{cases}
  \mu_F,&\delta_\epsilon=\epsilon^2,\quad m_1<\infty,\\
  \mu_H,&\delta_\epsilon=\epsilon^{2\eta},\quad
                   \overline H(s)\sim c_Hs^{-\eta}.
 \end{cases}
\end{equation}
\end{corollary}

\begin{proof}
For (i), the lifted process is concentrated on the
slice $\{t=m_1r\}$ at genealogical time $r$.  Its spatial marginal
at $r=t/m_1$, denoted by $X_t^\kappa$, satisfies
\eqref{eq:nc-SBM-martingale}--\eqref{eq:nc-SBM-bracket}.
Indeed, if $\widetilde X_r$ is the spatial marginal in genealogical
time and $\widetilde M_r(\psi)$ its martingale, then
\begin{align*}
 \ip{X_t^\kappa}\psi-\ip{X_0^\kappa}\psi
 &=\int_0^{t/m_1}\ip{\widetilde X_r}{\mathcal B\psi-\kappa\psi}\dd r
       +\widetilde M_{t/m_1}(\psi)\\
 &=\int_0^t\ip{X_s^\kappa}{\cA\psi-\kappa\psi/m_1}\dd s
       +\widetilde M_{t/m_1}(\psi),\\
 \left\langle\widetilde M(\psi)\right\rangle_{t/m_1}
 &=\int_0^{t/m_1}\ip{\widetilde X_r}{\psi^2}\dd r
   =\frac1{m_1}\int_0^t\ip{X_s^\kappa}{\psi^2}\dd s.
\end{align*}
Here $\mathcal B=\frac12\Sigma:D_x^2$ and
$\cA=\mathcal B/m_1$ by \eqref{eq:diffusion-coefficient}.
  Consequently,
\[
 \int_0^\infty\ip{\mathsf Z_r}{g}\dd r
     =\frac1{m_1}\int_0^\infty
                       \ip{X_t^\kappa}{g(t,\cdot)}\dd t.
\]
Its stationary occupation equation is
\[
m_1\partial_tu+\mathcal B u-\kappa u-\frac{u^2}{2}+g=0.
\]
The initial pairing is
\[
\ip{\mu_F}{u}=m_1\int_{\Rd}\lambda_0(x)u(0,x)\dd x,
\]
which agrees with \eqref{eq:nc-canonical-exponent}.
For bounded profiles of infinite mass, use increasing finite spatial
restrictions and the local first-moment bound, as in
Proposition~\ref{prop:ht-lifted}.  Part (ii) is that proposition.

For the first limit in \eqref{eq:preparation-clock-limits}, substitute
$t=\epsilon^2s$ and use
$\int_0^\infty\overline H(s)\dd s=m_1$ and the dominated convergence theorem.
For the second, the densities converge at every $t>0$ and are bounded
by $Ct^{-\eta}$ on compact time intervals by
\eqref{eq:ht-tail-bounds}.  The dominated convergence theorem applies because
$\eta<1$.  Spatial integration is over the compact support of the
test function in both cases.
\end{proof}

\begin{remark}[A common inverse clock does not preserve branching]
The mean representation \eqref{eq:ht-mean-profile} does not extend to
the whole law by applying one independent inverse stable clock to an
ordinary Dawson--Watanabe process.  A direct distinction is the
branching property in the initial profile:
\begin{equation}\label{eq:ht-profile-branching}
 \nu_{\lambda_1+\lambda_2}^{\eta,\kappa}
   \ \stackrel{\mathrm{law}}{=}\quad
 \nu_{\lambda_1}^{\eta,\kappa}
       +\widetilde\nu_{\lambda_2}^{\eta,\kappa},
\end{equation}
where the terms on the right are independent.  This follows from
linearity of \eqref{eq:ht-exponent} in $\lambda_0$.

Let $\lambda_0=p_1^\Sigma$, with $p_r^\Sigma$ defined in
\eqref{eq:ht-brownian-kernel}, and take an ordinary superprocess $Y$ with
spatial generator $\mathcal B$, branching variance one, linear drift
$-\kappa$, and initial measure $a\lambda_0(x)\dd x$.
Let one inverse clock $E$ be independent of $Y$, and set
$V_a(\dd t,\dd x)=Y_{E_t}(\dd x)\dd t$.
Conditional on $E$, its Laplace functional is
$\exp\{-aH_E(g)\}$ by the branching property of $Y$.
Choose $g(t,x)=q_0(t)p_1^\Sigma(x)$, with nonzero
$q_0\in C_c^+((0,\infty))$.  Its conditional mean is
\begin{equation}\label{eq:ht-common-clock-mean}
 \E[\ip{V_a}{g}\mid E]
    =a\int_0^\infty q_0(t)e^{-\kappa E_t}
                                  p_{2+E_t}^\Sigma(0)\dd t.
\end{equation}
The function $r\mapsto e^{-\kappa r}p_{2+r}^\Sigma(0)$ is strictly
decreasing for every $\kappa\geq0$.  If the support of $q_0$ lies
in $[t_1,t_2]\subset(0,\infty)$, the events
$E_{t_2}\leq r_1$ and $E_{t_1}\geq r_2$, with $0<r_1<r_2$,
both have positive probability: the stable subordinator at any
positive time has full support on $(0,\infty)$.
They give disjoint ranges for the integral in
\eqref{eq:ht-common-clock-mean}, so that integral is nonconstant.
These Gaussian spatial test functions are admissible by monotone compact
cutoffs and the finite initial mass.

Thus $H_E(\theta g)$ is nonconstant for some $\theta>0$;
otherwise its derivative at zero, the conditional mean per unit $a$,
would be deterministic.  The strict Cauchy--Schwarz inequality gives
\[
 \E\left[e^{-2H_E(\theta g)}\right]>
                   \left(\E\left[ e^{-H_E(\theta g)}\right]\right)^2,
\]
contradicting \eqref{eq:ht-profile-branching}.
A common random clock introduces dependence between clans, whereas
the space--time superprocess preserves their branching property.
\end{remark}

\subsubsection{Relation to temporal nearly unstable limits}

The joint resolvent symbol
$\left(c_\eta z^\eta+\kappa+\frac12\xi^{\top}\Sigma\xi\right)^{-1}$
separates the stable temporal clock from Brownian spatial motion.
At spatial frequency $\xi=0$ its time kernel is
\begin{equation}\label{eq:ht-Mittag-Leffler}
  k_{\eta,\kappa}(t)
       =\frac{t^{\eta-1}}{c_\eta}
          E_{\eta,\eta}\!\left(-\frac{\kappa}{c_\eta}t^\eta\right),
  \qquad
  E_{\eta,\beta}(z):=\sum_{n=0}^\infty
                              \frac{z^n}{\Gamma(\eta n+\beta)}.
\end{equation}
This is the Mittag--Leffler resolvent behind
\cite{jaisson2016rough}, with the constants expressed in our notation.
For a homogeneous profile $\lambda_0\equiv\lambda$,
\eqref{eq:ht-mean-profile} becomes
\[
  \rho^{\eta,\kappa}(t,x)
       :=\lambda E_{\eta,1}\!\left(-\frac{\kappa}{c_\eta}t^\eta\right).
\]
\begin{remark}[Comparison with temporal nearly unstable scaling]
\label{rem:formal-temporal-reduction}
Setting $\xi=0$ removes the spatial part of the resolvent symbol.
With $T=\epsilon^{-2}$ and $a_T=a_{T^{-1/2}}$,
\eqref{eq:ht-constants} reads $T^\eta(1-a_T)\to\kappa$ as
$T\to\infty$, the
temporal nearly unstable rate in \cite{jaisson2016rough} when $\kappa>0$.
The sources of activity differ: the prescribed past here provides
$c_Ht^{-\eta}\lambda_0$, while that temporal model uses positive-time
immigration.  Its temporal density result requires $\eta>1/2$;
our stochastic event-measure convergence \eqref{eq:ht-stochastic}
holds for $0<\eta<1$.
\end{remark}

\section{Conclusion}
\label{sec:conclusion}

We established scaling limits for the event measures of critical and
nearly critical space--time Hawkes processes initialized by a prescribed
Poisson field of past events.  The past intensity produces three
regimes.  At high density, we obtained a hydrodynamic limit which is a deterministic diffusion
profile.  At the threshold density, branching fluctuations survive and
the limit is the occupation measure of a Dawson--Watanabe process.  At
low density, the event measure vanishes locally, while its first-order
Laplace asymptotics retain the contribution of rare large clusters
through the canonical measure.

The temporal tail determines how the prescribed past appears in the
limit.  For finite-mean delays, the past collapses to the initial
measure of the limiting superprocess.  For heavy-tailed delays, it
persists at positive times as a memory profile, and the stochastic limit
is represented as the occupation measure of a quadratic superprocess on
space--time.  This representation applies in both the critical and
nearly critical cases.  It also shows that a common inverse stable clock
reproduces the mean behavior but does not preserve the branching law.

We also derived the phase diagram by matching the scaling of the Poisson
past to the observation scales and identified a one-dimensional local recovery phenomenon. A
shorter observation window can retain stochastic fluctuations that
disappear at the original macroscopic scale; depending on the scaling of
the prescribed past, the limit either keeps the full initial profile or freezes it at
the observation point.  Together, these results give a unified
description of how density, temporal memory, and observation scale
determine the macroscopic behavior of the Hawkes event measure.


\subsection*{Acknowledgements}

Lingjiong Zhu would like to thank Prof. S.R.S. Varadhan
for helpful discussions and suggestions.

\bibliographystyle{alpha}
\bibliography{bibtex}

\end{document}